\documentclass[12pt,a4paper]{article}

\usepackage{amsbsy,amsfonts,amsmath,amssymb,amsthm}
\usepackage{array}
\usepackage{bm}
\usepackage{color}
\usepackage{enumerate}
\usepackage{mathrsfs}
\usepackage{latexsym}
\usepackage{comment}
\usepackage[colorlinks=true]{hyperref}
\hypersetup{urlcolor=blue, citecolor=blue, linkcolor=blue}

\definecolor{wineRed}{rgb}{0.7,0,0.3}
\definecolor{grandBleu}{rgb}{0,0,0.8}
\definecolor{darkGreen}{rgb}{0,0.4,0}
\definecolor{blueViolet}{rgb}{0.4,0,1.0}
\definecolor{bloodOrange}{rgb}{0.85,0.05,0}
\definecolor{mycolor}{rgb}{0.8,0,0.2}
\definecolor{}{rgb}{0.8,0,0.2}

\usepackage[textsize=small]{todonotes}
\usepackage{cite}

\usepackage[colorlinks=true]{hyperref}
\hypersetup{urlcolor=blue, citecolor=red}

\usepackage{comment}

\DeclareMathAlphabet{\mathpzc}{OT1}{pzc}{m}{it}

\usepackage[textsize=small]{todonotes}
\numberwithin{equation}{section}

\theoremstyle{plain}
\newtheorem{theorem}{Theorem}[section]
\newtheorem{corollary}{Corollary}
\newtheorem{main}{Main Theorem}
\newtheorem{lemma}[theorem]{Lemma}
\newtheorem{proposition}{Proposition}

\theoremstyle{definition}
\newtheorem{definition}[theorem]{Definition}
\newtheorem{remark}{Remark}

\newtheorem{keyLem}{Key-Lemma}

\def\Sgn{\mathop{\mathrm{Sgn}}\nolimits}

\begin{document}
\thispagestyle{plain}
\begin{center}
    \textbf{\Large OPTIMAL CONTROL PROBLEMS \\ FOR 1D PARABOLIC STATE- \\ SYSTEMS OF {K}{W}{C} TYPES  WITH \\[1ex] DYNAMIC BOUNDARY CONDITIONS}\footnotemark[1]
\end{center}
    \bigskip
\begin{center}
    \textsc{Shodai Kubota}
    \\[1ex]
    {Department of Mathematics and Informatics, \\ Graduate School of Science and Engineering, Chiba University, \\ 1--33, Yayoi-cho, Inage-ku, 263--8522, Chiba, Japan}
    \\[0ex]
    ({\ttfamily skubota@chiba-u.jp})
\end{center}
\begin{center}
    \textsc{Ryota Nakayashiki} 
    \\[1ex]
    {Department of General Education, Salesian Polytechnic, \\ 4--6--8, Oyamagaoka, Machida-city, Tokyo, 194--0215, Japan}
    \\[0ex]
    ({\ttfamily nakayashiki.ryota@salesio-sp.ac.jp})
\end{center}
\begin{center}
    \textsc{and}
\end{center}
\begin{center}
    \textsc{Ken Shirakawa}
    \\[1ex]
    {Department of Mathematics, Faculty of Education, Chiba University, \\ 1--33, Yayoi-cho, Inage-ku, 263--8522, Chiba, Japan}
    \\[0ex]
    ({\ttfamily sirakawa@faculty.chiba-u.jp})
\end{center}

\footnotetext[1]{
The work of the third author is supported by Grant-in-Aid for Scientific Research (C) No. 16K05224 and No. 20K03672, JSPS.   
\\[1ex]
AMS Subject Classification: 
                    35K59, 
                    35K61, 
                    49J20, 
                    49K20, 
                    74N05, 
                    74N20. 
\\[1ex]
Keywords: optimal control problem, 1D parabolic state-system of KWC type, dynamic boundary condition, grain boundary motion, physically realistic problem, regularized approximating problems, the first order necessary optimality conditions; limiting optimality condition
}
\bigskip

\bigskip

\bigskip

\noindent
{\bf Abstract.} In this paper, we consider a class of optimal control problems governed by 1D parabolic state-systems of KWC types with dynamic boundary conditions. The state-systems are based on a phase-field model of grain boundary motion, proposed in [Kobayashi--Warren--Carter, Physica D, 140, 141--150, 2000], and in the context, the dynamic boundary conditions are supposed to reproduce the transmitted heat exchanges between interior and boundary of a  polycrystal body. Our optimal control problems are labeled by using a constant $ \varepsilon \geq 0 $, and roughly summarized, the case when $ \varepsilon = 0 $ and the cases when $ \varepsilon > 0 $ correspond to the physically realistic setting, and its regularized approximating ones, respectively. Under suitable assumptions, the mathematical results concerned with: the solvability and continuous dependence for the state-systems; the solvability and $ \varepsilon $-dependence of optimal control problems; and the first order necessary optimality conditions in the problems when $ \varepsilon > 0 $ and the limiting optimality condition as $ \varepsilon \downarrow 0 $; will be obtained in forms of three Main Theorems of this paper. 

\newpage

\section*{Introduction}
Let $ (0, T) $ be a time-interval with a constant $ 0 < T < \infty $, let $\Omega := (0, 1) \subset \mathbb{R} $ be a one-dimensional spatial domain, and let $ \Gamma := \{0, 1\} $ be the boundary of $ \Omega $. Besides, we set $ Q := (0, T) \times \Omega $ and $ \Sigma := (0, T) \times \Gamma $, and  we define:
\begin{equation*}
    \left\{ \parbox{12.5cm}{
        $ H := L^2(\Omega) $, ~ $ H_\Gamma := \left\{ \begin{array}{l|l}
            \tilde{w} & \tilde{w} : \Gamma \longrightarrow \mathbb{R}
        \end{array} \right\} ~ (\sim \mathbb{R}^2) $, ~ $ \mathbb{X} := H \times H_\Gamma $,
        \\[1ex]
        $ \mathscr{H} := L^2(0, T; H) $, ~ $ \mathscr{H}_\Gamma := L^2(0, T; H_\Gamma) $, ~ and ~ $ \mathfrak{X} := \mathscr{H} \times \mathscr{H}_\Gamma $,
    } \right.
\end{equation*}
as the base spaces for our problems. 
\medskip

In this paper, we consider a class of optimal control problems governed by the following 1D state-systems, which are denoted by (S)$_\varepsilon $, with $ \varepsilon \geq 0 $.

\begin{description}
    \item[\textmd{(S)$ _\varepsilon $}]Find a triplet of functions $ [\bm{\eta}, \theta] = [\eta, \eta_\Gamma, \theta] \in \mathfrak{X} \times \mathscr{H} $ with $ \bm{\eta} = [\eta, \eta_\Gamma]$, fulfilling:
\end{description}

\begin{subequations}\label{0.1}
    \hspace{4ex}
    $
    \left\{  \hspace{-2ex}\parbox{14.5cm}{
    \vspace{-2ex}
    \begin{align}
        &\partial_{t} \eta -\partial_{x}^{2} \eta +g(\eta) +\alpha'(\eta) \sqrt{\varepsilon^2 +|\partial_{x} \theta|^2} = L u \mbox{  in $ Q $,}~\hspace{18ex}~
        \label{0.1a}
        \\[1ex]
        &\partial_t \eta_\Gamma(t, \ell) +(-1)^{\ell -1} \partial_{x} \eta_{_{|_\Gamma}}(t, \ell) = L_\Gamma u_\Gamma(t, \ell), ~~ (t, \ell) \in \Sigma,
        \label{0.1b}
        \\[1ex]
        & \eta_{_{|_\Gamma}} = \eta_\Gamma \mbox{ on $ \Sigma $}
        \label{0.1c}
        \\[1ex]
        &\eta(0, x) = \eta_{0}(x), ~ x \in \Omega, ~~ \eta_\Gamma(0, \ell) = \eta_{\Gamma, 0}(\ell), ~ \ell \in \Gamma;
        \label{0.1d}
    \end{align}
    \vspace{-3ex}
}\right. 
$
\end{subequations}
\bigskip

\begin{subequations}\label{0.2}
    \hspace{4ex}
    $ \left\{ \hspace{-2ex} \parbox{14.5cm}{
        \vspace{-2ex}
        \begin{align}
            & \alpha_{0}(t, x) \partial_{t} \theta -\partial_{x} \left( \alpha(\eta) \frac{\partial_{x}\theta}{\sqrt{\varepsilon^2 +|\partial_{x} \theta|^2}} +\nu^{2}\partial_{x}\theta\right) = M v \mbox{ in $Q$,}~\hspace{8ex}~
            \label{0.2a}
            \\
            & \theta(t, \ell) = 0, ~~ (t, \ell) \in \Sigma , 
            \label{0.2c}
            \\[1ex]
            & \theta(0, x) = \theta_{0}(x), ~ x \in \Omega.
            \label{0.2d}
        \end{align}
    \vspace{-3ex}
}\right. 
$
\end{subequations}
\bigskip

\noindent
For each $ \varepsilon \geq 0 $, the optimal control problem is denoted by (OP)$_\varepsilon$, and is prescribed as follows.
\begin{description}
    \item[\textmd{(OP)$_\varepsilon$}]Find a triplet $ [\bm{u}^*, v^*] = [u^*, u_\Gamma^*, v^*] \in \mathfrak{X} \times \mathscr{H} $ with $ \bm{u}^* = [u^*, u_\Gamma^*] $, called \emph{optimal control}, which minimizes a cost functional $ \mathcal{J}_\varepsilon $, defined as:
        \begin{align}\label{J}
            \mathcal{J}_\varepsilon : [\bm{u}, &\,  v] = [u, u_\Gamma, v] \in \mathfrak{X} \times \mathscr{H} \mapsto \mathcal{J}_\varepsilon(\bm{u}, v) = \mathcal{J}_\varepsilon(u, u_\Gamma, v) 
            \nonumber
            \\
            := &~ \frac{K}{2} \int_0^T |(\eta -\eta_\mathrm{ad})(t)|_H^2 \, dt  + \frac{K_\Gamma}{2} \int_0^T |(\eta_\Gamma -\eta_{\Gamma, \mathrm{ad}})(t)|_{H_\Gamma}^2 \, dt
            \nonumber
            \\
            & ~+\frac{\Lambda}{2} \int_0^T |(\theta -\theta_\mathrm{ad})(t)|_H^2 \, dt  
            \nonumber 
            \\
            &~  +\frac{L}{2} \int_0^T |u(t)|_H^2 \, dt +\frac{L_\Gamma}{2} \int_0^T |u_\Gamma(t)|_{H_\Gamma}^2 \, dt
            \nonumber
            \\
            &~ +\frac{M}{2} \int_0^T |v(t)|_H^2 \, dt \in [0, \infty),  
        \end{align}
        where $[\bm{\eta}, \theta] = [\eta, \eta_\Gamma, \theta] \in \mathfrak{X} \times \mathscr{H} $ with $\bm{\eta} = [\eta, \eta_\Gamma] $ is a triplet solving the state-system (S)$_\varepsilon$. 
\end{description}
The state-system (S)$_\varepsilon$ is a type of \emph{KWC system}, i.e. it is based on a phase-field model of grain boundary motion, proposed by Kobayashi--Warren--Carter \cite{MR1752970, MR1794359}. In the original KWC system, the components $ \eta \in \mathscr{H} $ and $ \theta \in \mathscr{H} $ are order parameters which indicate the \emph{orientation order} and \emph{orientation angle} of a polycrystal body $ \Omega $, respectively. Also, the component $ \eta_\Gamma \in \mathscr{H}_\Gamma $ is an order parameter, which influences the dynamics of $\eta$, as an external factor exchanged via the boundary $ \Gamma $ of polycrystal.  $[\eta_0, \eta_{\Gamma, 0}] \in \mathbb{X}$ and $\theta_0 \in H$ are initial data of $ \bm{\eta} = [\eta, \eta_\Gamma] $ and $\theta$, respectively, and for simplicity, these initial data are written in a form of a \emph{initial triplet} $ [\bm{\eta}_0, \theta_0] = [\eta_0, \eta_{\Gamma, 0}, \theta_0] \in \mathbb{X} \times H $ with the \emph{initial pair} $ \bm{\eta}_0 = [\eta_0, \eta_{\Gamma, 0}] $. The forcing triplet $ [\bm{u}, v] \in \mathfrak{X} \times \mathscr{H} $ with the forcing pair $\bm{u} = [u, u_\Gamma]$, denotes the control variables that can control the profile of solution $ [\bm{\eta}, \theta] \in \mathfrak{X} \times \mathscr{H} $ to (S)$_\varepsilon$. $ 0 < \alpha_0 \in W^{1, \infty}(Q) $ and $ 0 < \alpha \in C^2(\mathbb{R}) $ are given functions to reproduce the mobilities of grain boundary motions. Besides, ``$|_\Gamma$'' denotes the trace on $\Gamma$ for a Sobolev function. Finally, $ g \in W_\mathrm{loc}^{1, \infty}(\mathbb{R}) $ is a perturbation for the orientation order $ \eta $, and $ \nu > 0 $ is a fixed constant to relax the diffusion of the orientation angle $ \theta $. 
\medskip

As a mathematical model of grain boundary motion, (S)$_\varepsilon$ can be said as a coupled system of an Allen--Cahn type equation \eqref{0.1a} subject to the dynamic boundary condition  \{\eqref{0.1b},\eqref{0.1c}\}, and a quasilinear diffusion equation \eqref{0.2a} subject to the homogeneous Dirichlet boundary condition \eqref{0.2c}. However, it should be noted that the PDE \eqref{0.1b} can be regarded as a governing equation on $ \Gamma $, and the initial-boundary value problem \eqref{0.1} also forms a kind of \emph{transmission problem} between the PDE \eqref{0.1a} in $\Omega$ and the PDE \eqref{0.1b} on $ \Gamma $, subject to the \emph{transmission condition} \eqref{0.1c}.  
\bigskip

Since $ \bm{u} = [u, u_\Gamma] \in \mathfrak{X} $ is a forcing term of the Allen-Cahn type equation \eqref{0.1a}, the components $ u \in \mathscr{H} $ and $ u_\Gamma \in \mathscr{H}_\Gamma $ can be understood as the temperature controls on the interior $ \Omega $ and the boundary $ \Gamma $ of polycrystal, respectively. Meanwhile, the quasilinear diffusion equation \eqref{0.2a} is to reproduce crystalline micro-structure of polycrystal, and the case when $ \varepsilon = 0 $ is the closest to the original setting adopted by Kobayashi--Warren--Carter \cite{MR1752970, MR1794359}. Indeed, when $ \varepsilon = 0 $, the quasi-linear diffusion as in \eqref{0.2a} is described in a singular form $ -\partial_x \bigl( \alpha(\eta) \frac{\partial_x \theta}{|\partial_x \theta|} +\nu^2 \partial_x \theta \bigr) $, and this type of singularity is said to be effective to reproduce the \emph{facet}, i.e. the locally uniform (constant) phase in each oriented grain (cf. \cite{MR1752970, MR1794359,MR2746654,MR2101878,MR2436794,MR2033382,MR2223383,MR3038131,MR3670006,MR3268865,MR2836557,MR1865089,MR1712447,MR3951294}). 
Hence, the systems (S)$_\varepsilon$, for positive $ \varepsilon $, can be said as regularized approximating systems, that are to approach to the physically realistic situation, reproduced by the limiting system (S)$_0$, as $ \varepsilon \downarrow 0 $. 
\bigskip

Furthermore, in the optimal control problems (OP)$_\varepsilon$ for $ \varepsilon \geq 0 $, the functions $ \eta_\mathrm{ad} \in \mathscr{H} $, $ \eta_{\Gamma, \mathrm{ad}} \in \mathscr{H}_\Gamma$, and $ \theta_\mathrm{ad} \in \mathscr{H} $ are given \emph{admissible target profile} of the order parameters $ \eta $, $ \eta_\Gamma $, and $ \theta $, respectively, and the coefficients $ K \geq 0 $, $ K_\Gamma \geq 0 $, $ \Lambda \geq 0 $, $ L \geq 0 $, $ L_\Gamma \geq 0 $, and $ M \geq 0 $ are fixed constants which are to adjust the meaning of optimality in the problems (OP)$ _\varepsilon $. 
\bigskip

This paper focuses on two issues:
\begin{itemize}
    \item[\textmd{~~$\sharp \, 1)$}]key-properties of the state-systems (S)$ _\varepsilon $, for $ \varepsilon \geq 0 $; 
        \vspace{0ex}
    \item[\textmd{~~$\sharp \, 2)$}]mathematical analysis of the optimal control problems (OP)$ _\varepsilon $, for $ \varepsilon \geq 0 $.
\end{itemize}
With regard to the first issue $\sharp \, 1)$, some kindred KWC type systems have been studied by several mathematicians, e.g. \cite{MR2469586,Equadiff_Nakayashiki17,MR3888636}, and in particular, the analytic ideas, as in \cite[Main Theorems 1 and 2]{MR3888636}, would be effective for the well-posedness and $ \varepsilon $-dependence of the system (S)$ _\varepsilon $. However, the previous works \cite{MR2469586,Equadiff_Nakayashiki17,MR3888636} adopted the homogeneous setting of forcing, and imposed different types of boundary conditions with this study. In this light, we need to enhance the existing mathematical method before we deal with the study of our optimal control problems (OP)$_\varepsilon$. Meanwhile, for issue $ \sharp \, 2) $, there are now a number of previous works \cite{MR3060198,MR3643644,MR1161486,MR3300403,MR2578571}, which dealt with optimal control problems, governed by PDE systems kindred to (S)$_\varepsilon$. Hence, by integrating the previous works \cite{MR2469586,Equadiff_Nakayashiki17,MR3888636,MR3060198,MR3643644,MR1161486,MR3300403,MR2578571}, we can expect to develop a certain mathematical control theory that enables to handle dynamically transmitted situations, as in the dynamic boundary condition of our state-system (S)$_\varepsilon$. 
\bigskip

Now, based on these, we set the goal of this paper to prove three Main Theorems, summarized as follows:
\begin{description}
    \item[{\boldmath Main Theorem 1.}]Mathematical results concerning the following items:
\item[{\boldmath~~(I-A)(Solvability of state-systems)}]Existence and uniqueness for the state-system (S)$_\varepsilon$, for any $ \varepsilon \geq 0 $.
\item[{\boldmath~~(I-B)(Continuous dependence among state-systems)}]Continuous dependence of solutions to the systems (S)$_\varepsilon$, with respect to the initial triplet, the forcing triplet, and the constant $\varepsilon \geq 0$. 
\end{description}
\begin{description}
\item[{\boldmath Main Theorem 2.}]Mathematical results concerning the following items:
\item[{\boldmath~~(II-A)(Solvability of optimal control problems)}]Existence for the optimal control problem (OP)$_\varepsilon$, for any $ \varepsilon \geq 0 $.
\item[{\boldmath~~(II-B)($\varepsilon$-dependence of optimal controls)}]Some semi-continuous association between the optimal controls, with respect to the initial triplet and the constant $ \varepsilon \geq 0$.    
\end{description}
\begin{description}
\item[{\boldmath Main Theorem 3.}]Mathematical results concerning the following items:
        \item[{\boldmath~~(III-A)(Necessary optimality conditions in cases of $ \varepsilon > 0 $)}] Derivation of the first order necessary optimality conditions for (OP)$_\varepsilon$ via adjoint method.
\item[{\boldmath~~(III-B)(Limiting optimality conditions as ~$\varepsilon \downarrow 0$)}]The limiting adjoint system as $ \varepsilon \downarrow 0 $ associated with the optimality condition for the problem (OP)$_0$.   
\end{description}

The Main Theorems are stated in Section 3, after the preliminaries in Section 1, and the auxiliary results in Section 2. The Main Theorems are proved in Sections 4--6, and in particular, the Main Theorem 3 is proved by means of three Theorems \ref{Prop01}--\ref{Prop03} which are stated as a part of the auxiliary results in Section 2. The proofs of Theorems \ref{Prop01}--\ref{Prop03} are given in the last Section 7 of appendix.

\section{Preliminaries} 

We begin by prescribing the notations used throughout this paper. 
\medskip

\noindent
\underline{\textbf{\textit{Abstract notations.}}}
For an abstract Banach space $ X $, we denote by $ |\cdot|_{X} $ the norm of $ X $, and denote by $ \langle \cdot, \cdot \rangle_X $ the duality pairing between $ X $ and its dual $ X^* $. In particular, when $ X $ is a Hilbert space, we denote by $ (\cdot,\cdot)_{X} $ the inner product of $ X $. 

For any subset $ A $ of a Banach space $ X $, let $ \chi_A : X \longrightarrow \{0, 1\} $ be the characteristic function of $ A $, i.e.:
    \begin{equation*}
        \chi_A: w \in X \mapsto \chi_A(w) := \begin{cases}
            1, \mbox{ if $ w \in A $,}
            \\[0.5ex]
            0, \mbox{ otherwise.}
        \end{cases}
    \end{equation*}

For two Banach spaces $ X $ and $ Y $,  we denote by $  \mathscr{L}(X; Y)$ the Banach space of bounded linear operators from $ X $ into $ Y $, and in particular, we let $ \mathscr{L}(X) := \mathscr{L}(X; X) $. 

For Banach spaces $ X_1, \dots, X_N $, with $ 1 < N \in \mathbb{N} $, let $ X_1 \times \dots \times X_N $ be the product Banach space endowed with the norm $ |\cdot|_{X_1 \times \cdots \times X_N} := |\cdot|_{X_1} + \cdots +|\cdot|_{X_N} $. However, when all $ X_1, \dots, X_N $ are Hilbert spaces, $ X_1 \times \dots \times X_N $ denotes the product Hilbert space endowed with the inner product $ (\cdot, \cdot)_{X_1 \times \cdots \times X_N} := (\cdot, \cdot)_{X_1} + \cdots +(\cdot, \cdot)_{X_N} $ and the norm $ |\cdot|_{X_1 \times \cdots \times X_N} := \bigl( |\cdot|_{X_1}^2 + \cdots +|\cdot|_{X_N}^2 \bigr)^{\frac{1}{2}} $. In particular, when all $ X_1, \dots,  X_N $ coincide with a Banach space $ Y $, we write:
\begin{equation*}
    [Y]^N := \overbrace{Y \times \cdots \times Y}^{\mbox{$N$ times}}.
\end{equation*}
Additionally, for any (possibly nonlinear) transform $ \mathcal{T} : X \longrightarrow Y $, we let:
\begin{equation*}
    \mathcal{T}[w_1, \dots, w_N] := \bigl[ \mathcal{T} w_1, \dots, \mathcal{T} w_N \bigl] \mbox{ in $ [Y]^N $, \quad for any $ [w_1, \dots, w_N] \in [X]^N $.}
\end{equation*}

\noindent
\underline{\textbf{\textit{Specific notations of this paper.}}} 
As is mentioned in the introduction, let $ (0, T) \subset \mathbb{R}$ be a bounded time-interval with a finite constant $ T > 0 $, and let  $ \Omega := (0, 1) \subset \mathbb{R} $ be a one-dimensional bounded spatial domain. We denote by $ \Gamma $ the boundary $ \partial \Omega = \{0, 1\} $ of $ \Omega $, and we define 
\begin{equation*}
    n_\Gamma(\ell) := (-1)^{\ell -1} \mbox{ for any $\ell \in \Gamma = \{0, 1\}$.}
\end{equation*}
Besides, we let $ Q := (0, T) \times \Omega $ and $ \Sigma := (0, T) \times \Gamma $. 

Throughout this paper, we denote by $ \partial_t $ and $ \partial_x $ the distributional time-derivative and the distributional spatial-derivative, respectively. Also, the measure theoretical phrases, such as ``a.e.'', ``$dt$'', ``$dx$'', and so on, are all with respect to the Lebesgue measure in each corresponding dimension. Additionally, ``$|_\Gamma$" denotes the trace on $\Gamma$ for a Sobolev function.

On this basis, we define  
\begin{equation*} 
    \left\{ \parbox{10.5cm}{
        $ H := L^2(\Omega) $, ~ $ H_\Gamma := \left\{ \begin{array}{l|l}
            \tilde{w} & \tilde{w} : \Gamma \longrightarrow \mathbb{R}
        \end{array} \right\} ~ (\sim \mathbb{R}^2) $, 
        \\[0.5ex]
        $ V := H^1(\Omega) $, ~ $V_0 := H_0^1(\Omega) $, 
    } \right.
\end{equation*}
\begin{equation*}
    \left\{ \parbox{10.5cm}{
        $ \mathbb{X} := H \times H_\Gamma $, ~ $ \mathbb{V} = V \times H_\Gamma $,
        \\[0.5ex]
        $ \mathbb{W} := \left\{ \begin{array}{l|l}
            [\tilde{w}, \tilde{w}_\Gamma] \in \mathbb{V} & \tilde{w}_{_{|_\Gamma}}(\ell) = \tilde{w}_\Gamma(\ell), ~ \ell \in \Gamma
        \end{array} \right\} ,$
    } \right. 
\end{equation*}
\begin{equation*}
    \left\{ \parbox{10.5cm}{
        $ \mathscr{H} := L^2(0, T; H) $, ~ $ \mathscr{H}_\Gamma := L^2(0, T; H_\Gamma) $,
        \\[1ex]
        $ \mathscr{V} := L^2(0, T; V) $, ~ $\mathscr{V}_0 := L^2(0, T; V_0)$,
    } \right.
\end{equation*}
and
\begin{equation*}
    \left. \parbox{10.5cm}{
        $ \mathfrak{X} := \mathscr{H} \times \mathscr{H}_\Gamma ~(= L^2(0, T; \mathbb{X})) $, ~ and ~ $\mathfrak{W} := L^2(0, T; \mathbb{W}) $.
    } \right.
\end{equation*}
Note that $ \mathbb{W} $ is a closed linear subspace in the Hilbert space $ \mathbb{V} $, so that $ \mathbb{W} $ is also a Hilbert space endowed with the inner product of $ \mathbb{V} $. 

In this paper, we identify the Hilbert spaces $ H $ and $ \mathscr{H} $ with their dual spaces. On this basis, we have the following relationships of continuous embeddings:
\begin{equation*}
    \left\{ \parbox{11cm}{$ V \subset H = H^* \subset V^* $, ~$ \mathscr{V} \subset \mathscr{H} = \mathscr{H}^* \subset \mathscr{V}^* $, 
    \\[0.5ex]
    $ \mathbb{W} \subset \mathbb{V} \subset \mathbb{X} = \mathbb{X}^* \subset \mathbb{V}^* \subset \mathbb{W}^* $, ~and~ $ \mathfrak{W} \subset \mathfrak{X} = \mathfrak{X}^* \subset \mathfrak{W}^* $,} \right.
\end{equation*}
among the Hilbert spaces $ H $, $ V $, $ \mathscr{H} $, $ \mathscr{V} $, $ \mathbb{X} $, $ \mathbb{V} $, $ \mathbb{W} $, $ \mathfrak{X} $, and $ \mathfrak{W} $, and the respective dual spaces $ H^* $, $ V^* $, $ \mathscr{H}^* $, $ \mathscr{V}^* $, $ \mathbb{X}^* $, $ \mathbb{V}^* $, $ \mathbb{W}^* $, $ \mathfrak{X}^* $, and $ \mathfrak{W}^* $. 

\begin{remark}\label{Rem.Prelim01}
    Due to the one-dimensional embeddings $ V \subset C(\overline{\Omega}) $ and $ V_0 \subset C(\overline{\Omega}) $, it is easily checked that: \begin{equation}\label{emb01}
        \begin{cases}
            \hspace{-2ex}
            \parbox{9cm}{
                \vspace{-2ex}
                \begin{itemize}
                    \item if $ \check{\mu} \in H $ and $ \check{p} \in V $, then $ \check{\mu} \check{p} \in H $, and $ |\check{\mu} \check{p}|_H \leq \sqrt{2} |\check{\mu}|_H |\check{p}|_V $,
                    \item if $ \hat{\mu} \in L^\infty(0, T; H) $ and $ \hat{p} \in \mathscr{V} $, then $ \hat{\mu} \hat{p} \in \mathscr{H} $, and $ |\hat{\mu} \hat{p}|_\mathscr{H} \leq \sqrt{2} |\hat{\mu}|_{L^\infty(0, T; H)} |\hat{p}|_\mathscr{V} $.
                \vspace{-2ex}
                \end{itemize}
            }
        \end{cases}
    \end{equation}
    Here, we note that the constant $ \sqrt{2} $ corresponds to the constant of embedding $ V \subset C(\overline{\Omega}) $. Moreover, under the setting $ \Omega := (0, 1) $, this $ \sqrt{2} $ can be used as a upper bound of the constants of embeddings $ V \subset L^q(\Omega) $ and $ V_0 \subset L^q(\Omega) $, for all $ 1 \leq q \leq \infty $.
\end{remark}
\begin{remark}\label{Rem.Prelim02}
Let us take any $ \tilde{a} \in W^{1, \infty}(Q) \cup L^\infty(0, T; W^{1, \infty}(\Omega)) $ and any $ w \in \mathscr{V}_0^* $. Then, we can say that $ \tilde{a} w $ $ (= w \tilde{a}) $ $\in \mathscr{V}_0^* $, via the following variational form:
    \begin{equation*}
        \langle \tilde{a} w, \psi \rangle_{\mathscr{V}_0} := \langle w, \tilde{a} \psi \rangle_{\mathscr{V}_0}, \mbox{ for any $ \psi \in \mathscr{V}_0 $,}
\end{equation*}
and can estimate that:
\begin{equation*}
    |\tilde{a} w|_{\mathscr{V}_0^*} \leq (1 +\sqrt{2}) \bigl( |\tilde{a}|_{L^\infty(Q)} +|\partial_x \tilde{a}|_{L^\infty(Q)} \bigr) |w|_{\mathscr{V}_0^*},
\end{equation*}
        by using the constant $ \sqrt{2} $ of the embedding $ V_0 \subset H $. 
        Also, if $ \{\tilde{a}_n\}_{n = 1}^\infty \subset W^{1, \infty}(Q) \cup L^\infty(0, T; W^{1, \infty}(\Omega)) $, $ \{w_n\}_{n = 1}^\infty \subset \mathscr{V}_0^* $, and
        \begin{equation*}
            \begin{cases}
                \tilde{a}_n \to \tilde{a} \mbox{ in $ L^\infty(Q) $,}
                \\[1ex]
                \partial_x \tilde{a}_n \to \partial_x \tilde{a} \mbox{ in $ L^\infty(Q) $,}
            \end{cases}
        \end{equation*}
        and
        \begin{equation*}
            w_n \to w \mbox{ weakly in $ \mathscr{V}_0^* $, as $ n \to \infty $,}
        \end{equation*}
        it holds that: 
        \begin{equation*}
            \tilde{a}_n w_n \to \tilde{a} w \mbox{ weakly in $ \mathscr{V}_0^* $, as $ n \to \infty $,}
        \end{equation*}
since
        \begin{equation*}
        \begin{array}{c}
            \mbox{$ \tilde{a} \psi \in \mathscr{V}_0 $, $ \{ \tilde{a}_n w_n \}_{n = 1}^\infty \subset \mathscr{V}_0 $, and } \tilde{a}_n \psi \to \tilde{a} \psi \mbox{ in $ \mathscr{V}_0 $ as $ n \to \infty $,}
            \\[1ex]
            \mbox{ for any $\psi \in \mathscr{V}_0$ .}
        \end{array}
        \end{equation*}
        In particular, if $ \tilde{a} \in W^{1, \infty}(Q) $ and $ w \in W^{1, 2}(0, T; V_0^*) $, then 
\begin{equation*}
    \mbox{$ \tilde{a} w \in W^{1, 2}(0, T; V_0^*) $, and } \partial_t (\tilde{a} w) = \tilde{a} \partial_t w +w \partial_t \tilde{a} \mbox{ in $ \mathscr{V}_0^* $.}
\end{equation*}
        Moreover, if $ \tilde{a} \in W^{1, \infty}(Q) \cup L^\infty(0, T; W^{1, \infty}(\Omega)) $, and $ \log \tilde{a} \in L^\infty(Q) $, then  it is estimated that:
\begin{equation*}
    |\tilde{a} w|_{\mathscr{V}_0^*} \geq \frac{\inf \tilde{a}(Q)^2 }{(1 +\sqrt{2})(\inf \tilde{a}(Q) +|\partial_x \tilde{a}|_{L^{\infty}(Q)})} |w|_{\mathscr{V}_0^*}. 
\end{equation*}
\end{remark}

\noindent
\underline{\textbf{\textit{Notations for the time-discretization.}}} 
Let $ \tau \in (0, 1) $ be a constant that denotes the time-step size, and let  $ \{ t_i \}_{i = 0}^{\infty} \subset [0, \infty) $ be a sequence of time defined as:
\begin{equation} \label{seqTime}
t_i := i \tau, ~ i = 0, 1, 2, \dots.
\end{equation}
Let $ X $ be a Banach space. Then, for any sequence $ \{ [t_i, \gamma_i] \}_{i = 0}^\infty \subset [0, \infty) \times X $, we define the \emph{forward time-interpolation} $ [\overline{\gamma}]_\tau \in L_\mathrm{loc}^\infty([0, \infty); X)$, the \emph{backward time-interpolation} $ [\underline{\gamma}]_\tau \in L_\mathrm{loc}^\infty([0, \infty); X) $, and the \emph{linear time-interpolation} $ [{\gamma}]_\tau \in W_\mathrm{loc}^{1, 2}([0, \infty); X) $, by letting:
\begin{equation}\label{timeInterp}
\begin{cases}
~ \displaystyle [\overline{\gamma}]_{\tau}(t) := \chi_{(-\infty, 0]}(t) \gamma_0 +\sum_{i = 1}^\infty \chi_{(t_{i -1}, t_i]}(t)\gamma_i,
\\
~ \displaystyle [\underline{\gamma}]_{\tau}(t) := \chi_{(-\infty, 0]}(t) \gamma_0 + \sum_{i = 0}^\infty \chi_{(t_{i}, t_{i +1}]}(t) \gamma_i,
\\
    ~ \displaystyle [{\gamma}]_{\tau}(t) := \sum_{i = 1}^\infty \chi_{[t_{i -1}, t_i)}(t) \left( \frac{t -t_{i -1}}{\tau} \gamma_i +\frac{t_i -t}{\tau} \gamma_{i -1} \right), 
\end{cases} \mbox{in $ X $, for $ t \geq 0 $,}
\end{equation}
respectively.

\begin{remark}\label{Rem.t-discrete}
    For an interval $ I \subset \mathbb{R} $,  a Banach space $ X $, and a constant $ q \in [1, \infty] $, we say that $ L^q(I; X) \subset L_\mathrm{loc}^q(\mathbb{R}; X) $ (resp. $ L_\mathrm{loc}^q(\mathbb{R}; X) \subset L^q(I; X) $) by identifying $ X $-valued functions on $ I $ (resp. on $ \mathbb{R} $) with the zero-extensions onto $ \mathbb{R} $ (resp. the restriction onto $ I $). Besides, under the notations as in \eqref{seqTime} and \eqref{timeInterp}, the following facts can be verified.
\begin{description}
    \item[{\textrm{\hypertarget{Fact0}{(Fact\,0)} $ \bullet $}}]If $ q \in [1, \infty) $, $ \gamma \in L^q(0, T; X) $, and the sequence $ \{ \gamma_i \}_{i = 0}^\infty \subset X $ is given by:
\begin{equation}\label{rApx}
\gamma_i := \frac{1}{\tau} \int_{t_{i -1}}^{t_{i}} \gamma(\varsigma) \, d \varsigma \mbox{ in $ X $,} ~ i = 0, 1, 2, \dots, \mbox{ with } t_{-1} := -\tau, 
\end{equation}
then
\begin{equation*}
    \begin{array}{c}
    [\overline{\gamma}]_\tau \to \gamma, ~ [\underline{\gamma}]_\tau \to \gamma, \mbox{ and } [\gamma]_\tau \to \gamma \mbox{ in  $ L^q(\mathbb{R}; X) $,}
    \\[1ex]
    \mbox{especially in $L^q(0, T; X),$ as $ \tau \downarrow 0 $,}
    \end{array}
\end{equation*}
and
\begin{equation*}    
    \begin{array}{c}
        [\overline{\gamma}]_\tau(t) \to \gamma(t), ~ [\underline{\gamma}]_\tau(t) \to \gamma(t), \mbox{ and }[\gamma]_\tau(t) \to \gamma(t) 
        \\[1ex]
        \mbox{ in $ X $, a.e. $ t \in \mathbb{R} $, as $ \tau \downarrow 0 $.}
    \end{array}
\end{equation*}
\item[~~~~$ \bullet $]If $ X $ is a reflexive Banach space, and $ \gamma \in L^\infty(0, T; X) $, then the sequence $ \{ \gamma_i \}_{i = 0}^\infty \subset X $ given by \eqref{rApx} fulfills that: 
\begin{equation*}
    \sup_{\tau \in (0, 1)} \bigl\{ |[\overline{\gamma}]_\tau|_{L^\infty(0, T; X)}, |[\underline{\gamma}]_\tau|_{L^\infty(0, T; X)}, |[\gamma]_\tau|_{L^\infty(0, T; X)} \bigr\} \leq |\gamma|_{L^\infty(0, T; X)},
\end{equation*}
\begin{equation*}
    \begin{cases}
        [\overline{\gamma}]_\tau \to \gamma, ~ [\underline{\gamma}]_\tau \to \gamma, \mbox{ and } [\gamma]_\tau \to \gamma 
        \\
        \qquad \mbox{in  $ L_\mathrm{loc}^q(\mathbb{R}; X) $, for any $ q \in [1, \infty) $,}
        \\
        \qquad \mbox{weakly-$ * $ in $ L^\infty(\mathbb{R}; X) $,}
    \end{cases} \mbox{as $ \tau \downarrow 0 $,}
\end{equation*}
\begin{equation*}
    \begin{array}{c}
        [\overline{\gamma}]_\tau(t) \to \gamma(t), ~ [\underline{\gamma}]_\tau(t) \to \gamma(t), \mbox{ and } [\gamma]_\tau(t) \to \gamma(t)
        \\[1ex]
        \mbox{in $ X $, a.e. $ t \in \mathbb{R} $, as $ \tau \downarrow 0 $.}
    \end{array}
\end{equation*}
\item[~~~~$ \bullet $]If $ \gamma \in W^{1, \infty}(Q) $, and the sequence $ \{ \gamma_i \}_{i = 0}^\infty \subset W^{1, \infty}(\Omega) $ is given as:
\begin{equation*}
    \gamma_i := \begin{cases}
        \gamma(t_i) \mbox{ in $ W^{1, \infty}(\Omega) $, if $ t_i \leq T $, } 
        \\
        \gamma(t_{i -1}) \mbox{ in $ W^{1, \infty}(\Omega) $, if $ t_{i -1} \leq T < t_i $,}
        \\
        0 \mbox{ in $ W^{1, \infty}(\Omega) $,  otherwise,}
    \end{cases}
    i = 0, 1, 2, \dots,
\end{equation*}
then
\begin{equation*}
    \begin{cases}
        \displaystyle \sup_{\tau \in (0, 1)} \bigl\{ |[\overline{\gamma}]_\tau|_{L^\infty(Q)}, |[\underline{\gamma}]_\tau|_{L^\infty(Q)}, |[\gamma]_\tau|_{C(\overline{Q})} \bigr\} \leq |\gamma|_{C(\overline{Q})},
        \\
        \displaystyle \sup_{\tau \in (0, 1)} \bigl\{ |\partial_x [\overline{\gamma}]_\tau|_{L^\infty(Q)}, |\partial_x [\underline{\gamma}]_\tau|_{L^\infty(Q)}, |\partial_x [\gamma]_\tau|_{L^\infty(Q)} \bigr\} \leq |\partial_x \gamma|_{L^\infty(Q)},
        \\
        \displaystyle \sup_{\tau \in (0, 1)} |\partial_t [\gamma]_\tau|_{L^\infty(Q)} \leq |\partial_t \gamma|_{L^\infty(Q)},
    \end{cases}
\end{equation*}
\begin{equation*}
\begin{cases}
[\overline{\gamma}]_\tau \to \gamma \mbox{ and } [\underline{\gamma}]_\tau \to \gamma, \mbox{ in $ L^\infty(0, T; C(\overline{\Omega})) $,}
\\[1ex]
[\gamma]_\tau \to \gamma \mbox{ in $ C(\overline{Q}) $,}
\end{cases}
\end{equation*}
\begin{equation*}
\begin{cases}
\partial_t [\gamma]_\tau \to \partial_t \gamma \mbox{ weakly-$*$ in $ L^\infty(Q) $,}
\\
\qquad \mbox{and in the pointwise sense a.e. in $ Q $,}
\end{cases}
\end{equation*}
and
\begin{equation*}
\begin{cases}
    \partial_x [\overline{\gamma}]_\tau \to \partial_x \gamma, ~ \partial_x [\underline{\gamma}]_\tau \to \partial_x \gamma, \mbox{ and } \partial_x [\gamma]_\tau \to \partial_x \gamma
\\
\qquad \mbox{weakly-$*$ in $ L^\infty(Q) $,}
\\
\qquad \mbox{and in the pointwise sense a.e. in $ Q $,}
\end{cases}
\mbox{as $ \tau \downarrow 0 $.}
\end{equation*}
\end{description}
\end{remark}

\noindent
\underline{\textbf{\textit{Notations in convex analysis. (cf. \cite[Chapter II]{MR0348562})}}} 
For a proper, lower semi-con- tinuous (l.s.c.), and convex function $ \Psi : X \to (-\infty, \infty] $ on a Hilbert space $ X $, we denote by $ D(\Psi) $ the effective domain of $ \Psi $. Also, we denote by $\partial \Psi$ the subdifferential of $\Psi$. The subdifferential $ \partial \Psi $ corresponds to a generalized derivative of $ \Psi $, and it is known as a maximal monotone graph in the product space $ X \times X $. The set $ D(\partial \Psi) := \bigl\{ z \in X \ |\ \partial \Psi(z) \neq \emptyset \bigr\} $ is called the domain of $ \partial \Psi $. We often use the notation ``$ [w_{0}, w_{0}^{*}] \in \partial \Psi $ in $ X \times X $\,'', to mean that ``$ w_{0}^{*} \in \partial \Psi(w_{0})$ in $ X $ for $ w_{0} \in D(\partial\Psi) $ '', by identifying the operator $ \partial \Psi $ with its graph in $ X \times X $.
\medskip

For Hilbert spaces $X_1, \cdots, X_N$, with $1<N \in \mathbb{N}$, let us consider a proper, l.s.c., and convex function on the product space $X_1 \times \dots \times X_N$:
\begin{equation*}
\tilde{\Psi}: w = [w_1,\cdots,w_N] \in X_1 \times\cdots\times X_N \mapsto \tilde{\Psi}(w)=\tilde{\Psi}(w_1,\cdots,w_N) \in (-\infty,\infty]. 
\end{equation*}
On this basis, for any $i \in \{1,\dots,N\}$, we denote by $\partial_{w_i} \tilde{\Psi}:X_1 \times \cdots \times X_N \to X_i$ a set-valued operator, which maps any $w=[w_1,\dots,w_i,\dots,w_N] \in X_1 \times \dots \times X_i \times \dots \times X_N$ to a subset $ \partial_{w_i} \tilde{\Psi}(w) \subset  X_i $, prescribed as follows:
\begin{equation*}
\begin{array}{rl}
\partial_{w_i}\tilde{\Psi}(w)&=\partial_{w_i}\tilde{\Psi}(w_1,\cdots,w_i,\cdots,w_N)
\\[2ex]
&:= \left\{\begin{array}{l|l}\tilde{w}^* \in X_i & \begin{array}{ll}\multicolumn{2}{l}{(\tilde{w}^*,\tilde{w}-w_i)_{X_i} \le \tilde{\Psi}(w_1,\cdots,\tilde{w},\cdots,w_N)}
\\[0.25ex] 
& \quad -\tilde{\Psi}(w_1,\cdots,w_i,\cdots,w_N), \mbox{ for any $\tilde{w} \in X_i$}\end{array}
\end{array}\right\}.
\end{array}
\end{equation*}
As is easily checked, 
\begin{align}\label{prodSubDif}
    \partial \tilde{\Psi}(w) \subset \partial_{w_1} \tilde{\Psi}(w) \times &\, \cdots \times \partial_{w_N} \tilde{\Psi}(w),
    \\
     \mbox{ for any  $w = [w_1, \dots, w_N]$} &\, \mbox{ $\in X_1 \times \cdots \times X_N $.} \nonumber
\end{align}
But, it should be noted that the converse inclusion of \eqref{prodSubDif} is not true, in general. 

\begin{remark}[Examples of the subdifferential]\label{exConvex}
    As one of the representatives of the subdifferentials, we exemplify the following set-valued function $ \Sgn^N: \mathbb{R}^N \rightarrow 2^{\mathbb{R}^N} $, with $ N \in \mathbb{N} $, which is defined as:
\begin{align*}
    \xi = [\xi_1, & \dots, \xi_N] \in \mathbb{R}^N \mapsto \Sgn^N(\xi) = \Sgn^N(\xi_1, \dots, \xi_N) 
    \\
    & := \left\{ \begin{array}{ll}
            \multicolumn{2}{l}{
                    \displaystyle \frac{\xi}{|\xi|} = \frac{[\xi_1, \dots, \xi_N]}{\sqrt{\xi_1^2 +\cdots +\xi_N^2}}, ~ } \mbox{if $ \xi \ne 0 $,}
                    \\[3ex]
            \mathbb{D}^N, & \mbox{otherwise,}
        \end{array} \right.
    \end{align*}
where $ \mathbb{D}^N $ denotes the closed unit ball in $ \mathbb{R}^N $ centered at the origin. Indeed, the set-valued function $ \Sgn^N $ coincides with the subdifferential of the Euclidean norm $ |{}\cdot{}| : \xi \in \mathbb{R}^N \mapsto |\xi| = \sqrt{\xi_1^2 + \cdots +\xi_N^2} \in [0, \infty) $, i.e.:
\begin{equation*}
\partial |{}\cdot{}|(\xi) = \Sgn^N(\xi), \mbox{ for any $ \xi \in D(\partial |{}\cdot{}|) = \mathbb{R}^N $,}
\end{equation*}
and furthermore, it is observed that:
\begin{equation*}
    \partial  |{}\cdot{}|(0) = \mathbb{D}^N \begin{array}{c} \subseteq_{\hspace{-1.25ex}\mbox{\tiny$_/$}}  
\end{array} [-1, 1]^N = \partial_{\xi_1}  |{}\cdot{}|(0) \times \cdots \times \partial_{\xi_N}  |{}\cdot{}|(0).
\end{equation*}
\end{remark}
\medskip

Finally, we mention about a notion of functional convergence, known as ``Mosco-convergence''. 
 
\begin{definition}[Mosco-convergence: cf. \cite{MR0298508}]\label{Def.Mosco}
    Let $ X $ be an abstract Hilbert space. Let $ \Psi : X \rightarrow (-\infty, \infty] $ be a proper, l.s.c., and convex function, and let $ \{ \Psi_n \}_{n = 1}^\infty $ be a sequence of proper, l.s.c., and convex functions $ \Psi_n : X \rightarrow (-\infty, \infty] $, $ n = 1, 2, 3, \dots $.  Then, it is said that $ \Psi_n \to \Psi $ on $ X $, in the sense of Mosco, as $ n \to \infty $, iff. the following two conditions are fulfilled:
\begin{description}
    \item[(\hypertarget{M_lb}{M1}) The condition of lower-bound:]$ \displaystyle \varliminf_{n \to \infty} \Psi_n(\check{w}_n) \geq \Psi(\check{w}) $, if $ \check{w} \in X $, $ \{ \check{w}_n  \}_{n = 1}^\infty \subset X $, and $ \check{w}_n \to \check{w} $ weakly in $ X $, as $ n \to \infty $. 
    \item[(\hypertarget{M_opt}{M2}) The condition of optimality:]for any $ \hat{w} \in D(\Psi) $, there exists a sequence \linebreak $ \{ \hat{w}_n \}_{n = 1}^\infty  \subset X $ such that $ \hat{w}_n \to \hat{w} $ in $ X $ and $ \Psi_n(\hat{w}_n) \to \Psi(\hat{w}) $, as $ n \to \infty $.
\end{description}
As well as, if the sequence of convex functions $ \{ \tilde{\Psi}_\varepsilon \}_{\varepsilon \in \Xi} $ is labeled by a continuous argument $\varepsilon \in \Xi$ with a infinite set $\Xi \subset \mathbb{R}$ , then for any $\varepsilon_{0} \in \Xi$, the Mosco-convergence of $\{ \tilde{\Psi}_\varepsilon \}_{\varepsilon \in \Xi}$, as $\varepsilon \to \varepsilon_{0}$, is defined by those of subsequences $ \{ \tilde{\Psi}_{\varepsilon_n} \}_{n = 1}^\infty $, for all sequences $\{ \varepsilon_n \}_{n=1}^{\infty} \subset \Xi$, satisfying $\varepsilon_{n} \to \varepsilon_{0}$ as $n \to \infty$.
\end{definition}

\begin{remark}\label{Rem.MG}
Let $ X $, $ \Psi $, and $ \{ \Psi_n \}_{n = 1}^\infty $ be as in Definition~\ref{Def.Mosco}. Then, the following facts hold.
\begin{description}
    \item[(\hypertarget{Fact1}{Fact\,1})](cf. \cite[Theorem 3.66]{MR0773850}, \cite[Chapter 2]{Kenmochi81}) Let us assume that
\begin{center}
$ \Psi_n \to \Psi $ on $ X $, in the sense of  Mosco, as $ n \to \infty $,
\vspace{0ex}
\end{center}
and
\begin{equation*}
\left\{ ~ \parbox{11cm}{
$ [w, w^*] \in X \times X $, ~ $ [w_n, w_n^*] \in \partial \Psi_n $ in $ X \times X $, $ n \in \mathbb{N} $,
\\[1ex]
$ w_n \to w $ in $ X $ and $ w_n^* \to w^* $ weakly in $ X $, as $ n \to \infty $.
} \right.
\end{equation*}
Then, it holds that:
\begin{equation*}
[w, w^*] \in \partial \Psi \mbox{ in $ X \times X $, and } \Psi_n(w_n) \to \Psi(w) \mbox{, as $ n \to \infty $.}
\end{equation*}
    \item[(\hypertarget{Fact2}{Fact\,2})](cf. \cite[Lemma 4.1]{MR3661429}, \cite[Appendix]{MR2096945}) Let $ N \in \mathbb{N} $ denote the dimension constant, and let $  S \subset \mathbb{R}^N $ be a bounded open set. Then, a sequence $ \{ \widehat{\Psi}_n^S \}_{n = 1}^\infty $ of proper, l.s.c., and convex functions on $ L^2(S; X) $, defined as:
        \begin{equation*}
            w \in L^2(S; X) \mapsto \widehat{\Psi}_n^S(w) := \left\{ \begin{array}{ll}
                    \multicolumn{2}{l}{\displaystyle \int_S \Psi_n(w(t)) \, dt,}
                    \\[1ex]
                    & \mbox{ if $ \Psi_n(w) \in L^1(S) $,}
                    \\[2.5ex]
                    \infty, & \mbox{ otherwise,}
                \end{array} \right. \mbox{for $ n = 1, 2, 3, \dots $;}
        \end{equation*}
        converges to a proper, l.s.c., and convex function $ \widehat{\Psi}^S $ on $ L^2(S; X) $, defined as:
        \begin{equation*}
            z \in L^2(S; X) \mapsto \widehat{\Psi}^S(z) := \left\{ \begin{array}{ll}
                    \multicolumn{2}{l}{\displaystyle \int_S \Psi(z(t)) \, dt, \mbox{ if $ \Psi(z) \in L^1(S) $,}}
                    \\[2ex]
                    \infty, & \mbox{ otherwise;}
                \end{array} \right. 
        \end{equation*}
        on $ L^2(S; X) $, in the sense of Mosco, as $ n \to \infty $. 
\end{description}
\end{remark}
\begin{remark}[Example of Mosco-convergence]\label{Rem.ExMG}
    For any $ \varepsilon \geq 0 $, let $ f_\varepsilon : \mathbb{R} \longrightarrow [0, \infty) $ be a continuous and convex function, defined as:
    \begin{equation}\label{f_eps}
        f_\varepsilon : \xi \in \mathbb{R} \mapsto f_\varepsilon(\xi) := \sqrt{\varepsilon^2 +|\xi|^2} \in [0, \infty).
    \end{equation}
    Then, due to the uniform estimate:
    \begin{equation}\label{ev.M00}
        \begin{array}{c}
            \bigl| f_\varepsilon(\xi) -f_{\tilde{\varepsilon}}(\xi) \bigr| = \bigl| \sqrt{\varepsilon^2 +|\xi|^2 } -\sqrt{\tilde{\varepsilon}^2 +|\xi|^2} \bigr| \leq |\varepsilon -\tilde{\varepsilon}|, 
            \\[1ex]
            \mbox{for all $ \xi \in \mathbb{R} $, and $ \varepsilon, \tilde{\varepsilon} \geq 0 $,}
        \end{array}
    \end{equation}
    we easily see that:
    \begin{equation*}
        f_\varepsilon \to f_0 ~ (= |{}\cdot{}|) \mbox{ on $ \mathbb{R} $, in the sense of Mosco, as $ \varepsilon \downarrow 0 $.}
    \end{equation*}
    In addition, for any $ \varepsilon > 0 $, it can be said that the subdifferential $ \partial f_\varepsilon $ coincides with the single-valued function of usual differential:
    \begin{equation*}
        f_\varepsilon' : \xi \in \mathbb{R} \mapsto f_\varepsilon'(\xi) = \frac{\xi}{\sqrt{\varepsilon^2 +|\xi|^2}} \in \mathbb{R}.
    \end{equation*}
\end{remark}

\section{{Auxiliary results}}

In this Section, we prepare some auxiliary results for our study. The auxiliary results are discussed through the following two Subsections.
\begin{description}
    \item[\textmd{$\S$\,2.1}]Abstract theory for the state-system (S)$_\varepsilon$;
    \item[\textmd{$\S$\,2.2}]Mathematical theory for the linearized system of (S)$_\varepsilon$. 
\end{description}

\subsection{Abstract theory for the state-system (S)$_\varepsilon$}

In this Subsection, we refer to \cite[Appendix]{AKSY2020} to overview the abstract theory of nonlinear evolution equation, which enables us to handle the state-systems (S)$_\varepsilon$, for all $ \varepsilon \geq 0 $, in a unified fashion. 
\medskip

The general theory consists of the following two Propositions. 

\begin{proposition}[{cf. \cite[Lemma 8.1]{AKSY2020}}]\label{Lem.CP}
    Let $ \{ \mathcal{A}_0(t) \, | \, t \in [0, T] \} \subset \mathscr{L}(X) $ be a class of time-dependent bounded linear operators, let $ \mathcal{G}_0 : X \longrightarrow X $ be a given nonlinear operator, and let $ \Psi_0 : X \longrightarrow [0, \infty] $ be a non-negative, proper, l.s.c., and convex function, fulfilling the following conditions: 
    \begin{itemize}
        \item[\textmd{(cp.0)}]$ \mathcal{A}_0(t) \in \mathscr{L}(X) $ is positive and selfadjoint, for any $ t \in [0, T] $, and it holds that
            \begin{equation*}
                (\mathcal{A}_0(t) w, w)_X \geq \kappa_0 |w|_X^2 ,\ \mbox{for any}\ w \in X,
            \end{equation*}
            with some constant $\kappa_{0} \in (0, 1)$, independent of $t \in [0, T]$ and $w \in X$.
        \item[\textmd{(cp.1)}]$ \mathcal{A}_0 : [0, T] \longrightarrow \mathscr{L}(X) $ is Lipschitz continuous, so that $ \mathcal{A}_0 $ admits the (strong) time-derivative $ \mathcal{A}_0'(t) \in \mathscr{L}(X) $ a.e. in $ (0, T) $, and  
            \begin{equation*}
                A_T^* := \mathrm{ess} \sup_{\hspace{-3ex}t \in (0, T)} \left\{ \max \{ |\mathcal{A}_0(t)|_{\mathscr{L}(X)}, |\mathcal{A}_0'(t)|_{\mathscr{L}(X)} \} \right\} < \infty;
            \end{equation*}
        \item[\textmd{(cp.2)}]$ \mathcal{G}_0 : X \longrightarrow X $ is a Lipschitz continuous operator with a Lipschitz constant $ L_0 $, and $ \mathcal{G}_0 $ has a $ C^1 $-potential functional $ \widehat{\mathcal{G}}_0 : X \longrightarrow \mathbb{R} $, so that the G\^{a}teaux derivative $ \widehat{\mathcal{G}}_0'(w) \in X^* $ $ (= X) $ at any $ w \in X $ coincides with $ \mathcal{G}_0(w) \in X $;  
        \item[\textmd{(cp.3)}]$ \Psi_0 \geq 0 $ on $ X $, and the sublevel set $ \bigl\{ w \in X \, \bigl| \, \Psi_0(w) \leq r \bigr\} $ is compact in $ X $, for any $ r \geq 0 $.
    \end{itemize}
    Then, for any initial data $ w_0 \in D(\Psi_0) $ and a forcing term $ \mathfrak{f}_0 \in L^2(0, T; X) $, the following Cauchy problem of evolution equation:
    \begin{equation*}
        (\mathrm{CP})~~
        \begin{cases}
            \mathcal{A}_0(t) w'(t) +\partial \Psi_0(w(t)) +\mathcal{G}_0(w(t)) \ni \mathfrak{f}_0(t) \mbox{ in $ X $, ~ $ t \in (0, T) $,}
            \\
            w(0) = w_0 \mbox{ in $ X $;}
        \end{cases}
    \end{equation*}
    admits a unique solution $ w \in L^2(0, T; X) $, in the sense that:
    \begin{equation*}
        w \in W^{1, 2}(0, T; X), ~  \Psi_0(w) \in L^\infty(0, T),
    \end{equation*}
    and
    \begin{equation*}
        \begin{array}{c}
            \displaystyle \bigl( \mathcal{A}_0(t)w'(t) +\mathcal{G}_0(w(t)) -\mathfrak{f}_0(t), w(t) -\varpi \bigr)_X +\Psi_0(w(t)) \leq \Psi_0(\varpi), 
            \\[1ex]
            \mbox{for any $ \varpi \in D(\Psi_0) $, a.e. $ t \in (0, T) $.} 
        \end{array}
    \end{equation*}
    Moreover, both $ t \in [0, T] \mapsto \Psi_0(w(t)) \in [0, \infty) $ and $ t \in [0, T] \mapsto \widehat{\mathcal{G}}_0(w(t)) \in \mathbb{R} $ are absolutely continuous functions in time, and
    \begin{equation*}
        \begin{array}{c}
            \displaystyle |\mathcal{A}_0(t)^{\frac{1}{2}}w'(t)|_X^2 +\frac{d}{dt} \left( \Psi_0(w(t)) +\widehat{\mathcal{G}}_0(w(t)) \right) = (\mathfrak{f}_0(t), w'(t))_X,
            \\[1ex]
            \mbox{for a.e. $ t \in (0, T) $.}
        \end{array}
    \end{equation*}
\end{proposition}
\begin{proposition}[{cf. \cite[Lemma 8.2]{AKSY2020}}]\label{Lem.CP02}
    Under the notations $ \mathcal{A}_0 $, $ \mathcal{G}_0 $, $ \Psi_0 $, and assumptions (cp.0)--(cp.3), as in the previous Proposition \ref{Lem.CP}, let us fix $ w_0 \in D(\Psi_0) $ and $ \mathfrak{f}_0 \in L^2(0, T; X) $, and take the unique solution $ w \in L^2(0, T; X) $ to the Cauchy problem (CP). Let $ \{ \Psi_n \}_{n = 1}^\infty $, $ \{ w_{0, n} \}_{n = 1}^\infty \subset X $, and $ \{ \mathfrak{f}_n \}_{n = 1}^\infty $ be, respectively, a sequence of proper, l.s.c., and convex functions on $ X $, a sequence of initial data in $ X $, and a sequence of forcing terms in $ L^2(0, T; X) $, such that:
    \begin{itemize}
        \item[\textmd{(cp.4)}]$ \Psi_n \geq 0 $ on $ X $, for $ n = 1, 2, 3, \dots $, and the union $ \bigcup_{n = 1}^\infty \bigl\{  w \in X \, \bigl| \, \Psi_n(w) \leq r \bigr\} $ of sublevel sets is relatively compact in $ X $, for any $ r \geq 0 $;
    \item[\textmd{(cp.5)}]$ \Psi_n $ converges to $ \Psi_0 $ on $ X $, in the sense of Mosco, as $ n \to \infty $; 
    \item[\textmd{(cp.6)}]$ \sup_{n \in \mathbb{N}} \Psi_n(w_{0, n}) < \infty $, and $ w_{0, n} \to w_0 $ in $ X $, as $ n \to \infty $;
    \item[\textmd{(cp.7)}]$ \mathfrak{f}_n \to \mathfrak{f}_0 $ weakly in $ L^2(0, T; X) $, as $ n \to \infty $.
    \end{itemize}
    Let $ w_n \in L^{2}(0, T; X)$ be the solution to the Cauchy problem (CP), for the initial data $ w_{0, n} \in D(\Psi_n) $ and forcing term $ \mathfrak{f}_n \in L^2(0, T; X) $. Then, 
    \begin{equation*}
        \begin{array}{c}
            \displaystyle w_n \to w\mbox{ in $ C([0, T]; X) $, weakly in $ W^{1, 2}(0, T; X) $,}
            \\[1ex]
            \displaystyle \int_0^T \Psi_n(w_n(t)) \, dt \to \int_0^T \Psi_0(w(t)) \, dt, 
            \mbox{ as $ n \to \infty $,}
        \end{array}
    \end{equation*}
    and
    \begin{equation*}
        \bigl| \Psi_0(w) \bigr|_{C([0, T])} \leq  \sup_{n \in \mathbb{N}} \, \bigl| \Psi_n(w_n) \bigr|_{C([0, T])} < \infty.
    \end{equation*}
\end{proposition}

In this paper, the readers are recommended  to see \cite[Appendix]{AKSY2020} for the detailed proofs of the above Propositions \ref{Lem.CP} and \ref{Lem.CP02}. Roughly summarized, these Propositions can be obtained by means of modified (mixed and reduced) methods of the existing theories, such as \cite{MR0348562,MR2582280,Kenmochi81}.


\subsection{Mathematical theory for the linearized system of (S)$_\varepsilon$}

In this Subsection, we set up auxiliary results for linearized systems of (S)$_\varepsilon$, which are associated with the first necessary optimality conditions in our optimal control problems (OP)$_\varepsilon $, for $ \varepsilon \geq 0$. The linearized systems are generally reduced to the following type of parabolic initial-boundary value problem, denoted by (\hyperlink{P}{P}). 
\pagebreak

(\hypertarget{P}{P}):
\vspace{0ex}
\begin{equation*}
\left\{\parbox{9.5cm}{
    $\partial_{t}p - \partial_{x}^{2} p + \mu(t, x)p + \omega(t, x)\partial_{x} z = h(t, x)$, $(t, x) \in Q,$
    \\[1ex]
    $ \partial_t p_\Gamma(t, \ell) +(-1)^{\ell -1} \partial_{x} p_{_{|_\Gamma}}(t, \ell) = h_\Gamma(t, \ell)$, ~ $ (t, \ell) \in \Sigma $,
    \\[1ex]
    $ p_{_{|_\Gamma}} = p_\Gamma $ on $ \Sigma $,
    \\[1ex]
    $p(0, x) = p_{0}(x)$, $x \in \Omega$;
}\right. 
\end{equation*}
\begin{equation*}
\left\{ \parbox{9.5cm}{
    $a(t, x)\partial_{t}z + b(t, x)z- \partial_{x} \bigl( A(t, x)\partial_{x}z + \nu^{2} \partial_{x}z + \omega(t, x)p \bigr) $ 
    \\
    $ ~~~~~~= k(t, x)$, $(t, x) \in Q $,
    \\[1ex]
    $ z(t, x) = 0 $, $ (t, x) \in \Sigma $,
    \\[1ex]
    $ z(0, x) = z_{0}(x) $, $ x \in \Omega $.
} \right. 
\end{equation*}
\noindent
This system is a key-problem for the G$\hat{\mbox{a}}$teaux differential of the cost functional $\mathcal{J}_{\varepsilon}$. 
In the context, $[a, b, \lambda, \omega, A] \in [\mathscr{H}]^{5}$ is a given quintet of functions which belongs to a class $\mathscr{S} \subset [\mathscr{H}]^{5}$, defined as:
\begin{equation}\label{P01}
    \mathscr{S} := \left\{ \begin{array}{l|l} [\tilde{a}, \tilde{b}, \tilde{\mu}, \tilde{\omega}, \tilde{A}] \in [\mathscr{H}]^{5} & 
        \parbox{6.75cm}{
            $ \tilde{a} \in W^{1, \infty}(Q) $ with $ \log \tilde{a} \in L^{\infty}(Q) $, $ [\tilde{b}, \tilde{\omega}]\in [L^{\infty}(Q)]^{2} $,  $ \tilde{\mu} \in L^{\infty}(0, T; H) $, and $  \tilde{A} \in L^{\infty}(Q) $  with $ \tilde{A} \geq 0 $ a.e. in $ Q $
}
\end{array} 
\right\}.
\end{equation}
Also, $[\bm{p}_{0}, z_{0}] = [p_0, p_{\Gamma, 0}, z_0] \in \mathbb{W} \times H $ with $ \bm{p}_{0} = [p_0, p_{\Gamma, 0}] $ and $[\bm{h}, k] = [h, h_\Gamma, k] \in \mathfrak{X} \times \mathscr{V}_0^* $ with $ \bm{h} = [h, h_\Gamma] $ are the initial triplet and forcing triplet in the system (\hyperlink{P}{P}), respectively.

\begin{remark}\label{Rem.delta(a)}
    If $ [a, b, \mu, \omega, A] \in \mathscr{S} $, then the condition:
    \begin{equation*}
        a \in W^{1, \infty}(Q) \mbox{ with } \log a \in L^\infty(Q),
    \end{equation*}
    brought by \eqref{P01}, implies the no degeneration property:
    \begin{equation}\label{delta(a)}
        \delta_*(a) := \inf a(Q) > 0,
    \end{equation}
    of the coefficient $ a $ in the system (\hyperlink{P}{P}).
\end{remark}

Now, as the key-properties of the system (\hyperlink{P}{P}), we can verify the following three Theorems. 

\begin{theorem}
    \label{Prop01}
    Let us assume $[a, b, \mu, \omega, A] \in \mathscr{S}$, $[\bm{p}_{0}, z_{0}] = [p_0, p_{\Gamma, 0}, z_0] \in \mathbb{W} \times H $ with $ \bm{p}_0 = [p_0, p_{\Gamma, 0}]  $, and $[\bm{h}, k] = [h, h_\Gamma, k] \in \mathfrak{X} \times \mathscr{V}_0^* $ with $ \bm{h} = [h, h_\Gamma] $. Then, the system (\hyperlink{P}{P}) admits a unique solution $ [\bm{p}, z] = [p, p_\Gamma, z] \in \mathfrak{X} \times \mathscr{H} $ with $ \bm{p} = [p, p_\Gamma] $, in the sense that: 
    \begin{equation*}
        \begin{cases}
            \bm{p} = [p, p_\Gamma] \in W^{1, 2}(0, T; \mathbb{X}) \cap L^\infty(0, T; \mathbb{W}) \subset C(\overline{Q}) \times C(\overline{\Sigma}), 
            \\
            z \in W^{1, 2}(0, T; V_{0}^{*})\cap \mathscr{V}_0 \subset C([0, T]; H);
        \end{cases}
    \end{equation*}
    \begin{align*}
        & \displaystyle ( \partial_{t} \bm{p}(t),  \bm{\varphi} )_{\mathbb{X}} 
        +(\partial_{x}p(t), \partial_{x}\varphi)_{H}  +(\mu(t)p(t), \varphi)_{H} +(\omega(t)\partial_{x}z(t), \varphi)_{H} 
        \nonumber
        \\
        & \qquad = (\bm{h}(t), \bm{\varphi})_{\mathbb{X}}, 
        \mbox{ for any $ \bm{\varphi} = [\varphi, \varphi_\Gamma] \in \mathbb{W} $, a.e. $ t \in (0, T) $,}
        \\
        & \mbox{subject to $ \bm{p}(0) = [p(0), p_\Gamma(0)] = \bm{p}_0 = [p_0, p_{\Gamma, 0}] $ in $ \mathbb{X} $;}
        \nonumber
    \end{align*}
    and
    \begin{align*}
        & \langle a(t)\partial_{t}z (t), \psi \rangle_{V_{0}} + (b(t)z(t), \psi)_{H} 
        \nonumber
        \\
        & \qquad +\bigl( A(t)\partial_{x}z(t) + \nu^{2}\partial_{x}z(t) +p(t)\omega(t), \partial_{x}\psi \bigr)_H = \langle k(t), \psi \rangle_{V_{0}},
        \\
        & \mbox{for any $ \psi \in V_{0} $, a.e. $ t \in (0, T) $, subject to $ z(0) = z_0 $ in $ H $.}
        \nonumber
        \end{align*}
\end{theorem}
\begin{theorem}
    \label{Prop02}
    Let us take arbitrary $ [a, b, \mu, \omega, A] \in \mathscr{S} $, $ [\bm{p}_0, z_0] = [p_0, p_{\Gamma, 0}, z_0] \in \mathbb{W} \times H $ with $ \bm{p}_0 = [p_0, p_{\Gamma, 0}] $, and $ [\bm{h}, k] = [h, h_\Gamma, k] \in \mathfrak{X} \times \mathscr{V}_0^* $ with $ \bm{h} = [h, h_\Gamma] $, and let us denote by $ [\bm{p}, z] = [p, p_\Gamma, z] \in \mathfrak{X} \times \mathscr{H} $ with $\bm{p} = [p, p_\Gamma]$ the solution to (\hyperlink{P}{P}). Additionally, let $ \delta_*(a) $ be the positive constant as in \eqref{delta(a)}. Then, the following two items hold.
    \begin{description}
        \item[\textmd{(I)}]Let $ C_0^* $ be a positive constant, 
            defined as:
    \begin{equation}\label{C0*}
        C_0^* := \frac{16 \bigl( 1 +| a|_{W^{1, \infty}(Q)} +|b|_{L^\infty(Q)} +|\mu|_{L^\infty(0, T; H)}^2 +|\omega|_{L^\infty(Q)}^2 \bigr)}{\min \{1, \nu^2, \delta_*(a)\}}.
    \end{equation}
            Then, it is estimated that:
    \begin{align}\label{est_I-B} 
        \frac{d}{dt} & \bigl( |\bm{p}(t)|_{\mathbb{X}}^2 +|\sqrt{a(t)}z(t)|_H^2 \bigr) +\bigl( |\bm{p}(t)|_{\mathbb{W}}^2 +\nu^2 |z(t)|_{V_0}^2 \bigr)
        \nonumber
        \\
        & \leq C_0^* \bigl( |\bm{p}(t)|_{\mathbb{X}}^2 +|\sqrt{a(t)}z(t)|_H^2 \bigr)+C_0^*  \bigl( |\bm{h}(t)|_{\mathbb{X}}^2 +|k(t)|_{V_0^*}^2 \bigr),
        \\
        & \qquad \mbox{for a.e. $ t \in (0, T) $.}
        \nonumber
    \end{align}
\item[\textmd{(II)}]Let $ C_0^* $ be the positive constant given in \eqref{C0*}, and let $ C_\ell^* $, $ \ell = 1, 2 $, be positive constants, 
    defined as:
        \begin{align}\label{C1*}
            & \begin{cases}
                C_1^* := 4 (C_0^*)^2 e^{\frac{3}{2}C_0^* T}, ~ 
                \\[1ex]
                \displaystyle C_2^* := 4(C_0^*)^6 e^{\frac{3}{2} C_0^* T}(1 + |a|_{W^{1, \infty}(Q)})^2\cdot 
                \\[1ex]
                \displaystyle \qquad \qquad \cdot(1 +\nu + |b|_{L^\infty(Q)} +|\omega|_{L^\infty(Q)} +|A|_{L^\infty(Q)})^2. 
            \end{cases}
        \end{align}
            Then, it is estimated that:
        \begin{align}\label{est_I-B01}
            & \begin{cases}
                |\partial_t \bm{p}|_{\mathfrak{X}}^2  
                +|p|_{L^\infty(0, T; V)}^2 
            \leq  C_1^* \bigl( |\bm{p}_0|_{\mathbb{W}}^2 +|\sqrt{a(0, \cdot)} z_0|_H^2 +|\bm{h}|_{\mathfrak{X}}^2 +|k|_{\mathscr{V}_0^*}^2 \bigr), 
                \\[1ex]
                |\partial_t z|_{\mathscr{V}_0^*}^2 \leq C_2^*\bigl( |\bm{p}_0|_{\mathbb{W}}^2 +|\sqrt{a(0, \cdot)} z_0|_H^2 +|\bm{h}|_{\mathfrak{X}}^2 +|k|_{\mathscr{V}_0^*}^2 \bigr).
            \end{cases}
        \end{align}
\end{description}
\end{theorem}
\begin{remark}\label{Rem.C_0*}
    By applying Gronwall's lemma to the inequality in Theorem \ref{Prop02} (I), we also estimate that: 
    \begin{align*}
        \bigl( |\bm{p}|_{C([0, T]; \mathbb{X})}^2 & +|\sqrt{a} z|_{C([0, T]; H)}^2 \bigr) +\bigl( |\bm{p}|_{\mathfrak{W}}^2 +\nu^2 |z|_{\mathscr{V}_0}^2 \bigr)
        \nonumber
        \\
        & \leq 2 C_0^* e^{C_0^* T} \bigl( |\bm{p}_0|_{\mathbb{X}}^2 +|\sqrt{a(0, \cdot)} z_0|_H^2 +|\bm{h}|_{\mathfrak{X}}^2 +|k|_{\mathscr{V}_0^*}^2 \bigr).
    \end{align*}
\end{remark}
\begin{theorem}
    \label{Prop03}
    Let us assume: 
\begin{subequations}\label{ASY01}
    \begin{equation}\label{ASY01-01}
        [a, b, \mu, \omega, A] \in \mathscr{S}, ~ \{ [a^n, b^n, \mu^n, \omega^n, A^n] \}_{n = 1}^\infty \subset \mathscr{S}, 
    \end{equation}
    \begin{align}\label{ASY01-02}
        [a^n, \partial_t a^n, & \partial_x a^n, b^n, \omega^n, A^n] \to [a, \partial_t a, \partial_x a, b, \omega, A] \mbox{ weakly-$*$ in $ [L^\infty(Q)]^6 $,}
        \nonumber
        \\
        & \mbox{and in the pointwise sense a.e. in  $ Q $, as $ n \to \infty $,}
    \end{align}
    and
    \begin{equation}\label{ASY01-03}
        \begin{cases}
            \mu^n \to \mu \mbox{ weakly-$*$ in $ L^\infty(0, T; H) $,}
            \\
            \mu^n(t) \to \mu(t) \mbox{ in $ H $, in the pointwise sense, }
            \\
            \qquad \qquad \qquad \qquad \mbox{ for a.e. $ t \in (0, T) $,}
        \end{cases}
        \mbox{as $ n \to \infty $.}
    \end{equation}
\end{subequations}
    Let us assume $ [\bm{p}_0, z_0] = [p_0, p_{\Gamma, 0}, z_0] \in \mathbb{W} \times H $ with $ \bm{p}_0 = [p_0, p_{\Gamma, 0}] $, $ [\bm{h}, k] = [h, h_\Gamma, k] \in \mathfrak{X} \times \mathscr{V}_0^* $ with $ \bm{h} = [h, h_\Gamma] $, and let us denote by $ [\bm{p}, z] = [p, p_\Gamma, z] \in \mathfrak{X} \times \mathscr{H} $ with $ \bm{p} = [p, p_\Gamma] $ the solution to (\hyperlink{P}{P}), for the initial triplet $ [\bm{p}_0, z_0] = [p_0, p_{\Gamma, 0}, z_0] $ and forcing triplet $ [\bm{h}, k]= [h, h_\Gamma, k] $. Also, for any $ n \in \mathbb{N} $, let us assume $ [\bm{p}_0^n, z_0^n] = [p_0^n, p_{\Gamma, 0}^n, z_0^n] \in \mathbb{W} \times H $ with $ \bm{p}_{0}^n = [p_0^n, p_{\Gamma, 0}^n] $, $ [\bm{h}^n, k^n] = [h^n, h_\Gamma^n, k^n] \in  \mathfrak{X} \times \mathscr{V}_0^* $ with $ \bm{h}^n = [h^n, h_\Gamma^n]  $, and let us denote by $ [\bm{p}^n, z^n] = [p^n, p_\Gamma^n, z^n] \in \mathfrak{X} \times \mathscr{H} $ with $ \bm{p}^n = [p^n, p_\Gamma^n] $ the solution to (\hyperlink{P}{P}), for the initial triplet $ [\bm{p}_0^n, z_0^n] = [p_0^n, p_{\Gamma, 0}^n, z_0^n] $ and forcing triplet $ [\bm{h}^n, k^n] = [h^n, h_\Gamma^n, k^n] $. Then, the convergences of given data:
            \begin{equation}\label{pr03-01}
                \begin{cases}
                    [\bm{p}_0^n, z_0^n] \to [\bm{p}_0, z_0]  \mbox{ weakly in $ \mathbb{W} \times H $},
                    \\
                    [\bm{h}^n, k^n] \to [\bm{h}, k] \mbox{ weakly in $ \mathfrak{X} \times \mathscr{V}_0^* $,}
                \end{cases}
                \mbox{as $ n \to \infty $,}
            \end{equation}
            implies the convergence of solutions, in the sense that:
            \begin{align}\label{pr03-02}
                [\bm{p}^n, z^n] &\, \to [\bm{p}, z]  \mbox{ in $ [C(\overline{Q}) \times C(\overline{\Sigma})] \times \mathscr{H} $, weakly in $ \mathfrak{W} \times \mathscr{V}_0 $,}
                \nonumber
                \\
                & \mbox{and weakly in $ W^{1, 2}(0, T; \mathbb{X}) \times W^{1, 2}(0, T; V_0^*) $, as $ n \to \infty $.}
            \end{align}
\end{theorem}

\begin{remark}[\textmd{Review of Theorems \ref{Prop01}--\ref{Prop03}}]\label{appendix}
Let us define:
\begin{equation*}
    \mathfrak{Y} := [ W^{1, 2}(0, T; \mathbb{X}) \cap L^\infty(0, T; \mathbb{W}) ] \times [ W^{1, 2}(0, T; V_0^*) \cap \mathscr{V}_0 ].
\end{equation*}
as a Banach space endowed with the following norm 
\begin{align*}
    \bigl| [\bm{\tilde{p}}, z] \bigr|_{\mathfrak{Y}} & := \bigl( \bigl| [\partial_t \tilde{p}, \partial_t \tilde{p}_\Gamma] \bigr|_{\mathfrak{X}}^2 +\bigl| [\tilde{p}, \tilde{p}_\Gamma] \bigl|_{L^\infty(0, T; \mathbb{W})}^2 +|\partial_t \tilde{z}|_{\mathscr{V}_0^*}^2 +|\tilde{z}|_{\mathscr{V}_0}^2 \bigr)^{\frac{1}{2}},
    \\
    & \mbox{for any $ [\bm{\tilde{p}}, \tilde{z}] = [\tilde{p}, \tilde{p}_\Gamma, \tilde{z}] \in \mathfrak{Y} $ with $ \bm{\tilde{p}} = [\tilde{p}, \tilde{p}_\Gamma] $.}
\end{align*}
The Banach space $ \mathfrak{Y} $ is to characterize the regularity of solution to the linearized system of  (S)$_\varepsilon$. Due to the compactness theory of Aubin's type (cf. \cite[Corollary 4]{MR0916688}), this Banach space $ \mathfrak{Y} $ is compactly embedded into the Banach space $ [ C(\overline{Q}) \times C(\overline{\Sigma}) ] \times \mathscr{H} $. 

    Now, for any quintet of functions $ [a, b, \mu, \omega, A] \in \mathscr{S} $, the first and second Theorems \ref{Prop01} and \ref{Prop02} will enable us to define a bounded linear operator $ \mathcal{P} = P(a, b, \mu, \omega, A) : [ \mathbb{W} \times H ] \times [\mathfrak{X} \times \mathscr{V}_0^*] \longrightarrow \mathfrak{Y} $, which maps any pair $ \bigl[ [\bm{p}_0, z_0], [\bm{h}, k] \bigr] \in [\mathbb{W} \times H] \times [\mathfrak{X} \times \mathscr{V}_0^*] $ of the initial triplet $ [\bm{p}_0, z_0] = [p_0, p_{\Gamma, 0}, z_0] \in \mathbb{W} \times H $ with $ \bm{p}_0 = [p_0, p_{\Gamma, 0}]  $ and the forcing triplet $ [\bm{h}, k] = [h, h_\Gamma, k] \in \mathfrak{X} \times \mathscr{V}_0^* $ with $ \bm{h} = [h, h_\Gamma] $, to the solution $ [\bm{p}, z] = [p, p_\Gamma, z] \in \mathfrak{Y} $ with $\bm{p} = [p, p_\Gamma]$ to the linear system (\hyperlink{P}{P}). Moreover, the third Theorem \ref{Prop03} will be to guarantee the continuous dependence of the solution operator $ \mathcal{P} = \mathcal{P}(a, b, \mu, \omega, A) $, in the following sense:
    \begin{align*}
        & \mathcal{P}(a^n, b^n, \mu^n, \omega^n, A^n) \bigl[ [\bm{p}_0^n, z_0^n], [\bm{h}^n, k^n] \bigr] \to \mathcal{P}(a, b, \mu, \omega, A) \bigl[ [\bm{p}_0, z_0], [\bm{h}, k] \bigr]
        \\
        & \quad \mbox{in the topologies as in \eqref{pr03-02}, whenever \eqref{ASY01} and \eqref{pr03-01} are fulfilled.
        }
        \end{align*}
\end{remark}

The proofs of the three Theorems \ref{Prop01}--\ref{Prop03} will be given in the appendix, that is assigned to the last Section 7 of this paper.

\section{Main Theorems}

We begin by setting up the assumptions needed in our Main Theorems. 
\begin{itemize}
\item[\textmd{(\hypertarget{A0l}{A0})}]
    $\nu>0$ is a fixed constant. Let $[\bm{\eta}_\mathrm{ad}, \theta_\mathrm{ad} ] = [\eta_\mathrm{ad}, \eta_{\Gamma, \mathrm{ad}}, \theta_\mathrm{ad}] \in \mathfrak{X} \times \mathscr{H} $ with $\bm{\eta}_\mathrm{ad} = [\eta_\mathrm{ad}, \eta_{\Gamma, \mathrm{ad}}] $ be a fixed triplet of functions, called the \emph{admissible target profile}. 
\item[\textmd{(\hypertarget{A1l}{A1)}}]
    $\alpha:\mathbb{R} \longrightarrow (0, \infty)$ and $\alpha_{0}: Q \longrightarrow (0, \infty)$ are Lipschitz continuous functions, such that:
        \begin{itemize}
            \item $\alpha \in C^2(\mathbb{R})$, with the first derivative $ \alpha' = \frac{d \alpha}{d \eta} \in C^1(\mathbb{R}) \cap L^\infty(\mathbb{R})$ and the second one $ \alpha'' = \frac{d^2 \alpha}{d\eta^2} \in C(\mathbb{R})$;
            \item $\alpha'(0) = 0$, $ \alpha'' \geq 0 $ on $ \mathbb{R} $, and $\alpha \alpha'$ is a Lipschitz continuous function on $\mathbb{R}$; 
            \item $ \alpha \geq \delta_* $ on $ \mathbb{R} $, and $ \alpha_0 \geq \delta_* $ on $ \overline{Q} $, for some constant $ \delta_* \in (0, 1) $.
        \end{itemize}
\item[\textmd{(\hypertarget{A2l}{A2})}]For any $ \varepsilon \geq 0 $, let $ f_\varepsilon : \mathbb{R} \longrightarrow [0, \infty) $ be the convex function, defined in \eqref{f_eps}.
\item[\textmd{(\hypertarget{A3l}{A3})}]
    $g : \mathbb{R} \longrightarrow  \mathbb{R}$ is a $C^{1}$-function, which is a Lipschitz continuous on $\mathbb{R}$. Also $g$ has a nonnegative primitive $ 0 \leq G \in C^{2}(\mathbb{R})$, i.e. the derivative $ G'= \frac{dG}{d\eta} $ coincides with $ g $ on $ \mathbb{R} $.
\end{itemize}
\bigskip

Now, the Main Theorems of this paper are stated as follows:
\begin{main}\label{mainTh01}
    Let us assume (\hyperlink{A0l}{A0})--(\hyperlink{A3l}{A3}). Let us fix a constant $ \varepsilon \geq 0 $,  an initial triplet $ [\bm{\eta}_0, \theta_0] = [\eta_0, \eta_{\Gamma, 0}, \theta_0] \in \mathbb{W} \times V_0 $ with $ \bm{\eta}_0 = [\eta_0, \eta_{\Gamma, 0}] $, and a forcing triplet $ [\bm{u}, v] = [u, u_\Gamma, v] \in \mathfrak{X} \times \mathscr{H} $ with $\bm{u} = [u, u_\Gamma]$. 
    Then, the following two items hold.
    \begin{itemize}
        \item[\textmd{(I-A)}]
            The state-system (S)$_{\varepsilon}$ admits a unique solution $[\bm{\eta}, \theta] = [\eta, \eta_\Gamma, \theta] \in \mathfrak{X} \times \mathscr{H} $ with $ \bm{\eta} = [\eta, \eta_\Gamma] $, in the sense that:
        \begin{equation*}
            \begin{cases}
            \bm{\eta} = [\eta, \eta_\Gamma] \in W^{1, 2}(0, T; \mathbb{X}) \cap L^\infty(0, T; \mathbb{W}) \subset C(\overline{Q}) \times C(\overline{\Sigma}),
                \\[0ex]
                \theta \in W^{1, 2}(0, T; H) \cap L^\infty(0, T; V_0) \subset C(\overline{Q});
        \end{cases}
        \end{equation*}
        \begin{equation*}
                \begin{array}{c}
                    \displaystyle \bigl( \partial_{t}\bm{\eta}(t), \bm{\varphi} \bigr)_{\mathbb{X}} + \bigl( \partial_{x}\eta(t), \partial_{x}\varphi \bigr)_{H} +\bigl( g(\eta(t)), \varphi \bigr)_{H} 
                    \\[1ex]
                    +\bigl( \alpha'(\eta(t))f_{\varepsilon}(\partial_{x}\theta(t)), \varphi \bigr)_{H}  = \bigl( L u(t), \varphi \bigr)_{H} + \bigl( L_\Gamma u_\Gamma(t), \varphi_\Gamma \bigr)_{H_\Gamma}, 
                    \\[1ex]
                    \mbox{for any $ \bm{\varphi} = [\varphi, \varphi_\Gamma] \in \mathbb{W} $, a.e. }  t \in (0, T) , 
                    \\[1ex]
                    \mbox{ subject to } \bm{\eta}(0) = [\eta(0), \eta_\Gamma(0)] = \bm{\eta}_0 = [\eta_0, \eta_{\Gamma, 0}] \mbox{ in $ \mathbb{X} $;}
                \end{array}
        \end{equation*}
and
        \begin{equation*}
                \begin{array}{c}
                    \displaystyle\bigl( \alpha_{0}(t)\partial_{t}\theta(t), \theta(t)-\psi \bigr)_{H} + \nu^2 \bigl( \partial_{x}\theta(t), \partial_{x}(\theta(t)-\psi) \bigr)_{H} 
                    \\[1ex]
                    \displaystyle  +\int_{\Omega}\alpha(\eta(t))f_{\varepsilon}(\partial_{x}\theta(t))\, dx
\leq \int_{\Omega}\alpha(\eta(t))f_{\varepsilon}(\partial_{x}\psi)\, dx
                    \\[1ex]
                    \displaystyle +\bigl( M v(t), \theta(t)-\psi \bigr)_{H}, \mbox{ for any } \psi  \in V_0,  
                    \\[1ex]
                     \mbox{ a.e. } t \in (0, T), \mbox{ subject to $ \theta(0) = \theta_0 $ in $ H $.} 
                \end{array}
        \end{equation*}
\item[\textmd{(I-B)}]Let $\{\varepsilon_{n} \}_{n=1}^{\infty} \subset [0, 1]$, $\{[\bm{\eta}_{0, n}, \theta_{0, n}] \}_{n=1}^{\infty} = \{[\eta_{0, n}, \eta_{\Gamma, 0, n}, \theta_{0, n}] \}_{n=1}^{\infty}  \subset \mathbb{W} \times V_0$ with \linebreak $ \{ \bm{\eta}_{0, n} \}_{n = 1}^\infty$ $ = \{[\eta_{0, n}, \eta_{\Gamma, 0, n}] \}_{n=1}^{\infty}  $, and $\{[\bm{u}_n, v_n] \}_{n=1}^{\infty} = \{ [u_{n}, u_{\Gamma, n}, v_n]\}_{n=1}^{\infty} \subset \mathfrak{X} \times \mathscr{H} $ with  $ \{ \bm{u}_n \}_{n=1}^\infty$ $ = \{ [u_{n}, u_{\Gamma, n}]\}_{n=1}^{\infty} $, be given sequences such that:
\begin{equation}\label{w.i}
    \begin{array}{c}
        \varepsilon_{n} \to \varepsilon,\ [\bm{\eta}_{0, n}, \theta_{0, n}] \to [\bm{\eta}_0, \theta_0]  \mbox{ weakly\ in } \mathbb{W} \times V_0,
        \\[1ex]
        \mbox{ and } [L u_{n}, L_\Gamma u_{\Gamma, n}, M v_{n}] \to [L u, L_\Gamma u_\Gamma , M v] \ \mbox{weakly\ in}\ \mathfrak{X} \times \mathscr{H},  \mbox{ as } n \to \infty.
    \end{array}
\end{equation}
            Let $[\bm{\eta}, \theta] = [\eta, \eta_\Gamma, \theta] \in \mathfrak{X} \times \mathscr{H} $ with $ \bm{\eta} = [\eta, \eta_\Gamma] $ be the unique solution to (S)$_{\varepsilon}$, for the initial triplet $ [\bm{\eta}_0, \theta_0] = [\eta_0, \eta_{\Gamma, 0}, \theta_0] $ and forcing triplet $[\bm{u}, v] = [u, u_\Gamma, v] $. Additionally,  for any $n \in \mathbb{N}$, let $[\bm{\eta}_n, \theta_n] = [\eta_n, \eta_{\Gamma, n}, \theta_n] \in \mathfrak{X} \times \mathscr{H}$ with $\bm{\eta}_n = [\eta_n, \eta_{\Gamma, n}]$ be the unique solution to (S)$_{\varepsilon_n}$, for the initial triplet $ [\bm{\eta}_{0, n}, \theta_{0, n} ] = [\eta_{0, n}, \eta_{\Gamma, 0, n}, \theta_{0, n}]$ and forcing triplet $[\bm{u}_n, v_n] = [u_n, u_{\Gamma, n}, v_n]$. Then, it holds that:
\begin{align}\label{mThConv}
    [\bm{\eta}_{n}, \theta_{n}] \to &\, [\bm{\eta}, \theta] \mbox{ in } [C(\overline{Q}) \times C(\overline{\Sigma})] \times C(\overline{Q}) , \nonumber
    \\
     & \mbox{ in } \mathfrak{W} \times \mathscr{V}_0 , \mbox{ weakly in $ W^{1, 2}(0, T; \mathbb{X}) \times W^{1, 2}(0, T; H)  $,} \nonumber
    \\
    & \mbox{ and weakly-$*$ in $ L^\infty(0, T; \mathbb{W}) \times L^\infty(0, T; V_0) $, as $ n \to \infty $,}
\end{align}
and in particular,
        \begin{align}\label{mThConv00}
            \alpha''(\eta_{n}) & f_{\varepsilon_{n}}(\partial_x \theta_{n}) \to \alpha''(\eta)f_\varepsilon (\partial_x \theta) \mbox{ in $ \mathscr{H} $,}
            \nonumber
            \\
            & \mbox{ and weakly-$*$ in $ L^\infty(0, T; H) $, as $ n \to \infty $.}
        \end{align}
    \end{itemize}
\end{main}
\begin{remark}\label{Rem.mTh01Conv}
    As a consequence of \eqref{mThConv} and \eqref{mThConv00}, we further find a subsequence $ \{ n_i \}_{i = 1}^\infty \subset \{n\} $, such that:
    \begin{equation*}
        \begin{array}{l}
            [\eta_{n_i}, \theta_{n_i}] \to [\eta, \theta], ~ [\partial_x \eta_{n_i}, \partial_x \theta_{n_i}] \to [\partial_x \eta, \partial_x \theta], 
            \\[0.5ex]
            \qquad \qquad \mbox{and } \alpha''(\eta_{n_i}) f_{\varepsilon_{n_i}}(\partial_x \theta_{n_i}) \to \alpha''(\eta) f_\varepsilon(\partial_x \theta), 
            \\[0.5ex]
            \qquad \qquad \qquad \qquad \mbox{in the pointwise sense a.e. in $ Q $, as $ i \to \infty $,}
        \end{array}
    \end{equation*}
    and
    \begin{equation*}
        \begin{array}{l}
            [\eta_{n_i}(t), \theta_{n_i}(t)] \to [\eta(t), \theta(t)] \mbox{ in $ V \times V_0 $,}
            \\[0.5ex]
            \qquad \qquad \mbox{and } \alpha''(\eta_{n_i}(t)) f_{\varepsilon_{n_i}}(\partial_x \theta_{n_i}(t)) \to \alpha''(\eta(t)) f_\varepsilon(\partial_x \theta(t)) \mbox{ in $ H $,}
            \\[0.5ex]
            \qquad \qquad \qquad \qquad \mbox{in the pointwise sense for a.e. $ t \in (0, T) $, as $ i \to \infty $.}
        \end{array}
    \end{equation*}
\end{remark}
\begin{main}\label{mainTh02}
    Under the assumptions (\hyperlink{A0l}{A0})--(\hyperlink{A3l}{A3}), let us fix any constant $ \varepsilon \geq 0 $, and any initial triplet $ [\bm{\eta}_0, \theta_0] = [\eta_0, \eta_{\Gamma, 0}, \theta_0] \in \mathbb{W} \times V_0 $ with $ \bm{\eta}_0 = [\eta_0, \eta_{\Gamma, 0}] $. Then, the following two items hold.
    \begin{itemize}
        \item[\textmd{(II-A)}]The problem (OP)$_{\varepsilon}$ has at least one optimal control $[\bm{u}^*, v^*] = [u^{*}, u^*_\Gamma, v^*] \in \mathfrak{X} \times \mathscr{H} $ with $ \bm{u}^* = [u^{*}, u^*_\Gamma] $, so that:
\begin{equation*}
\mathcal{J}_{\varepsilon}(\bm{u}^{*}, v^*) = \mathcal{J}_\varepsilon(u^*, u^*_\Gamma, v^*) = \min_{[\bm{u}, v] \in \mathfrak{X} \times \mathscr{H}}\mathcal{J}_{\varepsilon}(\bm{u}, v) = \min_{[u, u_\Gamma, v] \in \mathfrak{X} \times \mathscr{H}}\mathcal{J}_{\varepsilon}(u, u_\Gamma, v).
\end{equation*}
        \item[\textmd{(II-B)}]Let $\{\varepsilon_{n} \}_{n=1}^{\infty} \subset [0, 1]$ and $\{[\bm{\eta}_{0, n}, \theta_{0, n}] \}_{n=1}^\infty = \{[\eta_{0, n}, \eta_{\Gamma, 0, n}, \theta_{0, n}] \}_{n=1}^{\infty}  \subset \mathbb{W} \times V_0 $ with  $ \{\bm{\eta}_{0, n}\}_{n = 1}^\infty = \{[\eta_{0, n}, \eta_{\Gamma, 0, n}] \}_{n=1}^{\infty}  $ be given sequences such that: 
\begin{align}\label{ass.4}
 \varepsilon_{n} \to \varepsilon,\ \mbox{ and }\ [\bm{\eta}_{0, n}, \theta_{0, n} ] \to [\bm{\eta}_0, \theta_0 ]\ {\rm weakly\ in}\ \mathbb{W} \times V_0, \ {\rm as}\ n \to \infty.
\end{align}
            In addition, for any $n \in \mathbb{N}$, let $[\bm{u}_n^*, v_n^*] = [u_{n}^*, u^*_{\Gamma, n}, v^*_n] \in \mathfrak{X} \times \mathscr{H} $ with  $ \bm{u}_n^* = [u_{n}^*, u^*_{\Gamma, n}]  $ be the optimal control of (OP)$_{\varepsilon_{n}}$. Then, there exist a subsequence $ \{ n_{i} \}_{i = 1}^{\infty} \subset \{ n \} $ and a triplet of functions $[\bm{u}^{**}, v^{**} ] = [u^{**}, u^{**}_\Gamma, v^{**}] \in \mathfrak{X} \times \mathscr{H} $ with $ \bm{u}^{**} =[u^{**}, u^{**}_\Gamma] $, such that: 
\begin{align*}
     \varepsilon_{n_{i}} \to \varepsilon,  \mbox{ and } &\, [L u_{n_{i}}^*, L_\Gamma u^*_{\Gamma, n_i}, M v_{n_{i}}^*] \to [L u^{**}, L_\Gamma u^{**}_\Gamma, M v^{**}]  
    \\
    &\, \mbox{weakly in } \mathfrak{X} \times \mathscr{H}, \mbox{ as } i \to \infty,
\end{align*}
and 
\begin{equation*}
[\bm{u}^{**}, v^{**} ] = [u^{**}, u^{**}_\Gamma, v^{**}] \ \mbox{is an optimal control of (OP)}_{\varepsilon}.
\end{equation*}
\end{itemize}
\end{main}

\begin{main}
    \label{mainTh03}
    Under the assumptions (\hyperlink{A0l}{A0})--(\hyperlink{A3l}{A3}), let us fix any initial triplet $ [\bm{\eta}_0, \theta_0] $ $= [\eta_0, \eta_{\Gamma, 0}, \theta_0] \in \mathbb{W} \times V_0 $ with $ \bm{\eta}_0 = [\eta_0, \eta_{\Gamma, 0}] $. Then, the following two items hold.
    \begin{itemize}
        \item[\textmd{\textit{(III-A)}}](Necessary condition for (OP)$_\varepsilon$ when $ \varepsilon > 0 $) 
            For any $ \varepsilon > 0 $, let $[\bm{u}_{\varepsilon}^{*}, v_{\varepsilon}^{*}] = [u_\varepsilon^*, u_{\Gamma, \varepsilon}^*, v_\varepsilon^*]$ $ \in \mathfrak{X} \times \mathscr{H} $ with  $ \bm{u}_\varepsilon^* = [u_\varepsilon^*, u_{\Gamma, \varepsilon}^*] $ be an optimal control of (OP)$_{\varepsilon}$, and let $ [\bm{\eta}_{\varepsilon}^{*}, \theta_{\varepsilon}^{*}] = [\eta_\varepsilon^*, \eta_{\Gamma, \varepsilon}^*, \theta_\varepsilon^*] \in \mathfrak{X} \times \mathscr{H} $  with  $ \bm{\eta}_\varepsilon^* = [\eta_\varepsilon^*, \eta_{\Gamma, \varepsilon}^*] $ be the solution to (S)$_{\varepsilon}$, for the initial triplet $[\bm{\eta}_0, \theta_0] = [\eta_0, \eta_{\Gamma, 0}, \theta_0]$ and forcing triplet $ [\bm{u}_{\varepsilon}^{*}, v_{\varepsilon}^{*}] = [u_\varepsilon^*, u_{\Gamma, \varepsilon}^*, v_\varepsilon^*]$. Then, it holds that:
    \begin{equation}\label{Thm.5-00}
        \displaystyle [L (u_\varepsilon^{*} + p_{\varepsilon}^{*}), L_\Gamma (u_{\Gamma, \varepsilon}^* +  p_{\Gamma, \varepsilon}^*) , M (v_\varepsilon^{*} + z_{\varepsilon}^{*})] = [0, 0, 0] 
    \mbox{ in $ \mathfrak{X} \times \mathscr{H} $,}
\end{equation}
            where $[\bm{p}_{\varepsilon}^{*}, z_{\varepsilon}^{*}] \in \mathfrak{Y} $ is a unique solution to the following variational system:
\begin{align}\label{Thm.5-01}
    & \bigl( -\partial_{t} \bm{p}_{\varepsilon}^{*}(t), \varphi \bigr)_{\mathbb{X}}
    + \bigl( \partial_{x} p_{\varepsilon}^{*}(t), \partial_{x}\varphi \bigr)_{H} + \bigl( [\alpha''(\eta_{\varepsilon}^{*})f_{\varepsilon}(\partial_{x}\theta_{\varepsilon}^{*})](t) p_{\varepsilon}^{*}(t), \varphi \bigr)_{H}
    \nonumber
    \\
    & \qquad  +\bigl( g'(\eta_{\varepsilon}^{*}(t)) p_{\varepsilon}^{*}(t), \varphi \bigr)_{H} +\bigl( [\alpha'(\eta_{\varepsilon}^{*}) f_{\varepsilon}'(\partial_{x}\theta_{\varepsilon}^{*})](t)\partial_{x}z_{\varepsilon}^{*}(t), \varphi \bigr)_{H}
    \nonumber
    \\
    & \qquad = \bigl( K (\eta_{\varepsilon}^{*}-\eta_{\mbox{\scriptsize ad}})(t), \varphi \bigr)_{H} + \bigl( K_\Gamma (\eta_{\Gamma, \varepsilon}^{*}-\eta_{\Gamma, \mbox{\scriptsize ad}})(t), \varphi_\Gamma \bigr)_{H_\Gamma}, 
    \\
    & \qquad\qquad \mbox{ for any $ \bm{\varphi} = [\varphi, \varphi_\Gamma] \in \mathbb{W} $,  and a.e. $ t \in (0, T) $;}
    \nonumber
\end{align}
and
\begin{align}\label{Thm.5-02}
    \bigl< -\partial_{t} & \bigl( \alpha_{0}z_{\varepsilon}^{*} \bigr)(t), \psi \bigr>_{V_{0}}  +\bigl( [\alpha(\eta_{\varepsilon}^{*})f_{\varepsilon}''(\partial_{x}\theta_{\varepsilon}^{*})](t)\partial_{x}z_{\varepsilon}^{*}(t) + \nu^2 \partial_{x}z_{\varepsilon}^{*}(t), \partial_x \psi \bigl)_H
    \nonumber
    \\
    & +\bigl( [\alpha'(\eta_{\varepsilon}^{*})f_{\varepsilon}'(\partial_{x}\theta_{\varepsilon}^{*})](t)p_{\varepsilon}^{*}(t), \partial_{x}\psi \bigr)_{H} = \bigl( \Lambda (\theta_{\varepsilon}^{*}-\theta_{\mbox{\scriptsize ad}})(t), \psi \bigr)_{H},
    \\
    & \qquad \qquad \mbox{for any}\ \psi \in V_{0},\ \mbox{and a.e.}\ t \in (0, T);
    \nonumber
\end{align}
subject to the terminal condition:
\begin{equation}\label{Thm.5-03}
[\bm{p}_{\varepsilon}^{*}(T), z_{\varepsilon}^{*}(T)] = [p_\varepsilon^*(T), p_{\Gamma, \varepsilon}^*(T), z_\varepsilon^*(T)] = [0, 0, 0]\ \mbox{in}\ \mathbb{X} \times H.
\end{equation}
    \item[\textmd{\textit{(III-B)}}]Let us define a Hilbert space $ \mathscr{U}_0 $ as follows:
        \begin{equation*}
            \mathscr{U}_0 := \left\{ \begin{array}{l|l}
                \psi \in W^{1, 2}(0, T; H) \cap \mathscr{V}_0 & \psi(0) = 0 \mbox{ in $ H $}
            \end{array} \right\}.
        \end{equation*}
            Then, there exists an optimal control $ [\bm{u}^{\circ}, v^{\circ}] = [u^\circ, u^\circ_\Gamma, v^\circ] \in \mathfrak{X} \times \mathscr{H} $ with $ \bm{u}^\circ = [u^\circ, u^\circ_\Gamma] $ of the problem (OP)$_0$, together with the solution $ [\bm{\eta}^\circ, \theta^\circ] = [\eta^\circ, \eta^\circ_\Gamma, \theta^\circ] \in \mathfrak{X} \times \mathscr{H} $ to the system (S)$_0$ with $ \bm{\eta}^\circ = [\eta^\circ, \eta^\circ_\Gamma] $, for the initial triplet $[\bm{\eta}_0, \theta_0] = [\eta_0, \eta_{\Gamma, 0}, \theta_0]$ and forcing triplet $ [\bm{u}^\circ, v^\circ] = [u^\circ, u_\Gamma^\circ, v^\circ]$, and moreover, there exist a triplet of functions $ [\bm{p}^{\circ}, z^{\circ}] = [p^\circ, p^\circ_\Gamma, z^\circ] \in \mathfrak{X} \times \mathscr{H} $ with $ \bm{p}^\circ = [p^\circ, p^\circ_\Gamma] $, a pair of functions $ [\xi^\circ, \nu^\circ] \in \mathscr{H} \times L^\infty(Q) $, and a distribution $ \zeta^\circ \in \mathscr{U}_0^* $, such that:
    \begin{equation}\label{Thm.5-10}
        [L (u^\circ + p^\circ), L_\Gamma(u^\circ_\Gamma + p^\circ_\Gamma), M (v^\circ + z^\circ)] = [0, 0, 0] \mbox{ in $ \mathfrak{X} \times \mathscr{H}$;}
    \end{equation}
    \begin{equation}\label{Thm.5-13}
        \begin{cases}
            \bm{p}^\circ = [p^\circ, p_\Gamma^\circ] \in W^{1, 2}(0, T; \mathbb{X}) \cap L^\infty(0, T; \mathbb{W}) \subset C(\overline{Q}) \times C(\overline{\Sigma}),
            \\
            z^\circ \in L^\infty(0, T; H) \cap \mathscr{V}_0,
            \\
            \nu^\circ \in \Sgn^1(\partial_x \theta^\circ), \mbox{ a.e. in $ Q $;}
        \end{cases}
    \end{equation}
    \begin{align}\label{Thm.5-11}
        & \bigl( -\partial_t \bm{p}^\circ, \bm{\varphi} \bigr)_{\mathfrak{X}} 
        +\bigl( \partial_x p^\circ, \partial_x \varphi \bigr)_{\mathscr{H}} 
        +\bigl( \alpha''(\eta^\circ)|\partial_x \theta^\circ| p^\circ, \varphi \bigl)_\mathscr{H} 
        \nonumber
        \\
        & \qquad  +\bigl( g'(\eta^\circ)p^\circ, \varphi \bigr)_{\mathscr{H}} +\bigl(\alpha'(\eta^\circ) \xi^\circ, \varphi \bigr)_{\mathscr{H}} 
        \nonumber
        \\
        & = \bigl( K (\eta^\circ -\eta_\mathit{ad}), \varphi \bigr)_{\mathscr{H}} + \bigl( K_\Gamma (\eta^\circ_\Gamma -\eta_{\Gamma, \mathit{ad}}), \varphi_\Gamma \bigr)_{\mathscr{H}_\Gamma}, 
        \\
        & 
        \qquad \mbox{for any $ \bm{\varphi} = [\varphi, \varphi_\Gamma] \in \mathfrak{W} $,} 
        \nonumber
        \\
        & \mbox{subject to $ \bm{p}^\circ(T) = [p^\circ(T), p^\circ_\Gamma(T)] = [0, 0] $ in $ \mathbb{X} $;}
        \nonumber
    \end{align}
    and
    \begin{align}\label{Thm.5-12}
        \bigl( \alpha_0 z^\circ & , \partial_t \psi \bigr)_{\mathscr{H}} +\bigl< \zeta^\circ, \psi \bigr>_{\mathscr{U}_0} +\bigl( \nu^2 \partial_x z^\circ +\alpha'(\eta^\circ) \nu^\circ p^\circ, \partial_x \psi \bigr)_{\mathscr{H}}
        \nonumber
        \\
        & = \bigl( \Lambda (\theta^\circ -\theta_\mathit{ad}), \psi \bigr)_{\mathscr{H}}, \mbox{ for any $ \psi \in \mathscr{U}_0 $.}
    \end{align}
\end{itemize}
\end{main}
\begin{remark}\label{Rem.mTh03}
    Let $ \mathcal{R}_T \in \mathscr{L}(\mathscr{H}) $ be an isomorphism, defined as:
    \begin{equation*}
        \bigl( \mathcal{R}_T \varphi \bigr)(t) := \varphi(T -t) \mbox{ in $ H $, for a.e. $ t \in (0, T) $.}
    \end{equation*}
    Also, let us fix $ \varepsilon > 0 $, and denote by $ \mathcal{Q}_\varepsilon^* \in \mathscr{L}(\mathfrak{X} \times \mathscr{H}; \mathfrak{Y}) $ the restriction $\mathcal{P}|_{\{[0, 0, 0]\} \times [\mathfrak{X} \times \mathscr{H}]} $ of the bounded linear operator $ \mathcal{P} = \mathcal{P}(a, b, \mu, \omega, A) : [\mathbb{W} \times H] \times [\mathfrak{X} \times \mathscr{V}_0^*] \longrightarrow \mathfrak{Y} $, as in Remark \ref{appendix}, in the case when:
    \begin{equation}\label{setRem4}
        \begin{cases}
            [a, b] = \mathcal{R}_T [\alpha_0, -\partial_t \alpha_0] \mbox{ in $ W^{1, \infty}(Q) \times L^\infty(Q) $,}
            \\
            \mu = \mathcal{R}_T \bigl[ g'(\eta_\varepsilon^*) + \alpha''(\eta_\varepsilon^*) f_\varepsilon(\partial_x \theta_\varepsilon^*) \bigr] \mbox{ in $ L^\infty(0, T; H) $,}
            \\
            [\omega, A] = \mathcal{R}_T \bigl[ \alpha'(\eta_\varepsilon^*) f_\varepsilon'(\partial_x \theta_\varepsilon^*), \alpha(\eta_\varepsilon^*)f_\varepsilon''(\partial_x \theta_\varepsilon^*) \bigl] \mbox{ in $ [L^\infty(Q)]^2 $.}
        \end{cases}
    \end{equation}
    On this basis, let us define:
    \begin{equation*}
        \mathcal{P}_\varepsilon^* := \mathcal{R}_T \circ \mathcal{Q}_\varepsilon^* \circ \mathcal{R}_T \mbox{ in $ \mathscr{L}(\mathfrak{X} \times \mathscr{H}; \mathfrak{Y}) $.}
    \end{equation*}
    Then, having in mind:
    \begin{equation}\label{productrule}
        \partial_t (\alpha_0 \tilde{z}) = \alpha_0 \partial_t \tilde{z} +\tilde{z} \partial_t \alpha_0 \mbox{ in $ \mathscr{V}_0^* $, for any $ \tilde{z} \in W^{1, 2}(0, T; V^*_0) $,}
    \end{equation}
    we can obtain the unique solution $[\bm{p}_{\varepsilon}^{*}, z_{\varepsilon}^{*}] = [p_\varepsilon^*, p_{\Gamma, \varepsilon}^*, z_\varepsilon^*] \in \mathfrak{Y} $ with $ \bm{p}_\varepsilon^* = [p_\varepsilon^*, p_{\Gamma, \varepsilon}^*] 
    $ to the variational system \eqref{Thm.5-01}--\eqref{Thm.5-03} as follows:
\begin{equation*}
    [\bm{p}_\varepsilon^*, z_\varepsilon^*] =[p_\varepsilon^*, p_{\Gamma, \varepsilon}^*, z_\varepsilon^*] = \mathcal{P}_\varepsilon^* \bigl[ K  (\eta_\varepsilon^* -\eta_\mathrm{ad}), K_\Gamma(\eta_{\Gamma, \varepsilon}^* -\eta_{\Gamma, \mathrm{ad}}), \Lambda (\theta_\varepsilon^* -\theta_\mathrm{ad}) \bigr] \mbox{ in $ \mathfrak{Y} $.}
\end{equation*}
\end{remark}

\section{Proof of Main Theorem \ref{mainTh01}} 

In this Section, we give the proof of the first Main Theorem \ref{mainTh01}. Before the proof, we refer to the reformulation method as in \cite{MR3888636}, and reduce the state-system (S)$_{\varepsilon}$ to an evolution equation in the Hilbert space $\mathbb{X} \times H$.

Let us fix any $\varepsilon \geq 0$. Besides, for any $ R \geq 0 $, let us define a proper functional $\Phi_{\varepsilon}^R: \mathbb{X} \times H \longrightarrow [0, \infty]$, by setting: 
\begin{align}\label{Phi_eps}
    \Phi_{\varepsilon}^R: w & = [\bm{\eta}, \theta]  = [\eta, \eta_\Gamma, \theta] \in \mathbb{X} \times H \mapsto \Phi_{\varepsilon}^R(w) =  \Phi_{\varepsilon}^R(\bm{\eta}, \theta) = \Phi_{\varepsilon}^R(\eta, \eta_\Gamma, \theta)
    \nonumber
    \\
    &:= \left\{
        \begin{array}{ll}
            \multicolumn{2}{l}{\displaystyle\frac{1}{2}\int_{\Omega}|\partial_x \eta|^{2}\, dx +\frac{R}{2} \int_\Omega |\eta|^2 \, dx +\frac{1}{2}\int_{\Omega}\left(\nu f_{\varepsilon}(\partial_x \theta) + \frac{1}{\nu}\alpha(\eta) \right)^{2}\, dx,}
            \\[2ex]
            & \mbox{if $[\bm{\eta}, \theta] = [\eta, \eta_\Gamma, \theta] \in \mathbb{W} \times V_{0}$ with $\bm{\eta} = [\eta, \eta_\Gamma] $,}
            \\[2ex]
            \infty, & \mbox{otherwise.}
        \end{array}
    \right. 
\end{align}
Note that the assumptions (\hyperlink{A1l}{A1}) and (\hyperlink{A2l}{A2}) guarantee the lower semi-continuity and convexity of $\Phi_{\varepsilon}^R$ on $\mathbb{X} \times H$. 

\begin{remark}\label{Rem.mTh01-01}
    As consequences of standard variational methods, we easily check the following facts.
    \begin{description}
        \item[(\hypertarget{Fact3}{Fact\,3})]For the operator $ \partial_{\bm{\eta}} \Phi_\varepsilon^R : \mathbb{X} \times H \longrightarrow 2^{\mathbb{X}} $, 
            \begin{equation*}
                D(\partial_{\bm{\eta}} \Phi_\varepsilon^R) = \left\{ \begin{array}{l|l}
                    [\tilde{\bm{\eta}}, \tilde{\theta}] = [\tilde{\eta}, \tilde{\eta}_\Gamma, \tilde{\theta}] \in \mathbb{W} \times V_0 & \parbox{6.0cm}{$ \tilde{\eta} \in H^2(\Omega) $ with $\partial_x \tilde{\eta}\cdot n_\Gamma = 0$ on $\Gamma$}
                \end{array} \right\},
            \end{equation*}
            and  $\partial_{\bm{\eta}}\Phi_\varepsilon^R$ is a single-valued operator such that:
            \begin{equation*}
            \partial_{\bm{\eta}} \Phi_\varepsilon^R(w) = \rule{0pt}{18pt}^\mathrm{t} \hspace{-0.5ex} \begin{bmatrix}
                    -{\partial_x^2}\eta + R\eta + \alpha'(\eta)f_{\varepsilon}({\partial_x}\theta) + \nu^{-2}\alpha(\eta)\alpha'(\eta)  
                    \\
                    0
                \end{bmatrix} 
                \\[2ex] 
                 \mbox{ in } \mathbb{X} ,
            \end{equation*} 
            for any $ w = [\bm{\eta}, \theta] = [\eta, \eta_\Gamma, \theta] \in D(\partial_{\bm{\eta}} \Phi_\varepsilon^R) $.           
        \item[(\hyperlink{Fact4}{Fact\,4})]$\theta \in D(\partial_\theta \Phi_\varepsilon^R)$, and $ \theta^* \in \partial_\theta \Phi_\varepsilon^R(w) = \partial_\theta \Phi_\varepsilon^R(\bm{\eta}, \theta) = \partial_\theta \Phi_\varepsilon^R(\eta, \eta_\Gamma, \theta) $, iff. $\theta \in V_0$, and
            \begin{equation*}
                \begin{array}{c}
                    \displaystyle (\theta^*, \theta -\psi)_H \geq \nu^2 (\partial_x \theta, \partial_x(\theta -\psi))_{H} +\int_\Omega \alpha(\eta)f_\varepsilon(\partial_x \theta) \, dx -\int_\Omega \alpha(\eta) f_\varepsilon(\partial_x \psi), 
                    \\[2ex]
                    \mbox{for any $ \psi \in V_0 $.} 
                \end{array}
            \end{equation*}
    \end{description}
    In addition, let us define time-dependent operators $\mathcal{A}(t) \in \mathscr{L}(\mathbb{X} \times H)$, for $t \in [0, T]$, nonlinear operators $\mathcal{G}^R : \mathbb{X} \times H \longrightarrow \mathbb{X} \times H$, for $ R \geq 0 $, by setting:
\begin{align}\label{AA_0}
    \mathcal{A}(t): {w} & \, = [\bm{\eta}, {\theta}] = [\eta, \eta_\Gamma, \theta] \in \mathbb{X} \times H \nonumber 
    \\
     & \, \mapsto \mathcal{A}(t) {w} := [{\eta}, \, {\eta}_\Gamma, \, \alpha_0(t) {\theta}] \in \mathbb{X} \times H, \mbox{ for }t \in [0, T],
\end{align}
and
\begin{align}\label{GG}
    \mathcal{G}^R: w & \, = [\bm{\eta}, \theta] = [\eta, \eta_\Gamma, \theta] \in \mathbb{X} \times H \nonumber 
    \\
    & \, \mapsto \mathcal{G}^R(w) := \bigl[ g(\eta) -R \eta -\nu^{-2}\alpha(\eta)\alpha'(\eta), ~ 0 , ~ 0 \bigr] \in \mathbb{X} \times H,
\end{align}
    respectively. Then, based on the above (\hyperlink{Fact3}{Fact\,3}) and (\hyperlink{Fact4}{Fact\,4}), it is verified that the state-system (S)$_\varepsilon$ is equivalent to the following Cauchy problem:
    \begin{equation*}
        \left\{ \parbox{11cm}{
            $ \mathcal{A}(t) w'(t) +\bigl[ \partial_{\bm{\eta}} \Phi_\varepsilon^R \times \partial_\theta \Phi_\varepsilon^R \bigr](w(t)) +\mathcal{G}^R(w(t)) \ni \mathfrak{f}(t) $  in $ \mathbb{X} \times H $, 
            \\[1ex]
            \hspace*{4ex}a.e. $ t \in (0, T) $,
            \\[1.5ex]
            $ w(0) = w_0 $ in $ \mathbb{X} \times H $.
        } \right.
    \end{equation*}
    In the context, ``\,$'$\,'' is the time-derivative, and
    \begin{equation}\label{w0-f}
        \left\{ \hspace{-3ex} \parbox{11cm}{
            \vspace{-1ex}
            \begin{itemize}
                \item $ w_0 := [\bm{\eta}_0, \theta_0] = [\eta_0, \eta_{\Gamma, 0}, \theta_0] \in \mathbb{W} \times V_0 $ with $\bm{\eta}_0 = [\eta_0, \eta_{\Gamma, 0}] $ is the initial data of $ w = [\bm{\eta}, \theta] = [\eta, \eta_\Gamma, \theta] $, 
                \item $ \mathfrak{f} := [L u, L_\Gamma u_\Gamma, M v] \in \mathfrak{X} \times \mathscr{H} $ is the forcing term of the Cauchy problem.  
            \vspace{-1ex}
            \end{itemize}
        }
        \right.
    \end{equation}
\end{remark}

Now, before the proof of Main Theorem \ref{mainTh01}, we prepare the following Key-Lemma and Corollary.
\begin{keyLem}\label{Lem03-01}
Let us assume (\hyperlink{A0l}{A0})--(\hyperlink{A3l}{A3}), and let us fix any $\varepsilon \geq 0$. Then, there exists a positive constant $ R_0 > 0 $ such that: 
\begin{align*}
    \partial \Phi_\varepsilon^{R_0} = \bigl[ \partial_{\bm{\eta}} \Phi_\varepsilon^{R_0} \times \partial_\theta \Phi_\varepsilon^{R_0} \bigr] \mbox{ in $ [\mathbb{X} \times H] \times [\mathbb{X} \times H] $.}
\end{align*}
\end{keyLem}
\paragraph{Proof.}{
    We set:
    \begin{equation}\label{R0}
        R_0 := 1 +\frac{2}{\nu^{2}}|\alpha|^2_{L^\infty(\mathbb{R})},
    \end{equation}
    and prove this $ R_0 $ is the required constant. 

    In the light of \eqref{prodSubDif}, it is immediately verified that:
    \begin{equation*}
        \partial \Phi_\varepsilon^{R_0}  \subset \bigl[ \partial_{\bm{\eta}} \Phi_\varepsilon^{R_0} \times \partial_\theta \Phi_\varepsilon^{R_0} \bigr] \mbox{ in $ [\mathbb{X} \times H] \times [\mathbb{X} \times H] $.}
    \end{equation*}
    Hence, with the maximality of the monotone graph $ \partial \Phi_\varepsilon^{R_0} $ in $ [\mathbb{X} \times H] \times [\mathbb{X} \times H] $ in mind, we can reduce our task to show the monotonicity of $ \bigl[ \partial_{\bm{\eta}} \Phi_\varepsilon^{R_0} \times \partial_\theta \Phi_\varepsilon^{R_0} \bigr] $ in $ [\mathbb{X} \times H] \times [\mathbb{X} \times H] $. 
     
Let us assume:
\begin{equation*}
    \begin{array}{c}
        \displaystyle [w, w^*] \in \bigl[ \partial_{\bm{\eta}} \Phi_\varepsilon^{R_0} \times \partial_\theta \Phi_\varepsilon^{R_0} \bigr] \mbox{ and } [\tilde{w}, \tilde{w}^*] \in \bigl[ \partial_{\bm{\eta}} \Phi_\varepsilon^{R_0} \times \partial_\theta \Phi_\varepsilon^{R_0} \bigr] \mbox{ in $ [\mathbb{X} \times H] \times [\mathbb{X} \times H] $.}
    \end{array}
\end{equation*}
Then, by using \eqref{AA_0}, \eqref{GG}, (\hyperlink{Fact3}{Fact\,3}), (\hyperlink{Fact4}{Fact\,4}), and Young's inequality, we compute that:
\begin{subequations}\label{IA0-00}
\begin{align}
    & (w^* -\tilde{w}^*, w -\tilde{w})_{\mathbb{X} \times H} \geq I_A^1 + I_A^2 + I_A^3, 
\end{align}
with 
\begin{align}
I_A^1 := & |\partial_x(\eta - \tilde{\eta})|_{H}^2 + R_0|\eta - \tilde{\eta}|_H^2 +\nu^2|\partial_x (\theta - \tilde{\theta})|_{H}^2,
\end{align}
\begin{align}
    I_A^2 ~& := (\alpha'(\eta)f_\varepsilon(\partial_x \theta) - \alpha'(\tilde{\eta})f_\varepsilon(\partial_x \tilde{\theta}), \eta-\tilde{\eta})_H \nonumber \\
  &   = \int_{\Omega}f_\varepsilon(\partial_x \theta)(\alpha'(\eta) - \alpha'(\tilde{\eta}))(\eta - \tilde{\eta}))\, dx 
    \nonumber
    \\
    & \qquad + \int_{\Omega}\alpha'(\tilde{\eta})(f_\varepsilon(\partial_x \theta)- f_\varepsilon(\partial_x \tilde{\theta}))(\eta - \tilde{\eta})\, dx \nonumber \\
  &   \geq - |\alpha'|_{L^\infty(\mathbb{R})}|\eta-\tilde{\eta}|_H|\partial_x (\theta - \tilde{\theta})|_{H} \nonumber \\
  &   \geq - \frac{|\alpha'|_{L^\infty(\mathbb{R})}^2}{\nu^2}|\eta-\tilde{\eta}|_H^2 - \frac{\nu^2}{4}|\partial_x (\theta - \tilde{\theta})|_{H}^2,
\end{align}
and
\begin{align}
 I_A^3 & := \int_{\Omega} (\alpha(\eta) - \alpha(\tilde{\eta}))(f_\varepsilon(\partial_x \theta) - f_\varepsilon(\partial_x \tilde{\theta}))\, dx \nonumber \\ 
  &   \geq  - |\alpha'|_{L^\infty(\mathbb{R})}|\eta-\tilde{\eta}|_H|\partial_x (\theta - \tilde{\theta})|_{H} \nonumber \\
  &   \geq  - \frac{|\alpha'|_{L^\infty(\mathbb{R})}^2}{\nu^2}|\eta-\tilde{\eta}|_H^2 - \frac{\nu^2}{4}|\partial_x (\theta - \tilde{\theta})|_{H}^2.
\end{align}
\end{subequations}
Due to \eqref{R0}, the inequalities in \eqref{IA0-00} lead to:
\begin{align*}
 & (w^* -\tilde{w}^*, w -\tilde{w})_{\mathbb{X} \times H} \geq |\eta-\tilde{\eta}|_V^2 + \frac{\nu^2}{2}|\theta-\tilde{\theta}|_{V_0}^2 \geq 0,
\end{align*}
which implies the monotonicity of the operator $\bigl[ \partial_{\bm{\eta}} \Phi_\varepsilon^{R_0} \times \partial_\theta \Phi_\varepsilon^{R_0} \bigr]$ in $ [\mathbb{X} \times H] \times [\mathbb{X} \times H] $. 
\qed
}

\begin{corollary}\label{03-01}
Under the notations and assumptions as in the previous Key-Lemma \ref{Lem03-01}, it holds that
\begin{align*}
    \partial \Phi_\varepsilon^{R} = \bigl[ \partial_{\bm{\eta}} \Phi_\varepsilon^{R} \times \partial_\theta \Phi_\varepsilon^{R} \bigr] \mbox{ in $ [\mathbb{X} \times H] \times [\mathbb{X} \times H] $, for any $R \geq 0$.}
\end{align*}
\end{corollary}
\paragraph{Proof.}{
Let us take arbitrary two constants $0 \leq R, \tilde{R} < \infty$.
Then from (\hyperlink{Fact3}{Fact\,3}), we immediately have 
\begin{subequations}\label{sharp1}
\begin{align}
 D(\partial_{\bm{\eta}} \Phi_{\varepsilon}^{R}) = D(\partial_{\bm{\eta}} \Phi_{\varepsilon}^{\tilde{R}}) \mbox{ in } \mathbb{W},
\end{align}
and
\begin{align}
    \partial_{\bm{\eta}} \Phi_{\varepsilon}^{R}(w) ~& = \rule{0pt}{18pt}^\mathrm{t} \hspace{-0.5ex} \begin{bmatrix}
                    -{\partial_x^2}\eta + \tilde{R}\eta + (R - \tilde{R})\eta+ \alpha'(\eta)f_{\varepsilon}({\partial_x}\theta) + \nu^{-2}\alpha(\eta)\alpha'(\eta)  
                    \\
                    0
                \end{bmatrix} \nonumber
                \\[1ex]
                 & \hspace{8ex}= \partial_{\bm{\eta}} \Phi_{\varepsilon}^{\tilde{R}}(w) + (R -\tilde{R})[\eta, 0] \mbox{ in } \mathbb{X}, 
                 \\
                 & \mbox{for any } w = [\bm{\eta}, \theta] = [\eta, \eta_\Gamma, \theta] \in D(\partial_{\bm{\eta}} \Phi_{\varepsilon}^{R}) = D(\partial_{\bm{\eta}} \Phi_{\varepsilon}^{\tilde{R}}). \nonumber
\end{align}
\end{subequations}
Also, as a straightforward consequence of (\hyperlink{Fact4}{Fact\,4}), it is seen that:
\begin{align}\label{sharp2}
\partial_{\theta} \Phi_{\varepsilon}^{R} = \partial_{\theta} \Phi_{\varepsilon}^{\tilde{R}} \mbox{ in } H \times H.
\end{align}
In the meantime, invoking \eqref{Phi_eps}, \cite[Theorem 2.10]{MR2582280}, and \cite[Corollary 2.11]{MR0348562}, we will infer that
\begin{subequations}\label{sharp3}
\begin{align}
D(\partial \Phi_{\varepsilon}^{R}) = D(\partial \Phi_{\varepsilon}^{\tilde{R}}) \mbox{ in } \mathbb{W} \times V_0,
\end{align}
and
\begin{align}
\partial \Phi_{\varepsilon}^{R}(w) = \partial \Phi_{\varepsilon}^{\tilde{R}}(w) + (R - \tilde{R})[\eta, 0, 0] \mbox{ in } \mathbb{X} \times H.
\end{align}
\end{subequations}

Now, let us take the constant $R_0 > 0$ obtained in Key-Lemma \ref{Lem03-01}. 
Then, owing to \eqref{sharp1}--\eqref{sharp3}, and Key-Lemma \ref{Lem03-01}, we can compute that
\begin{align}\label{sharp4}
    \bigl[\partial_{\bm{\eta}} \Phi_{\varepsilon}^{R} & \times \partial_{\theta} \Phi_{\varepsilon}^{R} \bigr](w) = \bigl[\partial_{\bm{\eta}} \Phi_{\varepsilon}^{R_0} \times \partial_{\theta} \Phi_{\varepsilon}^{R_0} \bigr](w) + (R - R_0)[\eta, 0, 0]
\nonumber
\\
    & \hspace{-1ex} = \partial \Phi_{\varepsilon}^{R_0}(w) + (R -R_0)[\eta, 0, 0] = \partial \Phi_{\varepsilon}^{R}(w) \mbox{ in } \mathbb{X} \times H,
 \\
    & \hspace{-2ex}  \mbox{ for any } w \in D(\partial_{\bm{\eta}} \Phi_{~\varepsilon}^{R} \times \partial_{\theta} \Phi_{\varepsilon}^{R}) = D(\partial_{\bm{\eta}} \Phi_{\varepsilon}^{R}) \cap D(\partial_{\theta} \Phi_{\varepsilon}^{R}).\nonumber
\end{align}
In the light of \eqref{prodSubDif}, the above \eqref{sharp4} is sufficient to conclude this Corollary.
\qed
\medskip

}

\begin{remark}\label{Rem.sols}
Let $\varepsilon \geq 0$ be arbitrary constant.
Then, as a consequence of (\hyperlink{Fact3}{Fact\,3}), (\hyperlink{Fact4}{Fact\,4}), Key-Lemma \ref{Lem03-01}, and Corollary \ref{03-01}, we can say that the state-system (S)$_\varepsilon$ is equivalent to the following Cauchy problem of evolution equation, denoted by (E)$_\varepsilon$.
\begin{description}
 \item[\textmd{(E)$ _\varepsilon $}]:
\vspace{-4ex}
\end{description}
\begin{align*}
      &  \left\{ \parbox{13cm}{
            $ \mathcal{A}(t) w'(t) + \partial \Phi_\varepsilon^R (w(t)) +\mathcal{G}^R(w(t)) \ni \mathfrak{f}(t) $ in $ \mathbb{X} \times H $, a.e. $ t \in (0, T) $,
            \\[1.5ex]
            $ w(0) = w_0 $ in $ \mathbb{X} \times H $,
        } \right.
    \end{align*}
for any $ R \geq 0 $. 
\end{remark}

Now, we are ready to prove the Main Theorem \ref{mainTh01}.
\paragraph{\textbf{Proof of Main Theorem \ref{mainTh01} (I-A).}}{
    Let us fix any $R > 0$. Then, under the setting \eqref{Phi_eps}--\eqref{w0-f}, we immediately check that:
    \begin{itemize}
        \item[\textmd{(ev.0)}]for any $ t \in [0, T] $, $ \mathcal{A}(t) \in \mathscr{L}(\mathbb{X} \times H) $ is positive and selfadjoint, and  
            \begin{equation*}
                (\mathcal{A}(t) w, w)_{\mathbb{X} \times H} \geq \delta_* |w|_{\mathbb{X} \times H}^2, \mbox{ for any $ w \in \mathbb{X} \times H $,}
            \end{equation*}
            with the constant $ \delta_* \in (0, 1) $ as in (\hyperlink{A1l}{A1});
        \item[\textmd{(ev.1)}]$ \mathcal{A} \in W^{1, \infty}(0, T; \mathscr{L}(\mathbb{X} \times H)) $, and  
            \begin{equation*}
                A^* := \mathrm{ess} \sup_{\hspace{-3ex}t \in (0, T)} \left\{ \max \{ |\mathcal{A}(t)|_{\mathscr{L}(\mathbb{X} \times H)}, |\mathcal{A}'(t)|_{\mathscr{L}(\mathbb{X} \times H)} \} \right\} \leq 1 +|\alpha_0|_{W^{1, \infty}(Q)} < \infty;
            \end{equation*}
        \item[\textmd{(ev.2)}]$ \mathcal{G}^R : \mathbb{X} \times H \longrightarrow \mathbb{X} \times H $ is a Lipschitz continuous operator with a Lipschitz constant:
            \begin{equation*}
                L_* := R +|g'|_{L^\infty(\mathbb{R})} +\nu^{-2} {\textstyle{\bigl| (\alpha \alpha')' \bigr|_{L^\infty(\mathbb{R})}},}
            \end{equation*}
            and $ \mathcal{G}^R $ has a $ C^1 $-potential functional
            \begin{align*}
                \widehat{\mathcal{G}}^R: w ~& =  [\bm{\eta}, \theta ] = [\eta, \eta_\Gamma, \theta] \in \mathbb{X} \times H 
                \\
                & \mapsto \widehat{\mathcal{G}}^R(w) := \int_\Omega \left( G(\eta) - \frac{R\eta^2}{2} -\frac{\alpha(\eta)^2}{2 \nu^2} \right) \, dx \in \mathbb{R};
            \end{align*}
        \item[\textmd{(ev.3)}]$ \Phi_\varepsilon^R \geq 0 $ on $ \mathbb{X} \times H $, and the sublevel set $ \bigl\{  \tilde{w} \in \mathbb{X} \times H \, \bigl| \, \Phi_\varepsilon^R(\tilde{w}) \leq r \bigr\} $ is contained in a compact set $ K_\nu^R(r) $ in $ \mathbb{X} \times H $, defined as
            \begin{equation*}
               K_\nu^R(r) := \left\{ \begin{array}{l|l}
                    \tilde{w} = [\tilde{\bm{\eta}}, \tilde{\theta} ] = [\tilde{\eta}, \tilde{\eta}_\Gamma, \tilde{\theta}] \in \mathbb{W} \times V_0 &  |\tilde{\bm{\eta}}|_\mathbb{W}^2 + |\tilde{\theta}|_{V_0}^2 \leq \frac{2r}{\min\{1, R, \nu^2 \}}
            \end{array} \right\}, 
        \end{equation*}
        for any $ r \geq 0 $.
    \end{itemize}
        On account of \eqref{Phi_eps}--\eqref{w0-f} and (ev.0)--(ev.3), we can apply Proposition \ref{Lem.CP}, as the case when:
    \begin{align*}
        & X = \mathbb{X} \times H,  ~ \mathcal{A}_0 =  \mathcal{A} \mbox{ in $ W^{1, \infty}(0, T; \mathscr{L}(\mathbb{X} \times H)) $, }\\
         \mathcal{G}_0 =  \mathcal{G}^R & \mbox{ on $ \mathbb{X} \times H $,  } \Psi_0 = \Phi_\varepsilon^R \mbox{ on $ \mathbb{X} \times H $, and } \mathfrak{f}_0 = \mathfrak{f} \mbox{ in } \mathfrak{X} \times \mathscr{H},
    \end{align*}
        and we can find a solution $ w = [\bm{\eta}, \theta] = [\eta, \eta_\Gamma, \theta] \in \mathfrak{X} \times \mathscr{H} $ with $\bm{\eta} = [\eta, \eta_\Gamma]$ to the Cauchy problem (E)$_\varepsilon$. In the light of  Proposition \ref{Lem.CP} and Remark \ref{Rem.sols}, finding this $ w = [\bm{\eta}, \theta] = [\eta, \eta_\Gamma, \theta]$ directly leads to the existence and uniqueness of solution to the state-system (S)$_\varepsilon$.
    \medskip    
    \qed
}
\medskip

\paragraph{\textbf{Proof of Main Theorem \ref{mainTh01} (I-B).}}{
    Under the assumptions and notations as in Theorem 1 (I-A), we first fix a constant $R > 0$, and invoke Remark \ref{Rem.sols} to confirm that the solution $w := [\bm{\eta}, \theta] = [\eta, \eta_\Gamma, \theta] \in \mathfrak{X} \times \mathscr{H}$ with $\bm{\eta} = [\eta, \eta_\Gamma]$ to (S)$_\varepsilon$ coincides with the solution to the Cauchy problem (E)$_\varepsilon$, and as well as, the solutions $w_n:= [\bm{\eta}_n, \theta_n] = [\eta_n, \eta_{\Gamma, n}, \theta_n] \in \mathfrak{X} \times \mathscr{H}$ with $\bm{\eta}_n = [\eta_n, \eta_{\Gamma, n}]$ to (S)$_{\varepsilon_n}$, $n = 1, 2, 3, \ldots ,$ coincide with the solutions to the Cauchy problems (E)$_{\varepsilon_n}$ for the initial data $ w_{0, n}:= [\bm{\eta}_{0, n} \theta_{0, n}] = [\eta_{0, n}, \eta_{\Gamma, 0, n}, \theta_{0, n}] \in \mathbb{W} \times V_0 $ with $\bm{\eta}_{0, n} = [\eta_{0, n}, \eta_{\Gamma, 0, n}]$, and forcing terms $ \mathfrak{f_n} = [L u_n, L_\Gamma u_{\Gamma, n}, M v_n] \in \mathfrak{X} \times \mathscr{H}, n = 1, 2, 3,\ldots . $
\pagebreak

On this basis, we next verify:   
    
\begin{itemize}
    \item[\textmd{(ev.4)}]$ \Phi_{\varepsilon_n}^R \geq 0 $ on $ \mathbb{X} \times H $, for $ n = 1, 2, 3, \dots $, and the union $ \bigcup_{n = 1}^\infty \bigl\{  \tilde{w} \in \mathbb{X} \times H \, \bigl| \, \Phi_{\varepsilon_n}^R(\tilde{w}) \leq r \bigr\} $ of sublevel sets is contained in the compact set $ K_\nu^R(r) \subset \mathbb{X} \times H $, as in (ev.3), for any $ r > 0 $;
    \item[\textmd{(ev.5)}]$ \Phi_{\varepsilon_n}^R \to \Phi_\varepsilon^R $ on $ \mathbb{X} \times H $, in the sense of Mosco, as $ n \to \infty $, more precisely, the uniform estimate \eqref{ev.M00} 
        will lead to the corresponding lower bound condition and optimality condition, in the Mosco-convergence of $ \{ \Phi_{\varepsilon_n}^R \}_{n=1}^{\infty} $;
    \item[\textmd{(ev.6)}]$ \sup_{n \in \mathbb{N}} \Phi_{\varepsilon_n}^R(w_{0, n}) < \infty $, and 
        \begin{align*}
        w_{0, n} \to w_{0} \ \mbox{in}\ \mathbb{X} \times H,\ \mbox{as}\ n \to \infty,
        \end{align*}        
        more precisely, it follows from \eqref{w.i}, (\hyperlink{A0l}{A0}), and (\hyperlink{A1l}{A1}) that 
        \begin{align*}
            \sup_{n \in \mathbb{N}} \Phi_{\varepsilon_n}^R(w_{0, n}) & \leq \sup_{n \in \mathbb{N}} \left( \frac{1+R}{2}|\eta_{0, n}|_V^2 +\nu^2 ( 1 +|\theta_{0, n}|_{V_0}^2) +\frac{1}{\nu^2} |\alpha(\eta_{0, n})|_H^2 \right) < \infty,
        \end{align*}
        and the weak convergence of $\{w_{0, n} \}_{n=1}^{\infty}$ in $\mathbb{W} \times V_{0}$ and the compactness of embedding $\mathbb{W} \times V_{0} \subset \mathbb{X} \times H$ imply the strong convergence of $\{w_{0, n} \}_{n=1}^{\infty}$ in $ \mathbb{X} \times H$.
\end{itemize}

On account of \eqref{w.i} and (ev.0)--(ev.6), we can apply Proposition \ref{Lem.CP02}, to show that:
\begin{subequations}\label{convKS}
\begin{equation}\label{convKS01}
    \begin{cases}
        w_n \to w \mbox{ in $ C([0, T]; \mathbb{X} \times H) $} 
        \\
        \quad \mbox{ (i.e. in $ C([0, T]; \mathbb{X}) \times C([0, T]; H) $)}, 
        \\
        \quad \mbox{ weakly in $ W^{1, 2}(0, T; \mathbb{X} \times H) $}
        \\
        \quad \mbox{ (i.e. weakly in $ W^{1, 2}(0, T; \mathbb{X}) \times W^{1, 2}(0, T; H) $),}
        \\[1ex]
        \displaystyle \int_0^T \Phi_{\varepsilon_n}^R(w_n(t)) \, dt \to \int_0^T \Phi_{\varepsilon}^R(w(t)) \, dt,
    \end{cases}
    \mbox{as $ n \to \infty $,}
\end{equation}
\begin{align*}
    \displaystyle \sup_{n \in \mathbb{N}} & |w_n|_{L^\infty(0, T; \mathbb{W}) \times L^\infty(0, T; V_0)}^2 \leq 4 \sup_{n \in \mathbb{N}} |w_n|_{L^\infty(0, T; \mathbb{W} \times V_0)}^2 
    \\
    & \displaystyle \leq \frac{8}{\min \, \{1, \nu^2, R \}} \sup_{n \in \mathbb{N}} |\Phi_{\varepsilon_n}^R(w_n)|_{L^\infty(0, T)} < \infty, 
\end{align*}
and hence,
\begin{equation}\label{convKS02}
    w_n \to w \mbox{ weakly-$*$ in $ L^\infty(0, T; \mathbb{W}) \times L^\infty(0, T; V_0) $, as $ n \to \infty $.}
\end{equation}
\end{subequations}
Also, as a consequence of the one-dimensional compact embeddings $ V \subset C(\overline{\Omega}) $ and $ V_0 \subset C(\overline{\Omega}) $, the uniqueness of solution $ w $ to (E)$_\varepsilon$, and Ascoli's theorem (cf. \cite[Corollary 4]{MR0916688}), we can derive from \eqref{convKS01} that
\begin{equation}\label{convKS07s}
    w_n \to w \mbox{ in $ [C(\overline{Q}) \times C(\overline{\Sigma})] \times C(\overline{Q}) $,\ as $ n \to \infty $.}
\end{equation}
Furthermore, from \eqref{f_eps}, \eqref{ev.M00}, \eqref{convKS}, \eqref{convKS07s}, and the assumptions (\hyperlink{A0l}{A0})--(\hyperlink{A2l}{A2}), one can observe that:
\begin{subequations}\label{convKSs}
\begin{equation}\label{convKS03}
    \begin{cases}
        \displaystyle \varliminf_{n \to \infty} \frac{1}{2} | \partial_x \eta_n |_{\mathscr{H}}^2 \geq \frac{1}{2} |\partial_x\eta|_{\mathscr{H}}^2, \quad \varliminf_{n \to \infty} \frac{R}{2} | \eta_n |_{\mathscr{H}}^2 \geq \frac{R}{2} |\eta|_{\mathscr{H}}^2, 
        \\[2ex]
        \displaystyle \varliminf_{n \to \infty} \frac{\nu^2}{2} |\theta_n |_{\mathscr{V}_0}^2 \geq \frac{\nu^2}{2} |\theta|_{\mathscr{V}_0}^2, \quad \lim_{n \to \infty} \frac{1}{2 \nu^2} |\alpha(\eta_n)|_{\mathscr{H}}^2 = \frac{1}{2 \nu^2} |\alpha(\eta)|_{\mathscr{H}}^2,
    \end{cases}
\end{equation}
and
\begin{align}\label{convKS04}
    \displaystyle \varliminf_{n \to \infty} & \displaystyle \bigl| \alpha(\eta_n) f_{\varepsilon_n}(\partial_x \theta_n) \big|_{L^1(Q)} = \varliminf_{n \to \infty} \int_0^T \int_\Omega \alpha(\eta_n(t)) f_{\varepsilon_n}(\partial_x \theta_n(t)) \, dx dt
    \nonumber
    \\
    & \geq \varliminf_{n \to \infty} \int_0^T \int_\Omega \alpha(\eta(t)) f_{\varepsilon_n}(\partial_x \theta_n(t)) \, dx dt
    \nonumber
    \\
    & \qquad -\lim_{n \to \infty} |\alpha(\eta_n) -\alpha(\eta)|_{C(\overline{Q})} \cdot \sup_{n \in \mathbb{N}} \bigl( T \varepsilon_n +|\partial_x \theta_n|_{L^1(0, T; L^1(\Omega))} \bigr)
    \nonumber
    \\
    & \geq \varliminf_{n \to \infty} \int_0^T \int_\Omega \alpha(\eta(t)) f_{\varepsilon}(\partial_x \theta_n(t)) \, dx dt -|\alpha(\eta)|_{C(\overline{Q})} \cdot \lim_{n \to \infty} \bigl( T|\varepsilon_n -\varepsilon| \bigr)
    \nonumber
    \\
    & \geq \int_0^T \int_\Omega \alpha(\eta(t)) f_{\varepsilon}(\partial_x \theta(t)) \, dx dt = \bigl| \alpha(\eta) f_{\varepsilon}(\partial_x \theta) \big|_{L^1(Q)}.
\end{align}
\end{subequations}
Here, from \eqref{Phi_eps}, it is seen that:
\begin{align}\label{convKS'}
    & \displaystyle \int_0^T {\Phi}_{\tilde{\varepsilon}}^R(\tilde{w}(t)) \, dt =  \int_0^T {\Phi}_{\tilde{\varepsilon}}^R(\tilde{\eta}(t), \tilde{\eta}_\Gamma(t), \tilde{\theta}(t)) \, dt
    \nonumber
    \\
    & ~ = \frac{1}{2} |\partial_x\tilde{\eta}|_{\mathscr{H}}^2 + \frac{R}{2}|\tilde{\eta} |_{\mathscr{H}}^2 +\frac{\nu^2}{2} |\tilde{\theta}|_{\mathscr{V}_0}^2 +\bigl| \alpha(\tilde{\eta}) f_{\tilde{\varepsilon}}(\partial_x \tilde{\theta}) \bigr|_{L^1(Q)} +\frac{1}{2 \nu^2} |\alpha(\tilde{\eta})|_{\mathscr{H}}^2 +\frac{\nu^2 \tilde{\varepsilon}^2}{2} T
    \nonumber
    \\[1ex]
    & \qquad \qquad \mbox{for all $ \tilde{\varepsilon} > 0 $ and $ \tilde{w} = [\tilde{\bm{\eta}}, \tilde{\theta} ] = [\tilde{\eta}, \tilde{\eta}_\Gamma, \tilde{\theta}] \in D(\Phi_{\tilde{\varepsilon}}) = \mathfrak{W} \times \mathscr{V}_0 $.}
\end{align}
Taking into account \eqref{convKS01}, \eqref{convKSs}, and \eqref{convKS'}, we deduce that:
\begin{align}\label{convKS05}
        |\partial_x\eta_n|_\mathscr{H}^2 & \, + R|\eta_n|_\mathscr{H}^2 +\nu^2 |\theta_n|_{\mathscr{V}_0}^2 \to |\partial_x\eta|_\mathscr{H}^2 + R|\eta|_\mathscr{H}^2 + \nu^2 |\theta|_{\mathscr{V}_0}^2, \nonumber
        \\
        & \, \mbox{and hence, } |[\eta_n, \theta_n]|_{\mathscr{V} \times \mathscr{V}_0} \to |[\eta, \theta]|_{\mathscr{V} \times \mathscr{V}_0}, \mbox{ as $ n \to \infty $.}
\end{align}

Since the norm of Hilbert space $ \mathscr{V} \times \mathscr{V}_0 $ is uniformly convex, the convergences \eqref{convKS02} and \eqref{convKS05} imply the strong convergences:
\begin{subequations}\label{convKS06-07}
\begin{equation}\label{convKS06}
    w_n \to w \mbox{ in $ \mathfrak{W} \times \mathscr{V}_0 $, as $ n \to \infty $,}
\end{equation}
and furthermore, it follows from \eqref{ev.M00} and \eqref{convKS06} that:
\begin{align}\label{convKS07}
    |f_{\varepsilon_n} (\partial_x \theta_n) & -f_\varepsilon(\partial_x \theta)|_\mathscr{H} \leq  |f_{\varepsilon_n}(\partial_x \theta_n) -f_{\varepsilon_n}(\partial_x \theta)|_\mathscr{H} +|f_{\varepsilon_n}(\partial_x \theta) -f_{\varepsilon}(\partial_x \theta)|_\mathscr{H} 
    \nonumber
    \\
    \leq &  |\theta_n -\theta|_{\mathscr{V}_0} +\sqrt{T} |\varepsilon_n -\varepsilon| \to 0, \mbox{ as $ n \to \infty $.}
\end{align}
\end{subequations}

The convergences \eqref{convKS}, \eqref{convKS07s}, and \eqref{convKS06-07} 
are sufficient to verify the conclusions \eqref{mThConv} and \eqref{mThConv00} of Main Theorem \ref{mainTh01} (I-B). 
\qed
}

\section{Proof of Main Theorem \ref{mainTh02}}

In this Section, we prove the second Main Theorem \ref{mainTh02}. Let $[\bm{\eta}_0, \theta_0] = [\eta_0, \eta_{\Gamma, 0}, \theta_0] \in \mathbb{W} \times V_0$ with $\bm{\eta_0} = [\eta_0, \eta_{\Gamma, 0}] $ be the initial triplet. Also, let us fix arbitrary forcing triplet $ [\bar{\bm{u}}, \bar{v}] = [\bar{u}, \bar{u}_\Gamma, \bar{v}] \in \mathfrak{X} \times \mathscr{H} $ with $\bar{\bm{u}} = [\bar{u}, \bar{u}_\Gamma] $ , and let us invoke the definition of the cost functional \eqref{J}, to estimate that:
\begin{align}\label{mTh02-00}
    0 &~ \leq \underline{J}_{\varepsilon}:= \inf_{[\bm{u}, v] \in \mathfrak{X} \times \mathscr{H}} \mathcal{J}_{\varepsilon}(\bm{u}, v) \leq \overline{J}_{\varepsilon}:= \mathcal{J}_\varepsilon(\bar{\bm{u}}, \bar{v}) = \mathcal{J}_\varepsilon(\bar{u}, \bar{u}_\Gamma, \bar{v}) < \infty, \mbox{ for all } \varepsilon \geq 0.
\end{align}
Also, for any $ \varepsilon \geq 0 $, we denote by $[\bar{\bm{\eta}}_{\varepsilon}, \bar{\theta}_{\varepsilon}] = [\bar{\eta}_\varepsilon, \bar{\eta}_{\Gamma, \varepsilon}, \bar{\theta}_\varepsilon] \in \mathfrak{X} \times \mathscr{H}$ with $\bar{\bm{\eta}}_\varepsilon = [\bar{\eta}_\varepsilon, \bar{\eta}_{\Gamma, \varepsilon}] $ the solution to (S)$_{\varepsilon}$, for the initial triplet $[\bm{\eta}_0, \theta_0] = [\eta_0, \eta_{\Gamma, 0}, \theta_0] $ and forcing triplet $ [\bar{\bm{u}}, \bar{v}] = [\bar{u}, \bar{u}_\Gamma, \bar{v}]$.
\bigskip

Based on these, the Main Theorem \ref{mainTh02} is proved as follows.

\paragraph{\textbf{Proof of Main Theorem \ref{mainTh02} (II-A).}}
Let us fix any $\varepsilon \geq 0$. 
Then, from the estimate \eqref{mTh02-00}, we immediately find a sequence of forcing triplets $\{[\bm{u}_{n}, v_{n}] \}_{n=1}^{\infty} = \{ [u_{n}, u_{\Gamma, n}, v_n] \}_{n=1}^{\infty}$ $\subset \mathfrak{X} \times \mathscr{H}$ with $\{\bm{u}_n \}_{n=1}^\infty = \{ [u_{n}, u_{\Gamma, n}] \}_{n=1}^{\infty} $, such that:
\begin{subequations}\label{mTh02-0102}
\begin{equation}\label{mTh02-01}
\mathcal{J}_{\varepsilon}(\bm{u}_{n}, v_{n}) = \mathcal{J}_{\varepsilon}(u_{n}, u_{\Gamma, n}, v_{n}) \downarrow \underline{J}_{\varepsilon},\ \mbox{as}\ n \to \infty,
\end{equation}
and
\begin{align}\label{mTh02-02}
    \frac{1}{2} \sup_{n \in \mathbb{N}}& \bigl| [\sqrt{L} u_{n}, \sqrt{L_\Gamma} u_{\Gamma, n}, \sqrt{M} v_{n}] \bigr|_{\mathfrak{X} \times \mathscr{H}}^{2} \leq \mathcal{J}_\varepsilon(\bar{u}, \bar{u}_\Gamma, \bar{v}) < \infty.
\end{align}
\end{subequations}
Also, the estimate \eqref{mTh02-02} enables us to take a subsequence of $\{[\bm{u}_{n}, v_{n}] \}_{n=1}^{\infty} = $ \linebreak $ \{ [u_{n}, u_{\Gamma, n}, v_n] \}_{n=1}^{\infty} $ (not relabeled), and to find a triplet of functions $[\bm{u}^{*}, v^{*}] = [u^*, u^*_\Gamma, v^*] \in \mathfrak{X} \times \mathscr{H}$ with $\bm{u}^* = [u^*, u^*_\Gamma]$, such that:
\begin{subequations}\label{mTh02-03}
\begin{align}\label{mTh02-03a}
[\sqrt{L} u_{n}, &\, \sqrt{L_\Gamma} u_\Gamma, \sqrt{M} v_{n}] \to [\sqrt{L} u^{*}, \sqrt{L_\Gamma} u_\Gamma^*, \sqrt{M} v^{*}]
\nonumber
\\
& \mbox{ weakly in } \mathfrak{X} \times \mathscr{H}, \mbox{ as } n \to \infty, 
\end{align}
and as well as,
\begin{align}\label{mTh02-03b}
[L u_{n}, L_\Gamma u_{\Gamma, n}, M v_{n}] & \to [L u^{*}, L_\Gamma u^*_\Gamma, M_v v^{*} ] \mbox{ weakly in } \mathfrak{X} \times \mathscr{H}, \mbox{ as } n \to \infty.
\end{align}
\end{subequations}

Let $[\bm{\eta}^*, \theta^* ] = [\eta^{*}, \eta^*_{\Gamma}, \theta^*] \in \mathfrak{X} \times \mathscr{H}$ with $ \bm{\eta}^* = [\eta^{*}, \eta^*_{\Gamma}] $ be the solution to (S)$_{\varepsilon}$, for the initial triplet $[\bm{\eta}_0, \theta_0 ] = [\eta_0, \eta_{\Gamma, 0}, \theta_0] $ and forcing triplet $[\bm{u}^*, v^* ] = [u^*, u_\Gamma^*, v^*]$. 
As well as, for any $n \in \mathbb{N}$, let $[\bm{\eta}_n, \theta_n ] = [\eta_{n}, \eta_{\Gamma, n}, \theta_n] \in \mathfrak{X} \times \mathscr{H}$ with $\bm{\eta}_n = [\eta_{n}, \eta_{\Gamma, n}] $ be the solution to (S)$_{\varepsilon}$, for the initial triplet $[\bm{\eta}_0, \theta_0] = [\eta_0, \eta_{\Gamma, 0}, \theta_0]$ and the forcing triplet $[\bm{u}_n, v_n ] = [u_n, u_{\Gamma, n}, v_n]$.
Then, having in mind \eqref{mTh02-03} and the initial condition:
\begin{align*}
     [\bm{\eta}_{n}(0),&\, \theta_n(0)] = [\eta_{n}(0), \eta_{\Gamma, n}(0), \theta_{n}(0)]\\
    & = [\bm{\eta}^{*}(0), \theta^{*}(0)]  = [\eta^*(0), \eta^*_\Gamma(0), \theta^*(0) ]\\
    & = [\bm{\eta}_0, \theta_0 ] = [\eta_{0}, \eta_{\Gamma, 0}, \theta_{0}] \mbox{ in $\mathbb{X} \times H$, for $n = 1, 2, 3, \dots $,}
\end{align*}
we can apply Main Theorem \ref{mainTh01} (I-B), to see that:
\begin{align}\label{mTh02-04}
[\bm{\eta}_n, \theta_n ] \to [\bm{\eta}^{*}, \theta^{*}] \mbox{ in } [C(\overline{Q}) \times C(\overline{\Sigma})] \times C(\overline{Q}), \mbox{ as } n \to \infty.
\end{align}
On account of \eqref{mTh02-01}, \eqref{mTh02-03a}, and \eqref{mTh02-04}, it is computed that:
\begin{align*}
    \mathcal{J}_{\varepsilon}&\,(\bm{u}^{*}, v^{*}) = \mathcal{J}_{\varepsilon}(u^*, u^*_\Gamma, v^*) \\
                                            & = \frac{1}{2} \bigl| [\sqrt{K} (\eta^{*}-\eta_\mathrm{ad}), \sqrt{K_\Gamma} (\eta^{*}_{\Gamma} - \eta_{\Gamma, \mathrm{ad}}), \sqrt{\Lambda} (\theta^{*}-\theta_\mathrm{ad})] \bigr|_{\mathfrak{X} \times \mathscr{H}}^{2}\\
                                            &\qquad + \frac{1}{2} \bigl| [\sqrt{L} u^{*}, \sqrt{L_\Gamma} u^*_\Gamma, \sqrt{M} v^{*}] \bigr|_{\mathfrak{X} \times \mathscr{H}}^{2}\\
                                            & \leq \frac{1}{2}\lim_{n \to \infty} \bigl| [\sqrt{K} (\eta_{n}-\eta_\mathrm{ad}), \sqrt{K_\Gamma} (\eta_{\Gamma, n} - \eta_{\Gamma, \mathrm{ad}}), \sqrt{\Lambda} (\theta_{n}-\theta_\mathrm{ad}) ] \bigr|_{\mathfrak{X} \times \mathscr{H}}^{2}\\
                                            &\qquad + \frac{1}{2} \varliminf_{n \to \infty} \bigl| [\sqrt{L} u_{n}, \sqrt{L_\Gamma} u_{\Gamma, n}, \sqrt{M} v_{n}] \bigr|_{\mathfrak{X} \times \mathscr{H}}^{2}
    \\
    & = \lim_{n \to \infty}\mathcal{J}_{\varepsilon}(\bm{u}_{n}, v_{n}) = \lim_{n \to \infty}\mathcal{J}_{\varepsilon}(u_{n}, u_{\Gamma, n}, v_{n}) = \underline{J}_{\varepsilon}~ (\leq \mathcal{J}_{\varepsilon}(\bm{u}^{*}, v^{*})),
\end{align*}
and it implies that
\begin{equation*}
    \mathcal{J}_\varepsilon(\bm{u}^{*}, v^{*}) = \mathcal{J}_\varepsilon(u^{*}, u^*_\Gamma, v^{*}) = \min_{[\bm{u}, v] \in \mathfrak{X} \times \mathscr{H}} \mathcal{J}_{\varepsilon}(\bm{u}, v) = \min_{[u, u_\Gamma, v] \in \mathfrak{X} \times \mathscr{H}} \mathcal{J}_{\varepsilon}(u, u_\Gamma, v).
\end{equation*}

Thus, we conclude the item (II-A). 
\qed

\paragraph{\textbf{Proof of Main Theorem \ref{mainTh02} (II-B).}}
Let $ \varepsilon  \in [0, 1] $ and $\{\varepsilon_{n} \}_{n=1}^{\infty} \subset [0, 1]$ be as in \eqref{ass.4}. Let $[\bm{\bar{\eta}}_{\varepsilon}, \bar{\theta}_{\varepsilon}] = [\bar{\eta}_\varepsilon, \bar{\eta}_{\Gamma, \varepsilon}, \bar{\theta}_\varepsilon] \in \mathfrak{X} \times \mathscr{H}$ with $\bar{\bm{\eta}}_\varepsilon = [\bar{\eta}_\varepsilon, \bar{\eta}_{\Gamma, \varepsilon}] $  be the solution to the system (S)$_\varepsilon$, for the initial triplet $ [\bm{\eta}_0, \theta_0] = [\eta_0, \eta_{\Gamma, 0}, \theta_0] $ and forcing triplet $ [\bm{\bar{u}}, \bar{v}] = [\bar{u}, \bar{u}_\Gamma, \bar{v}]$, and let $[\bm{\bar{\eta}}_{\varepsilon_{n}}, \bar{\theta}_{\varepsilon_{n}}] = [\bar{\eta}_{\varepsilon_{n}}, \bar{\eta}_{\Gamma, \varepsilon_n}, \bar{\theta}_{\varepsilon_n}] \in \mathfrak{X} \times \mathscr{H}$ with $\bar{\bm{\eta}}_{\varepsilon_n} = [\bar{\eta}_{\varepsilon_{n}}, \bar{\eta}_{\Gamma, \varepsilon_n}] $, $n=1, 2, 3,\ldots ,$ be the solutions to (S)$_{\varepsilon_n}$, for the respective initial triplets $ [\bm{\eta}_{0, n}, \theta_{0, n}] = [\eta_{0, n}, \eta_{\Gamma, 0, n}, \theta_{0, n}] \in \mathbb{W} \times V_0$ with $\bm{\eta}_{0, n} = [\eta_{0, n}, \eta_{\Gamma, 0, n}] $, $ n = 1, 2, 3, \dots $, and the fixed forcing triplet $ [\bar{\bm{u}}, \bar{v}] = [\bar{u}, \bar{u}_\Gamma, \bar{v}] $. On this basis, let us first apply Main Theorem \ref{mainTh01} (I-B) to the solutions $[\bm{\bar{\eta}}_{\varepsilon}, \bar{\theta}_{\varepsilon}] = [\bar{\eta}_\varepsilon, \bar{\eta}_{\Gamma, \varepsilon}, \bar{\theta}_\varepsilon]$ and $[\bm{\bar{\eta}}_{\varepsilon_{n}}, \bar{\theta}_{\varepsilon_{n}}] = [\bar{\eta}_{\varepsilon_n}, \bar{\eta}_{\Gamma, \varepsilon_n}, \bar{\theta}_{\varepsilon_n}]$, $n=1, 2, 3,\ldots .$ Then, we have 
\begin{equation}\label{mTh02-09}
    \begin{cases}
        [\bm{\bar{\eta}}_{\varepsilon_{n}}, \theta_{\varepsilon_{n}}] \to [\bm{\bar{\eta}}_{\varepsilon}, \theta_{\varepsilon}] \mbox{ in } [C(\overline{Q}) \times C(\overline{\Sigma})] \times C(\overline{Q}),
        \\[2ex]
        [\bm{\bar{\eta}}_{n}(0), \bar{\theta}_{n}(0)]  = [\bm{\eta}_{0, n}, \theta_{0, n}] \to [\bm{\bar{\eta}}_\varepsilon(0), \bar{\theta}_\varepsilon(0)] = [\bm{\eta}_0, \theta_0] 
        \\[1ex]
        \hspace{4ex}\mbox{in $ [C(\overline{\Omega}) \times C(\Gamma)] \times C(\overline{\Omega}) $,}
        \mbox{ as } n \to \infty,
    \end{cases}
\end{equation}
and hence, 
\begin{equation}\label{mTh02-10}
    \overline{J}_\mathrm{sup} := \sup_{n \in \mathbb{N}} \mathcal{J}_{\varepsilon_n}(\bm{\bar{u}}, \bar{v}) = \sup_{n \in \mathbb{N}} \mathcal{J}_{\varepsilon_n}(\bar{u}, \bar{u}_\Gamma, \bar{v}) < \infty.
\end{equation}

Next, for any $n \in \mathbb{N}$, let us denote by $[\bm{\eta}^{*}_{n}, \theta^{*}_{n}] = [\eta^*_n, \eta^*_{\Gamma, n}, \theta^*_n] \in \mathfrak{X} \times \mathscr{H}$ with $\bm{\eta}^*_n = [\eta^*_n, \eta^*_{\Gamma, n}] $ the solution to (S)$_{\varepsilon_{n}}$, for the initial triplet $[\bm{\eta}_{0, n}, \theta_{0, n}] = [\eta_{0, n}, \eta_{\Gamma, 0, n}, \theta_{0, n}]$ and forcing triplet $[\bm{u}^{*}_{n}, v^{*}_{n}] = [u^*_n, u^*_{\Gamma, n}, v^*_n] \in \mathfrak{X} \times \mathscr{H}$ with $\bm{u}^*_n = [u^*_n, u^*_{\Gamma, n}] $. Then, in the light of \eqref{mTh02-00} and \eqref{mTh02-10}, we can see that:
\begin{align*}
    0 & \leq \frac{1}{2}|[\sqrt{L} u^{*}_{n}, \sqrt{L_\Gamma} u^*_{\Gamma, n}, \sqrt{M} v^{*}_{n}] |_{\mathfrak{X}\times \mathscr{H}}^{2} \leq \underline{J}_{\varepsilon_{n}} \leq \overline{J}_\mathrm{sup} < \infty,  \mbox{ for } n = 1, 2, 3, \dots.
\end{align*}
Therefore, we can find a subsequence $\{n_{i} \}_{i=1}^{\infty} \subset \{n\}$, together with a triplet of functions $[\bm{u}^{**}, v^{**}] = [u^{**}, u^{**}_{\Gamma, n}, v^{**}] \in \mathfrak{X} \times \mathscr{H}$ with $\bm{u}^{**} = [u^{**}, u^{**}_{\Gamma, n}]$, such that:
\begin{align}\label{mTh02-11}
[\sqrt{L} u^{*}_{n_i}, &\, \sqrt{L_\Gamma} u^*_{\Gamma, n_i}, \sqrt{M} v^{*}_{n_i}]  \to [\sqrt{L} u^{**}, \sqrt{L_\Gamma} {u}^{**}_\Gamma, \sqrt{M} v^{**}]
\nonumber
\\
& \mbox{ weakly in } \mathfrak{X} \times \mathscr{H}, \mbox{ as } i \to \infty, 
\end{align}
and as well as,
\begin{align*}
    [L u^{*}_{n_i}, L_\Gamma u^*_{\Gamma, n_i}, M v^{*}_{n_i}] & \to [L u^{**}, L_\Gamma u^{**}_\Gamma, M v^{**}] \mbox{ weakly in } \mathfrak{X} \times \mathscr{H}, \mbox{ as } i \to \infty.
\end{align*}
Here, let us denote by $[\bm{\eta}^{**}, \theta^{**}] = [\eta^{**}, \eta^{**}_\Gamma, \theta^{**}] \in \mathfrak{X} \times \mathscr{H}$ with $\bm{\eta}^{**} = [\eta^{**}, \eta^{**}_\Gamma] $ the solution to (S)$_{\varepsilon}$, for the initial triplet $ [\bm{\eta}_0, \theta_0] = [\eta_{0, n}, \eta_{\Gamma, 0, n}, \theta_{0, n}]$ and forcing triplet $[\bm{u}^{**}, v^{**}] = [u^{**}, u^{**}_\Gamma, v^{**}]$.
Then, applying Main Theorem \ref{mainTh01} (I-B), again, to the solutions $[\bm{\eta}^{**}, \theta^{**}] = [\eta^{**}, \eta^{**}_\Gamma, \theta^{**}] $ and $[\bm{\eta}^{*}_{n_i}, \theta^{*}_{n_i}] = [\eta^*_{n_i}, \eta^*_{\Gamma, n_i}, \theta^*_{n_i}]$, $i = 1, 2, 3, \dots$, we can observe that:
\begin{align}\label{mTh02-12}
[\bm{\eta}^{*}_{n_i}, \theta^{*}_{n_i}] \to [\bm{\eta}^{**}, \theta^{**}] \mbox{ in } [C(\overline{Q}) \times C(\overline{\Sigma})] \times C(\overline{Q}), \mbox{ as } i \to \infty.
\end{align}

As a consequence of \eqref{mTh02-09}, \eqref{mTh02-11}, and \eqref{mTh02-12}, it is verified that:
\begin{align*}
    \mathcal{J}_{\varepsilon} (\bm{u}^{**},&\, v^{**}) = \mathcal{J}_{\varepsilon} (u^{**}, u^{**}_\Gamma, v^{**}) \\
                                                   & = \frac{1}{2} \bigl| [\sqrt{K} (\eta^{**}-\eta_\mathrm{ad}), \sqrt{K_\Gamma} (\eta^{**}_\Gamma - \eta_{\Gamma, \mathrm{ad}}), \sqrt{\Lambda} (\theta^{**}-\theta_\mathrm{ad}))] \bigr|_{\mathfrak{X} \times \mathscr{H}}^{2}
    \\
    &\qquad + \frac{1}{2} \bigl| [\sqrt{L} u^{**}, \sqrt{L_\Gamma} u^{**}_\Gamma, \sqrt{M} v^{**}] \bigr|_{\mathfrak{X} \times \mathscr{H}}^{2}
    \\
    & \leq \frac{1}{2} \lim_{i \to \infty} \bigl| [\sqrt{K} (\eta^{*}_{n_i} -\eta_\mathrm{ad}), \sqrt{K_\Gamma} (\eta^*_{\Gamma, n_i} - \eta_{\Gamma, \mathrm{ad}}), \sqrt{\Lambda} (\theta^{*}_{n_i} -\theta_\mathrm{ad})] \bigr|_{\mathfrak{X} \times \mathscr{H}}^{2}
    \\
    &\qquad  + \frac{1}{2}\varliminf_{i \to \infty} \bigl| [\sqrt{L} u^{*}_{n_i}, \sqrt{L_\Gamma} u^{*}_{\Gamma, n_i}, \sqrt{M} v^{*}_{n_i}] \bigr|_{\mathfrak{X} \times \mathscr{H}}^{2}
    \\
    & = \varliminf_{i \to \infty}\mathcal{J}_{\varepsilon_{n_i}}(\bm{u}^{*}_{n_i}, v^{*}_{n_i})  = \varliminf_{i \to \infty}\mathcal{J}_{\varepsilon_{n_i}}(u^{*}_{n_i}, u^{*}_{\Gamma, n_i}, v^{*}_{n_i})
    \\
    & \leq \lim_{i \to \infty}\mathcal{J}_{\varepsilon_{n_i}}(\bm{\bar{u}}, \bar{v}) = \lim_{i \to \infty}\mathcal{J}_{\varepsilon_{n_i}}(\bar{u}, \bar{u}_\Gamma, \bar{v})
    \\
    & = \frac{1}{2}\lim_{i \to \infty} \bigl| [\sqrt{K} (\bar{\eta}_{\varepsilon_{n_i}}-\eta_\mathrm{ad}), \sqrt{K_\Gamma} (\bar{\eta}_{\Gamma, \varepsilon_{n_i}} - \eta_{\Gamma, \mathrm{ad}}), \sqrt{\Lambda} (\bar{\theta}_{\varepsilon_{n_i}}-\theta_\mathrm{ad})] \bigr|_{\mathfrak{X} \times \mathscr{H}}^{2}
    \\
    &\qquad  + \frac{1}{2}\bigl| [\sqrt{L} \bar{u}, \sqrt{L_\Gamma} \bar{u}_\Gamma, \sqrt{M} \bar{v}] \bigr|_{\mathfrak{X} \times \mathscr{H}}^{2}
    \\
    & = \mathcal{J}_{\varepsilon}(\bm{\bar{u}}, \bar{v}) = \mathcal{J}_{\varepsilon}(\bar{u}, \bar{u}_\Gamma, \bar{v}).
\end{align*}
Since the choice of $[\bm{\bar{u}}, \bar{v}] = [\bar{u}, \bar{u}_\Gamma, \bar{v}] \in \mathfrak{X} \times \mathscr{H}$ is arbitrary, we conclude that:
\begin{equation*}
    \mathcal{J}_\varepsilon(\bm{u}^{**}, v^{**}) = \mathcal{J}_\varepsilon(u^{**}, u^{**}_\Gamma, v^{**}) = \min_{[\bm{u}, v] \in \mathfrak{X} \times \mathscr{H}} \mathcal{J}_{\varepsilon}(\bm{u}, v) = \min_{[u, u_\Gamma, v] \in \mathfrak{X} \times \mathscr{H}} \mathcal{J}_{\varepsilon}(u, u_\Gamma, v),
\end{equation*}
and complete the proof of the item (II-B).
\qed

\section{Proof of Main Theorem \ref{mainTh03}}

This Section is devoted to the proof of the third Main Theorem \ref{mainTh03}. To this end, we need to start with the case of $ \varepsilon > 0 $, and prepare some Lemmas, associated with the G\^{a}teaux differential of the regular cost functional $ \mathcal{J}_\varepsilon $. 
\medskip

Let $ \varepsilon > 0 $ be a fixed constant, and let $ [\bm{\eta}_0, \theta_0] = [\eta_0, \eta_{\Gamma, 0}, \theta_0] \in \mathbb{W} \times V_0 $ with $\bm{\eta}_0 = [\eta_0, \eta_{\Gamma, 0}] $ be the initial triplet. Let us take any  forcing triplet $ [\bm{u}, v] = [u, u_\Gamma, v] \in \mathfrak{X} \times \mathscr{H} $ with $\bm{u} = [u, u_\Gamma]$, and consider the unique solution $ [\bm{\eta}, \theta] = [\eta, \eta_\Gamma, \theta] \in \mathfrak{X} \times \mathscr{H} $ with $\bm{\eta} = [\eta, \eta_\Gamma] $ to the state-system (S)$_\varepsilon$. Also, let us take any constant $ \delta \in (-1, 1) \setminus \{0\} $ and any triplet of functions $ [\bm{h}, k] = [h, h_\Gamma, k] \in \mathfrak{X} \times \mathscr{H} $ with $\bm{h} = [h, h_\Gamma] $, and consider another solution $ [\bm{\eta}^\delta, \theta^\delta] = [\eta^\delta, \eta^\delta_\Gamma, \theta^\delta] \in \mathfrak{X} \times \mathscr{H} $ with $\bm{\eta}^\delta = [\eta^\delta, \eta^\delta_\Gamma] $ to the system (S)$_\varepsilon$, for the initial triplet $ [\bm{\eta}_0, \theta_0] = [\eta_0, \eta_{\Gamma, 0}, \theta_0] $ and a perturbed forcing triplet $ [\bm{u} +\delta \bm{h}, v +\delta k] = [u + \delta h, u_\Gamma + \delta h_\Gamma, v + \delta k] \in \mathfrak{X} \times \mathscr{H} $ with $\bm{u} + \delta\bm{h} = [u + \delta h, u_\Gamma + \delta h_\Gamma] $. On this basis, we consider a sequence of triplets of functions $ \{ [\bm{\chi}^\delta, \gamma^\delta] \}_{\delta \in (-1, 1) \setminus \{0\}} = \{[\chi^\delta, \chi^\delta_\Gamma, \gamma^\delta ]\}_{\delta \in (-1, 1)\setminus \{0\}} \subset \mathfrak{X} \times \mathscr{H} $ with $\{\bm{\chi}^\delta \}_{\delta \in (-1, 1)\setminus \{0\}} = \{[\chi^\delta, \chi^\delta_\Gamma ]\}_{\delta \in (-1, 1)\setminus \{0\}}$, defined as:
\begin{align}\label{set00}
    [\bm{\chi}^\delta, \gamma^\delta] & \, = [\chi^\delta, \chi_\Gamma^\delta, \gamma^\delta] := \left[ \frac{\bm{\eta}^\delta -\bm{\eta}}{\delta}, ~ \frac{\theta^\delta -\theta}{\delta} \right] = \left[\frac{\eta^\delta - \eta}{\delta}, ~ \frac{\eta^\delta_\Gamma - \eta_\Gamma}{\delta}, ~ \frac{\theta^\delta - \theta}{\delta} \right] \in \mathfrak{X} \times \mathscr{H}\nonumber
    \\
    &\, \mbox{ with } \bm{\chi}^\delta = [\chi^\delta, \chi_\Gamma^\delta] = \left[\frac{\eta^\delta - \eta}{\delta}, ~ \frac{\eta^\delta_\Gamma - \eta_\Gamma}{\delta} \right] \mbox{ for $ \delta \in (-1, 1) \setminus \{0\} $.}
\end{align}
This sequence acts a key-role in the computation of G\^{a}teaux differential of the cost functional $ \mathcal{J}_\varepsilon $, for $ \varepsilon > 0 $. 
\begin{remark}\label{Rem.GD01}
    Note that for any  $ \delta \in (-1, 1) \setminus \{0\} $, the triplet of functions $ [\bm{\chi}^\delta, \gamma^\delta] = [\chi^\delta, \chi^\delta_\Gamma, \gamma^\delta] \in \mathfrak{X} \times \mathscr{H} $ with $\bm{\chi}^\delta = [\chi^\delta, \chi^\delta_\Gamma]$
 fulfills the following variational forms:
\begin{align*}
    (\partial_{t} & \bm{\chi}^\delta(t), \bm{\varphi})_{\mathbb{X}}  + (\partial_{x} \chi^\delta(t), \partial_{x}\varphi)_{H} 
    \\
    & +\int_{\Omega} \left(\int_{0}^{1} g'(\eta(t)+\varsigma \delta \chi^\delta(t)) \, d \varsigma \right) \chi^\delta(t) \varphi \, dx
    \\
    &  + \int_{\Omega} \left(f_{\varepsilon}(\partial_{x} \theta(t)) \int_{0}^{1} \alpha''(\eta(t) +\varsigma \delta \chi^\delta(t)) \, d\varsigma \right) \chi^\delta(t) \varphi \, dx
    \\
    & +\int_{\Omega}\left(\alpha'(\eta^\delta(t)) \int_{0}^{1} f_{\varepsilon}'(\partial_{x}\theta(t) +\varsigma \delta \partial_{x} \gamma^\delta(t)) \, d\varsigma \right) \partial_{x}\gamma^\delta(t) \varphi \, dx
    \\
    = & (L h(t), \varphi)_{H} + (L_\Gamma h_\Gamma(t), \varphi_\Gamma)_{H_\Gamma}, \mbox{ for any $ \bm{\varphi} = [\varphi, \varphi_\Gamma] \in \mathbb{W} $}, 
    \\
     & \qquad \mbox{a.e. $ t \in (0, T) $, subject to $ \bm{\chi}^\delta(0) = [\chi^\delta(0), \chi^\delta_\Gamma(0)] = [0, 0] $ in $ \mathbb{X} $,}
\end{align*}
and
\begin{align*}
    (\alpha_{0}(t) & \partial_{t} \gamma^\delta(t), \psi)_{H} + \nu^{2}(\partial_{x} \gamma^\delta(t), \partial_{x}\psi)_{H} 
    \nonumber 
    \\
 & + \int_{\Omega} \left(\alpha(\eta^\delta(t)) \int_{0}^{1} f_{\varepsilon}''(\partial_{x}\theta(t) + \varsigma \delta \partial_{x} \gamma^\delta(t)) \, d \varsigma \right) \partial_{x} \gamma^\delta(t) \partial_{x}\psi \, dx 
    \nonumber 
    \\
& + \int_{\Omega}\left(f'_{\varepsilon}(\partial_{x}\theta(t)) \int_{0}^{1}\alpha'(\eta(t) +\varsigma \delta \chi^\delta(t)) \, d \varsigma \right) \chi^\delta(t) \partial_{x} \psi \, dx 
    \nonumber
    \\
    = & (M k(t), \psi)_{H}, \mbox{ for any $ \psi \in V_{0} $,  a.e. $ t \in (0, T) $, subject to $ \gamma^\delta(0) = 0 $ in $ H $.}
\end{align*}
    In fact, these variational forms are obtained by taking the difference between respective two variational forms for $ [\bm{\eta}^\delta, \theta^\delta] = [\eta^\delta, \eta_\Gamma^\delta, \theta^\delta] \in \mathfrak{X} \times \mathscr{H}$ and $ [\bm{\eta}, \theta] = [\eta, \eta_\Gamma, \theta] \in \mathfrak{X} \times \mathscr{H}$, as in Main Theorem \ref{mainTh01} (I-A), and by using the following linearization formulas: 

\begin{align*}
    & \frac{1}{\delta} \bigl( g(\eta^{\delta})-g(\eta) \bigr) = \left( \int_{0}^{1}g'(\eta + \varsigma \delta \chi^\delta) \, d\varsigma \right) \chi^{\delta} \mbox{ in $ \mathscr{H} $,}
\end{align*}
\begin{align*}
    & \frac{1}{\delta} \bigl( \alpha' (\eta^{\delta})f_{\varepsilon}(\partial_{x}\theta^{\delta}) - \alpha'(\eta)f_{\varepsilon}(\partial_{x}\theta) \bigr)
    \\
    & \qquad = \frac{1}{\delta} \bigl( \alpha' (\eta^{\delta})-\alpha'(\eta)\bigr)f_{\varepsilon}(\partial_{x} \theta) + \frac{1}{\delta}  \alpha'(\eta^\delta)\bigl(f_{\varepsilon}(\partial_{x}\theta^{\delta})-f_{\varepsilon}(\partial_{x}\theta) \bigr)
    \\
    & \qquad = \left( f_{\varepsilon}(\partial_{x}\theta) \int_{0}^{1}\alpha''(\eta+\varsigma \delta \chi^\delta) \, d\varsigma \right) \chi^{\delta}
    \\
    &  \qquad \qquad + \left( \alpha'(\eta^{\delta}) \int_{0}^{1}f_{\varepsilon}'(\partial_{x}\theta +\varsigma \delta \partial_{x} \gamma^\delta) \, d\varsigma \right) \partial_{x} \gamma^{\delta} \mbox{ in $ \mathscr{H} $,}
\end{align*}
and
\begin{align*}
    & \frac{1}{\delta} \bigl( \alpha(\eta^{\delta}) f'_{\varepsilon}(\partial_{x}\theta^{\delta}) - \alpha(\eta)f'_{\varepsilon}(\partial_{x}\theta) \bigr)
    \\
    & \qquad = \frac{1}{\delta} \alpha(\eta^{\delta})\bigl(f'_{\varepsilon}(\partial_{x}\theta^{\delta})-f'_{\varepsilon}(\partial_{x}\theta)\bigr)  + \frac{1}{\delta} \bigl( \alpha (\eta^{\delta})-\alpha(\eta)\bigr)f'_{\varepsilon}(\partial_{x}\theta) 
    \\
    & \qquad = \left( \alpha(\eta^{\delta}) \int_{0}^{1}f_{\varepsilon}''(\partial_{x}\theta +\varsigma \delta \partial_{x} \gamma^\delta) \, d\varsigma \right) \partial_{x}\gamma^{\delta}
    \\
    & \qquad \qquad +\left( f'_{\varepsilon}(\partial_{x}\theta) \int_{0}^{1}\alpha'(\eta +\varsigma \delta \chi^\delta) \, d \varsigma \right) \chi^{\delta} \mbox{ in $ \mathscr{H} $.}
\end{align*}
Incidentally, the above linearization formulas can be verified as consequences of the assumptions (\hyperlink{A0l}{A0})--(\hyperlink{A3l}{A3}) and the mean-value theorem (cf. \cite[Theorem 5 in p. 313]{lang1968analysisI}).
\end{remark}

Now, we verify the following two Lemmas. 
\begin{lemma}\label{Lem.GD01}
    Let us fix $ \varepsilon > 0 $, and assume (\hyperlink{A0l}{A0})--(\hyperlink{A3l}{A3}). Then, for any $ [\bm{u}, v] = [u, u_\Gamma, v] \in \mathfrak{X} \times \mathscr{H}$ with $\bm{u} = [u, u_\Gamma]$, the cost functional $ \mathcal{J}_\varepsilon $ admits the G\^{a}teaux derivative $ \mathcal{J}_\varepsilon'(\bm{u}, v) = \mathcal{J}_\varepsilon'(u, u_\Gamma, v) \in \mathfrak{X} \times \mathscr{H} $ $ (= [\mathfrak{X} \times \mathscr{H} ]^*) $, such that:
    \begin{align}\label{GD01}
        \bigl( \mathcal{J}_\varepsilon'(\bm{u}, v), &\, [\bm{h}, k] \bigr)_{\mathfrak{X} \times \mathscr{H}} = \bigl( \mathcal{J}_\varepsilon'(u, u_\Gamma, v), [h, h_\Gamma, k] \bigr)_{\mathfrak{X} \times \mathscr{H}}
        \nonumber
        \\
        & = \bigl( [K (\eta -\eta_\mathrm{ad}), K_\Gamma (\eta_\Gamma -\eta_{\Gamma, \mathrm{ad}}), \Lambda (\theta -\theta_\mathrm{ad})], \bar{\mathcal{P}}_\varepsilon [L h, L_\Gamma u_\Gamma, M k] \bigr)_{\mathfrak{X}\times \mathscr{H}}
        \nonumber
        \\
        & \qquad +\bigl( [L u, L_\Gamma u_\Gamma, M v], [h, h_\Gamma, k] \bigr)_{\mathfrak{X} \times \mathscr{H}},
        \nonumber
        \\
        &  \mbox{ for any $ [\bm{h}, k] = [h, h_\Gamma, k] \in \mathfrak{X} \times \mathscr{H} $ with $\bm{h} = [h, h_\Gamma]$.}
    \end{align}
    In the context, $ [\bm{\eta}, \theta] = [\eta, \eta_\Gamma, \theta] \in \mathfrak{X} \times \mathscr{H} $ with $ \bm{\eta} = [\eta, \eta_\Gamma]$ is the solution to the state-system (S)$_\varepsilon$, for the initial triplet $[\bm{\eta}_0, \theta_0] = [\eta_0, \eta_{\Gamma, 0}, \theta_0] \in \mathbb{W} \times V_0 $ with $\bm{\eta}_0 = [\eta_0, \eta_{\Gamma, 0}]$ and forcing triplet $ [\bm{u}, v] = [u, u_\Gamma, v] $, and $ \bar{\mathcal{P}}_\varepsilon \in \mathscr{L}(\mathfrak{X} \times \mathscr{H}; \mathfrak{Y})$ is the restriction $ \mathcal{P}|_{\{[0, 0, 0]\} \times [\mathfrak{X} \times \mathscr{H}]} $ of the bounded linear operator $ \mathcal{P} = \mathcal{P}(a, b, \mu, \omega, A) : [\mathbb{W} \times H] \times [\mathfrak{X} \times \mathscr{V}_0^* ] \longrightarrow \mathfrak{Y} $, as in Remark \ref{appendix}, in the case when:
    \begin{equation}\label{set01}
        \begin{cases}
            [a, b]= [\alpha_0, 0] \mbox{ in $ W^{1, \infty}(Q) \times L^\infty(Q) $,}
            \\
            \mu = \bar{\mu}_\varepsilon := g'(\eta) + \alpha''(\eta) f_\varepsilon(\partial_x \theta) \mbox{ in $ L^\infty(0, T; H) $,} 
            \\
            [\omega, A] = [\bar{\omega}_\varepsilon, \bar{A}_\varepsilon] := \bigl[\alpha'(\eta) f_\varepsilon'(\partial_x \theta), \alpha(\eta) f_\varepsilon''(\partial_x \theta) \bigr] \mbox{ in $ [L^\infty(Q)]^2 $.}
            \end{cases}
    \end{equation}
\end{lemma}
\paragraph{\textbf{Proof.}}{Let us fix any $ [\bm{u}, v] = [u, u_\Gamma, v] \in \mathfrak{X} \times \mathscr{H} $ with $\bm{u} = [u, u_\Gamma] $, and take any $ \delta \in (-1, 1) \setminus \{0\} $ and any $ [\bm{h}, k] = [h, h_\Gamma, k] \in \mathfrak{X} \times \mathscr{H}$ with $\bm{h} = [h, h_\Gamma]$. Then, it is easily seen that:
\begin{align}\label{GD02}
    \frac{1}{\delta} \bigl(  \mathcal{J}_\varepsilon &( \bm{u} +\delta \bm{h}, v +\delta k) -\mathcal{J}_\varepsilon(\bm{u}, v) \bigr) 
    \nonumber
    \\
    = & \frac{1}{\delta} \bigl(  \mathcal{J}_\varepsilon ( u +\delta h, u_\Gamma + \delta h_\Gamma, v +\delta k) -\mathcal{J}_\varepsilon(u, u_\Gamma, v) \bigr)
    \nonumber
    \\
    = & \left( \frac{K}{2} (\eta^\delta +\eta -2 \eta_\mathrm{ad}), \chi^\delta \right)_{\hspace{-0.5ex}\mathscr{H}} + \left( \frac{K_\Gamma}{2} (\eta^\delta_\Gamma +\eta_\Gamma -2 \eta_{\Gamma, \mathrm{ad}}), \chi^\delta_\Gamma \right)_{\hspace{-0.5ex}\mathscr{H}_\Gamma}
    \\ 
     & +\left( \frac{\Lambda}{2} (\theta^\delta +\theta -2 \theta_\mathrm{ad}), \gamma^\delta \right)_{\hspace{-0.5ex}\mathscr{H}} + \left( \frac{L}{2}(2u +\delta h), h \right)_{\hspace{-0.5ex}\mathscr{H}} \nonumber
     \\
     & +\left( \frac{L_\Gamma}{2}(2u_\Gamma +\delta h_\Gamma), h_\Gamma \right)_{\hspace{-0.5ex}\mathscr{H}_\Gamma} +\left( \frac{M}{2} (2v +\delta k), k \right)_{\hspace{-0.5ex}\mathscr{H}}.
    \nonumber
\end{align}
Here, let us set:
\begin{subequations}\label{sets02}
\begin{equation}\label{set02}
    \begin{cases}
        \displaystyle \bar{\mu}_\varepsilon^\delta := \int_0^1 g'(\eta +\varsigma \delta \chi^\delta) \, d\varsigma + f_\varepsilon(\partial_x \theta) \int_0^1 \alpha''(\eta +\varsigma \delta \chi^\delta) \, d \varsigma \mbox{ in $ L^\infty(0, T; H) $,} 
        \\
        \displaystyle \bar{\omega}_\varepsilon^\delta := \alpha'(\eta^\delta) \int_0^1 f_\varepsilon'(\partial_x \theta +\varsigma \delta \partial_x \gamma^\delta) \, d \varsigma \mbox{ in $ L^\infty(Q) $,}
        \\
        \displaystyle \bar{A}_\varepsilon^\delta := \alpha(\eta^\delta) \int_0^1 f_\varepsilon''(\partial_x \theta +\varsigma \delta \partial_x \gamma^\delta) \, d \varsigma \mbox{ in $ L^\infty(Q) $,}
    \end{cases} 
\end{equation}
and
\begin{align}\label{set03}
    \displaystyle \bar{k}_\varepsilon^\delta := M k + \partial_x & \left[ \rule{-1pt}{16pt} \right. 
    \chi^\delta f_\varepsilon'(\partial_x \theta) \int_0^1 \alpha'(\eta +\varsigma \delta \chi^\delta) \, d \varsigma 
    \nonumber
    \\
    & \qquad \displaystyle  -\chi^\delta \alpha'(\eta^\delta) \int_0^1 f_\varepsilon'(\partial_x \theta +\varsigma \delta \partial_x \gamma^\delta) \, d \varsigma 
    \left. \rule{-1pt}{16pt} \right] \mbox{ in $ \mathscr{V}_0^* $,}
    \\
    & \mbox{for all $ \delta \in (-1, 1) \setminus \{0\} $.}
    \nonumber
\end{align}
\end{subequations}
Then, in the light of Remark \ref{Rem.GD01}, one can say that:
\begin{equation*}
    [\bm{\chi}^\delta, \gamma^\delta] = [\chi^\delta, \chi^\delta_\Gamma, \gamma^\delta ] = \bar{\mathcal{P}}_\varepsilon^\delta [L h, L_\Gamma h_\Gamma, \bar{k}_\varepsilon^\delta] \mbox{ in $ \mathfrak{Y} $, for $ \delta \in (-1, 1) \setminus \{0\} $,}
\end{equation*}
by using the restriction $ \bar{\mathcal{P}}_\varepsilon^\delta := \mathcal{P}|_{\{[0, 0, 0]\} \times [\mathfrak{X} \times \mathscr{V}_0^*]} : \mathfrak{X} \times \mathscr{V}_0^* \longrightarrow \mathfrak{Y} $ of the bounded linear operator $ \mathcal{P} = \mathcal{P}(a, b, \mu, \omega, A) : [\mathbb{W} \times H] \times [\mathfrak{X} \times \mathscr{V}_0^*] \longrightarrow \mathfrak{Y} $, as in Remark \ref{appendix}, in the case when:
\begin{equation*}
    \begin{cases}
        [a, b, \omega, A]= [\alpha_0, 0, \bar{\omega}_\varepsilon^\delta, \bar{A}_\varepsilon^\delta] \mbox{ in $ W^{1, \infty}(Q) \times [L^\infty(Q)]^3 $,}
        \\
        \mu = \bar{\mu}_\varepsilon^\delta \mbox{ in $ L^\infty(0, T; H) $, for $ \delta \in (-1, 1) \setminus \{0\} $.} 
    \end{cases}
\end{equation*}
Besides, taking into account \eqref{f_eps}, \eqref{sets02}, (\hyperlink{A0l}{A0})--(\hyperlink{A3l}{A3}), 
we have:
\begin{subequations}\label{est3-01}
    \begin{align}\label{est3-01a}
    \bar{C}_0^* := & \frac{16}{\min\{ 1, \nu^2, \delta_* \}} \bigl( 1 +|\alpha_0|_{W^{1, \infty}(Q)} +2|g'|_{L^\infty(\mathbb{R})}^2
    \nonumber
    \\
    & \qquad \qquad + 2|\alpha''|_{L^\infty(\mathbb{R})}^2|f_\varepsilon(\partial_x \theta)|_{L^\infty(0, T; H)}^2  +|\alpha'|_{L^\infty(\mathbb{R})}^2 \bigr)
        \\
        \geq & \frac{16}{\min\{ 1, \nu^2, \delta_* \}} \sup_{0 < |\delta| < 1} \left\{ 1 +|\alpha_0|_{W^{1, \infty}(Q)} +|\bar{\mu}_\varepsilon^\delta|_{L^\infty(0, T; H)}^2 +|\bar{\omega}_\varepsilon^\delta|_{L^\infty(Q)}^2 \right\},
        \nonumber
    \end{align}
and
    \begin{align}\label{est3-01b}
    \bigl| \bigl< [L h(t)&\, , L_\Gamma h_\Gamma(t), \bar{k}_\varepsilon^\delta(t)], [\varphi, \varphi_\Gamma, \psi] \bigr>_{\mathbb{X} \times V_0} \bigr|
    \nonumber
    \\
     & \leq \bigl|\bigl( L h(t), \varphi \bigr)_{H}\bigr| + \bigl|\bigl(L_\Gamma h_\Gamma(t), \varphi_\Gamma \bigr)_{H_\Gamma}\bigr| +|\langle \bar{k}_\varepsilon^\delta(t), \psi \rangle_{V_0}| 
    \nonumber
        \\
    & \leq L|h(t)|_H |\varphi|_H  + L_\Gamma|h_\Gamma(t)|_{H_\Gamma}|\varphi_\Gamma|_{H_\Gamma}
    \nonumber
    \\
    & \qquad +M| k(t)|_H |\psi|_H +2|\alpha'|_{L^\infty(\mathbb{R})}|\chi^\delta(t)|_H|\partial_x \psi|_H
    \nonumber
        \\
       & \leq L |h(t)|_H |\varphi|_V + L_\Gamma |h_\Gamma(t)|_{H_\Gamma}|\varphi_\Gamma|_{H_\Gamma}
        \nonumber
        \\
        & \qquad +\bigl( \sqrt{2} M|k(t)|_H +2|\alpha'|_{L^\infty(\mathbb{R})}|\chi^\delta(t)|_H \bigr) |\psi|_{V_0},
        \\
    & \mbox{ for a.e. $ t \in (0, T) $, any $ [\bm{\varphi}, \psi] = [\varphi, \varphi_\Gamma, \psi] \in \mathbb{X} \times V_0 $ }
    \nonumber
    \\
    & \qquad \mbox{with $\bm{\varphi} = [\varphi, \varphi_\Gamma] $, and any $  \delta \in (-1, 1) \setminus \{0\} $,}
    \nonumber
\end{align} 
so that 
    \begin{align}\label{est3-01c}
    \bigl| [L h(t), L_\Gamma & h_\Gamma(t), \bar{k}_\varepsilon^\delta(t)] \bigr|_{\mathbb{X} \times V_0^*}^2 \leq \bar{B}_0^* \bigl( \bigl| [\bm{h}(t), k(t)] \bigr|_{\mathbb{X} \times H}^2 +|\chi^\delta(t)|_H^2 \bigr),
    \nonumber
        \\
    & \mbox{for a.e. $ t \in (0, T) $, and any $ \delta \in (-1, 1) \setminus \{ 0 \} $,}
\end{align}
\end{subequations}
with a positive constant $ \bar{B}_0^* := 4 \bigl(  L^2+ L_\Gamma^2 + M^2 + |\alpha'|_{L^\infty(\mathbb{R})}^2  \bigr) $.
\medskip

Now, having in mind \eqref{est3-01}, let us apply Theorem \ref{Prop02} (I) to the case when:
\begin{equation*}
    \begin{array}{c}
    \begin{cases}
        [a, b, \mu, \omega, A] = [\alpha_0, 0, \overline{\mu}_\varepsilon^\delta, \bar{\omega}_\varepsilon^\delta, \bar{A}_\varepsilon^\delta],
        \\[0.5ex]
        [\bm{h}, k] = [h, h_\Gamma, k] = [L h, L_\Gamma h_\Gamma, \bar{k}_\varepsilon^\delta],
        \\[0.5ex]
        [\bm{p}, z] = [p, p_\Gamma, z] = [\bm{\chi}^\delta, \gamma^\delta ] = [\chi^\delta, \chi^\delta_\Gamma, \gamma^\delta] = \bar{\mathcal{P}}_\varepsilon^\delta [L h, L_\Gamma h_\Gamma, \bar{k}_\varepsilon^\delta],
    \end{cases}
        \ \\
        \ \\[-1.5ex]
    \quad\mbox{for $ \delta \in (-1, 1) \setminus \{0\} $.}
    \end{array}
\end{equation*}
Then, we estimate that:
\begin{align*}
    \frac{d}{dt} & \bigl( |\bm{\chi}^\delta(t)|_{\mathbb{X}}^2 +|\sqrt{\alpha_0(t)} \gamma^\delta(t)|_H^2 \bigr) +\bigl( |\bm{\chi}^\delta(t)|_{\mathbb{W}}^2 +\nu^2 |\gamma^\delta(t)|_{V_0}^2 \bigr)
    \\
    & \leq  \bar{C}_0^* \bigl( |\bm{\chi}^\delta(t)|_\mathbb{X}^2 +|\sqrt{\alpha_0(t)} \gamma^\delta(t)|_H^2 \bigr) +\bar{C}_0^* \big( |L h(t)|_{V^*}^2 + |L_\Gamma h_\Gamma(t)|^2_{H_\Gamma^*} +|\bar{k}_\varepsilon^\delta(t)|_{V_0^*}^2 \bigr)
    \\
    & \leq  \bar{C}_0^* (1 +\bar{B}_0^*) \bigl( |\bm{\chi}^\delta(t)|_\mathbb{X}^2 +|\sqrt{\alpha_0(t)} \gamma^\delta(t)|_H^2 \bigr) +\bar{C}_0^* \bar{B}_0^* \big( |\bm{h}(t)|_{\mathbb{X}}^2 +|k(t)|_{H}^2 \bigr), 
    \\
    & \qquad\mbox{for a.e. $ t \in (0, T) $,}
\end{align*}
and subsequently, by using (\hyperlink{A1l}{A1}) and Gronwall's lemma, we observe that:
\begin{itemize}
    \item[\textmd{$(\hypertarget{star1}{\star\,1})$}]the sequence $ \{ [\bm{\chi}^\delta, \gamma^\delta] \}_{\delta \in (-1, 1) \setminus \{0\}} = \{ [\chi^\delta, \chi^\delta_\Gamma, \gamma^\delta] \}_{\delta \in (-1, 1) \setminus \{0\}}$ is bounded in \linebreak $ [C([0, T]; \mathbb{X})$ $ \times C([0, T]; H)] \cap [\mathfrak{W} \times \mathscr{V}_0] $. 
\end{itemize}

Meanwhile, as consequences of \eqref{set00}, \eqref{set01}--\eqref{est3-01}, $(\hyperlink{star1}{\star\,1})$, (\hyperlink{A0l}{A0})--(\hyperlink{A3l}{A3}), Main Theorem \ref{mainTh01} \linebreak (I-B), Remark \ref{Rem.mTh01Conv}, and Lebesgue's dominated convergence theorem, one can find a sequence $ \{ \delta_n \}_{n = 1}^\infty \subset \mathbb{R} $, such that:
\begin{subequations}\label{convs3-01}
    \begin{align}\label{convs3-00a}
        0 < |\delta_n| < 1, \mbox{ and } \delta_n \to 0, \mbox{ as $ n \to \infty $,}
    \end{align}
    \begin{align}\label{convs3-01a}
        &
        \begin{cases}
            [\delta_n \bm{\chi}^{\delta_n}, \delta_n \gamma^{\delta_n}] = [\bm{\eta}^{\delta_n} -\eta, \theta^{\delta_n} -\theta] \to [0, 0, 0]
            \\
            \qquad  \mbox{ in $ [C(\overline{Q}) \times C(\overline{\Sigma})]\times C(\overline{Q}) $, and in $ \mathfrak{W} \times \mathscr{V}_0 $,}
            \\[1ex]
            [\delta_n \partial_x \chi^{\delta_n}, \delta_n \partial_x \gamma^{\delta_n}] = [\partial_x (\eta^{\delta_n} -\eta), \partial_x (\theta^{\delta_n} -\theta)] \to [0, 0]
            \\
            \qquad \mbox{ in $ [\mathscr{H}]^2 $, and in the pointwise sense a.e. in $ Q $,}
        \end{cases}
        \mbox{as $ n \to \infty $,}
    \end{align}
    \begin{align}\label{convs3-01b}
        [\bar{\omega}_\varepsilon^{\delta_n}, & \bar{A}_\varepsilon^{\delta_n}] \to [\bar{\omega}_\varepsilon, \bar{A}_\varepsilon]  \mbox{ weakly-$ * $ in $ [L^\infty(Q)]^2 $,}
        \nonumber
        \\
        & \mbox{and in the pointwise sense a.e. in $ Q $, ~ as $ n \to \infty $,}
    \end{align}
    \begin{equation}\label{convs3-01c}
        \begin{cases}
            \bar{\mu}_\varepsilon^{\delta_n} \to \bar{\mu}_\varepsilon \mbox{ weakly-$ * $ in $ L^\infty(0, T; H) $,}
            \\
            \bar{\mu}_\varepsilon^{\delta_n}(t) \to \bar{\mu}_\varepsilon(t) \mbox{ in $ H $,}
            \mbox{ for a.e. $ t \in (0, T) $,}
        \end{cases}
        \mbox{as $ n \to \infty $,}
    \end{equation}
    and
    \begin{align}\label{convs3-01d}
        \langle \bar{k}_\varepsilon^{\delta_n} - &\, M k, \psi \rangle_{\mathscr{V}_0} = -\left( \chi^{\delta_n}, ~ f_\varepsilon'(\partial_x \theta) \left( \rule{-1pt}{14pt} \right. \int_0^1 \alpha'(\eta +\varsigma \delta_n \chi^{\delta_n}) \, d \varsigma  \left. \rule{-1pt}{14pt} \right) \partial_x \psi\right)_{\hspace{-0.5ex}\mathscr{H}}
        \nonumber
        \\
        &~ +\left( \chi^{\delta_n}, ~ 
        \alpha'(\eta^{\delta_n}) \left( \rule{-1pt}{14pt} \right. \int_0^1 f_\varepsilon'(\partial_x \theta +\varsigma \delta_n \partial_x \gamma^{\delta_n}) \, d \varsigma \left. \rule{-1pt}{14pt} \right) \partial_x \psi \right)_{\hspace{-0.5ex}\mathscr{H}} \to  0, \mbox{ as $ n  \to \infty $.}
    \end{align}
\end{subequations}

On account of \eqref{set00}, \eqref{set01}--\eqref{convs3-01}, and Remark \ref{appendix}, we can apply Theorem \ref{Prop03}, and can see that:
\begin{align}\label{conv3-02}
    [\bm{\chi}^{\delta_n}, & \gamma^{\delta_n}] = [\chi^{\delta_n}, \chi^{\delta_n}_\Gamma, \gamma^{\delta_n}]= \bar{\mathcal{P}}_\varepsilon^{\delta_n}[L h, L_\Gamma h_\Gamma, \bar{k}_\varepsilon^{\delta_n}]\nonumber
    \\
    & \to [\bm{\chi}, \gamma] = [\chi, \chi_\Gamma, \gamma ] := \bar{\mathcal{P}}_\varepsilon[L h, L_\Gamma h_\Gamma, M k] 
    \nonumber
    \\
    & \mbox{ in $ [C(\overline{Q}) \times C(\overline{\Sigma})] \times \mathscr{H} $, as $ n \to \infty $.}
\end{align}
Since the uniqueness of the solution $ [\bm{\chi}, \gamma] = [\chi, \chi_\Gamma, \gamma ] = \bar{\mathcal{P}}_\varepsilon[L h, L_\Gamma h_\Gamma, M k] $ is guaranteed in Theorem \ref{Prop01}, the observations \eqref{GD02}, \eqref{convs3-01}, and \eqref{conv3-02} enable us to compute the directional derivative $ D_{[\bm{h}, k]}\mathcal{J}_\varepsilon(\bm{u}, v) = D_{[h, h_\Gamma, k]}\mathcal{J}_\varepsilon(u, u_\Gamma, v) \in \mathbb{R} $, as follows:
\begin{align*}
    D_{[\bm{h}, k]} & \mathcal{J}_\varepsilon(\bm{u}, v) = D_{[h, h_\Gamma, k]} \mathcal{J}_\varepsilon(u, u_\Gamma, v):= \lim_{\delta \to 0} \frac{1}{\delta} \bigl( \mathcal{J}_\varepsilon(\bm{u} +\delta \bm{h}, v +\delta k) -\mathcal{J}_\varepsilon(\bm{u}, v) \bigr) 
    \\
    & = \lim_{\delta \to 0} \frac{1}{\delta} \bigl( \mathcal{J}_\varepsilon(u +\delta h, u_\Gamma + \delta h_\Gamma, v +\delta k) -\mathcal{J}_\varepsilon(u, u_\Gamma, v) \bigr) 
    \\
    & = \bigl( [K (\eta -\eta_\mathrm{ad}), K_\Gamma (\eta_\Gamma - \eta_{\Gamma, \mathrm{ad}}), \Lambda (\theta -\theta_\mathrm{ad})], \bar{\mathcal{P}}_\varepsilon [L h, L_\Gamma h_\Gamma, M k] \bigr)_{\mathfrak{X} \times \mathscr{H}}
    \\
    & \qquad 
    +\bigl( [L u, L_\Gamma u_\Gamma, M v], [h, h_\Gamma, k] \bigr)_{\mathfrak{X} \times \mathscr{H}}, 
    \\
    & \mbox{ for any $ [\bm{u}, v] = [u, u_\Gamma, v] \in \mathfrak{X} \times \mathscr{H} $ with $\bm{u} = [u, u_\Gamma]$},
    \\
    & \mbox{ and any direction $ [\bm{h}, k] = [h, h_\Gamma, k] \in \mathfrak{X} \times \mathscr{H} $ with $\bm{h} = [h, h_\Gamma] $.}
\end{align*}
Moreover, with Remark \ref{appendix} and Riesz's theorem in mind, we deduce the existence of the G\^{a}teaux derivative $ \mathcal{J}_\varepsilon'(\bm{u}, v) = \mathcal{J}_\varepsilon'(u, u_\Gamma, v) \in [\mathfrak{X} \times \mathscr{H}]^* $ $ (= \mathfrak{X} \times \mathscr{H}) $ at $ [\bm{u}, v] = [u, u_\Gamma, v] \in \mathfrak{X} \times \mathscr{H} $ with $\bm{u} = [u, u_\Gamma]$, i.e.:
\begin{align*}
    \bigl( &\,\mathcal{J}_\varepsilon'(u, u_\Gamma, v), [h, h_\Gamma, k] \bigr)_{\mathfrak{X} \times \mathscr{H}} =  D_{[h, h_\Gamma, k]} \mathcal{J}_\varepsilon(u, u_\Gamma, v),
    \\
    & \mbox{for every $ [\bm{u}, v] = [u, u_\Gamma, v] $ and $ [\bm{h}, k] = [h, h_\Gamma, k] \in \mathfrak{X} \times \mathscr{H}, $}
    \\
    & \mbox{with $\bm{u} =[u, u_\Gamma]$ and $ \bm{h} = [h, h_\Gamma]$, respectively.}
\end{align*}
Thus, we conclude this lemma with the required property \eqref{GD01}. 
\hfill \qed
}
\begin{lemma}\label{Lem.GD02}
    Under the assumptions (\hyperlink{A0l}{A0})--(\hyperlink{A3l}{A3}), let $ [\bm{u}_\varepsilon^*, v_\varepsilon^*] = [u_\varepsilon^*, u^*_{\Gamma, \varepsilon}, v^*_\varepsilon] \in \mathfrak{X} \times \mathscr{H} $ with $\bm{u}^*_\varepsilon = [u_\varepsilon^*, u^*_{\Gamma, \varepsilon}] $ be an optimal control of the problem (OP)$_\varepsilon$, and let $ [\bm{\eta}_\varepsilon^*, \theta_\varepsilon^*] = [\eta_\varepsilon^*, \eta_{\Gamma, \varepsilon}^*, \theta_\varepsilon^*] \in \mathfrak{X} \times \mathscr{H}$ with $\bm{\eta}_\varepsilon^* = [\eta_\varepsilon^*, \eta_{\Gamma, \varepsilon}^*] $ be the solution to the system (S)$_\varepsilon$, for the initial triplet $[\bm{\eta}_0, \theta_0] = [\eta_0, \eta_{\Gamma, 0}, \theta_0] \in \mathbb{W} \times V_0 $ with $\bm{\eta}_0 = [\eta_0, \eta_{\Gamma, 0}] $ and forcing triplet $ [\bm{u}_\varepsilon^*, v_\varepsilon^*] = [u_\varepsilon^*, u_{\Gamma, \varepsilon}^*, v_\varepsilon^*]$. Also, let $ \mathcal{P}_\varepsilon^* : \mathfrak{X} \times \mathscr{H} \longrightarrow \mathfrak{Y} $ be the bounded linear operator, defined in Remark \ref{Rem.mTh03}, with use of the solution $ [\bm{\eta}_\varepsilon^*, \theta_\varepsilon^*] = [\eta_\varepsilon^*, \eta_{\Gamma, \varepsilon}^*, \theta_\varepsilon^*]$. Let $ \mathcal{P}_\varepsilon \in \mathscr{L}( \mathfrak{X} \times \mathscr{H}; \mathfrak{Y})$ be the restriction $ \mathcal{P}|_{\{[0, 0, 0]\} \times [\mathfrak{X} \times \mathscr{H}]} $ of the bounded linear operator $ \mathcal{P} = \mathcal{P}(a, b, \mu, \omega, A) : [\mathbb{W} \times H] \times [\mathfrak{X} \times \mathscr{V}_0^*] \longrightarrow \mathfrak{Y} $, as in Remark \ref{appendix}, in the case when:
    \begin{equation}\label{set_GD02-01}
        \begin{cases}
            [a, b]= [\alpha_0, 0] \mbox{ in $ W^{1, \infty}(Q) \times L^\infty(Q) $,}
            \\
            \mu = g'(\eta_\varepsilon^*) + \alpha''(\eta_\varepsilon^*) f_\varepsilon(\partial_x \theta_\varepsilon^*) \mbox{ in $ L^\infty(0, T; H) $,} 
            \\
            [\omega, A] = \bigl[\alpha'(\eta_\varepsilon^*) f_\varepsilon'(\partial_x \theta_\varepsilon^*), \alpha(\eta_\varepsilon^*) f_\varepsilon''(\partial_x \theta_\varepsilon^*) \bigr] \mbox{ in $ [L^\infty(Q)]^2 $.}
            \end{cases}
    \end{equation}
    Then, the operators $ \mathcal{P}_\varepsilon^* $ and $ \mathcal{P}_\varepsilon $  have a conjugate relationship, in the following sense:
    \begin{align*}
        \bigl( \mathcal{P}_\varepsilon^*&\, [\bm{u}, v] [\bm{h}, k] \bigr)_{\mathfrak{X} \times \mathscr{H}} = \bigl( \mathcal{P}_\varepsilon^*[u, u_\Gamma, v] , [h, h_\Gamma, k] \bigr)_{\mathfrak{X} \times \mathscr{H}}
        \\
        & = \bigl( [\bm{u}, v], \mathcal{P}_\varepsilon [\bm{h}, k] \bigr)_{\mathfrak{X} \times \mathscr{H}} = \bigl( [u, u_\Gamma, v], \mathcal{P}_\varepsilon [h, h_\Gamma, k] \bigr)_{\mathfrak{X} \times \mathscr{H}}, 
        \\
        & \mbox{ for all $ [\bm{h}, k] = [h, h_\Gamma, k] $ and $ [\bm{u}, v] = [u, u_\Gamma, v] \in \mathfrak{X} \times \mathscr{H} $}
        \\
        & \mbox{ with $\bm{h} = [h, h_\Gamma] $ and $ \bm{u} = [u, u_\Gamma] $, respectively.}
    \end{align*}
\end{lemma}
\paragraph{\textbf{Proof.}}{
    Let us fix arbitrary triplets of functions $ [\bm{h}, k] = [h, h_\Gamma, k] \in \mathfrak{X} \times \mathscr{H}$ with $\bm{h} = [h, h_\Gamma]$ and $ [\bm{u}, v] = [u, u_\Gamma, v] \in \mathfrak{X} \times \mathscr{H} $ with $\bm{u} = [u, u_\Gamma] $, and let us put:
    \begin{align*}
    & [\bm{\chi}_\varepsilon, \gamma_\varepsilon] = [\chi_\varepsilon, \chi_{\Gamma, \varepsilon}, \gamma_\varepsilon]:= \mathcal{P}_\varepsilon [\bm{h}, k] = \mathcal{P}_\varepsilon [h, h_\Gamma, k] \mbox{ and }
    \\
    & [\bm{p}_\varepsilon, z_\varepsilon] = [p_\varepsilon, p_{\Gamma, \varepsilon}, z_\varepsilon]:= \mathcal{P}_\varepsilon^* [\bm{u}, v] = \mathcal{P}_\varepsilon^* [u, u_\Gamma, v], \mbox{ in $ \mathfrak{X} \times \mathscr{H} $.} 
\end{align*}
Then, invoking Theorem \ref{Prop01}, and the settings as in \eqref{setRem4} and \eqref{set_GD02-01}, we compute that:
\begin{align*}
    \bigl( \mathcal{P}_\varepsilon^* & [\bm{u}, v], [\bm{h}, k] \bigr)_{\mathfrak{X} \times \mathscr{H}} = \int_0^T \bigl( \bm{p}_\varepsilon(t), \bm{h}(t) \bigr)_{\mathbb{X}} \, dt  +\int_0^T \bigl( z_\varepsilon(t), k(t) \bigr)_{H} \, dt
    \\
    = & \int_0^T \bigl( \bm{h}(t), \bm{p}_\varepsilon(t) \bigr)_{\mathbb{X}} \, dt  +\int_0^T \bigl( k(t), z_\varepsilon(t) \bigr)_{H} \, dt
    \\
    = & \int_0^T \left[ \rule{-1pt}{16pt} \right. \bigl( \partial_t \bm{\chi}_\varepsilon(t), \bm{p}_\varepsilon(t) \bigr)_{\mathbb{X}} +\bigl( \partial_x \chi_\varepsilon(t), \partial_x p_\varepsilon(t) \bigr)_H 
    \\
    & \qquad +\bigl( \alpha''(\eta_\varepsilon^*(t)) f_\varepsilon(\partial_x \theta_\varepsilon^*(t)) \chi_\varepsilon(t),  p_\varepsilon(t) \bigr)_H 
    \\
    & \qquad +\bigl( g'(\eta_\varepsilon^*(t)) \chi_\varepsilon(t), p_\varepsilon(t) \bigr)_H +\bigl( \alpha'(\eta_\varepsilon^*(t))f_\varepsilon'(\partial_x \theta_\varepsilon^*(t)) \partial_x \gamma_\varepsilon(t), p_\varepsilon(t) \bigr)_{H} \left. \rule{-1pt}{16pt} \right] \, dt
    \\
    & +\int_0^T \left[ \rule{-1pt}{16pt} \right. \bigl< \alpha_0(t) \partial_t \gamma_\varepsilon(t), z_\varepsilon(t) \bigr>_{V_0} +\bigl( \alpha'(\eta_\varepsilon^*(t)) f_\varepsilon'(\partial_x \theta_\varepsilon^*(t)) \chi_\varepsilon(t), \partial_x z_\varepsilon(t) \bigr)_{H} 
   \\
    & \qquad +\bigl( \alpha(\eta_\varepsilon^*(t)) f_\varepsilon''(\partial_x \theta_\varepsilon^*(t)) \partial_x \gamma_\varepsilon(t), \partial_x z_\varepsilon(t) \bigr)_H +\nu^2 \bigl( \partial_x \gamma_\varepsilon(t), \partial_x z_\varepsilon(t) \bigr)_H \left. \rule{-1pt}{16pt} \right] \, dt
\end{align*}
\begin{align*}
    = & \bigl( \bm{p}_\varepsilon(T), \bm{\chi}_\varepsilon(T) \bigr)_{\mathbb{X}} -\bigl( \bm{p}_\varepsilon(0), \bm{\chi}_\varepsilon(0) \bigr)_{\mathbb{X}} +\int_0^T \left[ \rule{-1pt}{16pt} \right. \bigl( -\partial_t \bm{p}_\varepsilon(t), \bm{\chi}_\varepsilon(t) \bigr)_{\mathbb{X}} 
    \\
    & \qquad +\bigl( \partial_x p_\varepsilon(t), \partial_x \chi_\varepsilon(t) \bigr)_H +\bigl( \alpha''(\eta_\varepsilon^*(t)) f_\varepsilon(\partial_x \theta_\varepsilon^*(t)) p_\varepsilon(t), \chi_\varepsilon(t) \bigr)_H 
    \\
& \qquad +\bigl( g'(\eta_\varepsilon^*(t)) p_\varepsilon(t), \chi_\varepsilon(t) \bigr)_H +\bigl( \alpha'(\eta_\varepsilon^*(t))f_\varepsilon'(\partial_x \theta_\varepsilon^*(t)) \partial_x z_\varepsilon(t), \chi_\varepsilon(t) \bigr)_{H} \left. \rule{-1pt}{16pt} \right] \, dt
    \\
    & +\bigl( \alpha_0(T) z_\varepsilon(T), \gamma_\varepsilon(T) \bigr)_H -\bigl( \alpha_0(0) z_\varepsilon(0), \gamma_\varepsilon(0) \bigr)_H 
    \\
    & \qquad +\int_0^T \left[ \rule{-1pt}{16pt} \right. \bigl< -\partial_t \bigl(\alpha_0 z_\varepsilon)(t),  \gamma_\varepsilon(t) \bigr>_{V_0} +\bigl( \alpha'(\eta_\varepsilon^*(t)) f_\varepsilon'(\partial_x \theta_\varepsilon^*(t)) p_\varepsilon(t), \partial_x \gamma_\varepsilon(t) \bigr)_{H} 
    \\
    & \qquad +\bigl( \alpha(\eta_\varepsilon^*(t)) f_\varepsilon''(\partial_x \theta_\varepsilon^*(t)) \partial_x z_\varepsilon(t), \partial_x \gamma_\varepsilon(t) \bigr)_H +\nu^2 \bigl( \partial_x z_\varepsilon(t), \partial_x \gamma_\varepsilon(t) \bigr)_H \left. \rule{-1pt}{16pt} \right] \, dt
    \\
    = & ( \bm{u}, \bm{\chi}_\varepsilon )_{\mathfrak{X}}  + ( v, \gamma_\varepsilon )_{\mathscr{H}} = \bigl( [\bm{u}, v], \mathcal{P}_\varepsilon [\bm{h}, k]  \bigr)_{\mathfrak{X} \times \mathscr{H}}.
\end{align*}
\qed
}
\begin{remark}\label{Rem.GD02}
    Note that the operator $ \mathcal{P}_\varepsilon \in \mathscr{L}(\mathfrak{X} \times \mathscr{H}; \mathfrak{Y}) $, as in Lemma \ref{Lem.GD02}, corresponds to the operator $ \bar{\mathcal{P}}_\varepsilon \in \mathscr{L}(\mathfrak{X} \times \mathscr{H}; \mathfrak{Y}) $, as in the previous Lemma \ref{Lem.GD01}, under the special setting \eqref{set_GD02-01}. 
\end{remark}
\bigskip

Now, we are ready to prove the Main Theorem \ref{mainTh03}.

\paragraph{\textbf{Proof of (III-A) of Main Theorem \ref{mainTh03}.}}{
    Let $ [\bm{u}_\varepsilon^*, v_\varepsilon^*] = [u_\varepsilon^*, u_{\Gamma, \varepsilon}^*, v_\varepsilon^*] \in \mathfrak{X} \times \mathscr{H} $ with $\bm{u}_\varepsilon^* = [u_\varepsilon^*, u_{\Gamma, \varepsilon}^*] $ be the optimal control of (OP)$_\varepsilon$, let $ [\bm{\eta}_\varepsilon^*, \theta_\varepsilon^*] = [\eta_\varepsilon^*, \eta_{\Gamma, \varepsilon}^*, \theta_\varepsilon^*] \in \mathfrak{X} \times \mathscr{H} $ with $\bm{\eta}_\varepsilon^* = [\eta_\varepsilon^*, \eta_{\Gamma, \varepsilon}^*] $ be the solution to the system (S)$_\varepsilon$ for the initial triplet $[\bm{\eta}_0, \theta_0] = [\eta_0, \eta_{\Gamma, 0}, \theta_0] \in \mathbb{W} \times V_0$ with $\bm{\eta}_0 = [\eta_0, \eta_{\Gamma, 0}] $, and forcing triplet $ [\bm{u}_\varepsilon^*, v_\varepsilon^*] = [u_\varepsilon^*, u_{\Gamma, \varepsilon}^*, v_\varepsilon^*]$, and let $ \mathcal{P}_\varepsilon \in \mathscr{L}(\mathfrak{X} \times \mathscr{H}; \mathfrak{Y}) $ and $ \mathcal{P}_\varepsilon^* \in \mathscr{L}(\mathfrak{X} \times \mathscr{H}; \mathfrak{Y}) $ be the two operators as in Lemma \ref{Lem.GD02}. Then, on the basis of the previous Lemmas \ref{Lem.GD01} and \ref{Lem.GD02}, Main Theorem \ref{mainTh03} (III-A) will be demonstrated as follows:
    \begin{align*}
        0 = &~  \left(\mathcal{J}_\varepsilon'(\bm{u}_\varepsilon^*, v_\varepsilon^*), [\bm{h}, k]\right)_{\mathfrak{X} \times \mathscr{H}} = \left(\mathcal{J}_\varepsilon'(u_\varepsilon^*, u_{\Gamma, \varepsilon}^*, v_\varepsilon^*), [h, h_\Gamma, k]\right)_{\mathfrak{X} \times \mathscr{H}}
        \\
        = &~ \lim_{\delta \to 0} \frac{1}{\delta} \bigl( \mathcal{J}_\varepsilon(\bm{u}_\varepsilon^* +\delta \bm{h}, v_\varepsilon^* +\delta k) -\mathcal{J}_\varepsilon(\bm{u}_\varepsilon^*, v_\varepsilon^*) \bigr) 
        \\
        = &~ \lim_{\delta \to 0} \frac{1}{\delta} \bigl( \mathcal{J}_\varepsilon(u_\varepsilon^* +\delta h, u_{\Gamma, \varepsilon}^* + \delta h_\Gamma, v_\varepsilon^* +\delta k) -\mathcal{J}_\varepsilon(u_\varepsilon^*, u_{\Gamma, \varepsilon}, v_\varepsilon^*) \bigr)
        \\
        = &~ \bigl( [K (\eta_\varepsilon^* -\eta_\mathrm{ad}), K_\Gamma (\eta_{\Gamma, \varepsilon}^* - \eta_{\Gamma, \mathrm{ad}}), \Lambda (\theta_\varepsilon^* -\theta_\mathrm{ad})], \mathcal{P}_\varepsilon [L h, L_\Gamma h_\Gamma, M k] \bigr)_{\mathfrak{X} \times \mathscr{H}}
        \\
        &~ +\bigl( [L u_\varepsilon^*, L_\Gamma u_{\Gamma, \varepsilon}^*, M v_\varepsilon^*], [h, h_\Gamma, k] \bigr)_{\mathfrak{X} \times \mathscr{H}}
        \\
        = &~ \bigl( \mathcal{P}_\varepsilon^* [K (\eta_\varepsilon^* -\eta_\mathrm{ad}), K_\Gamma(\eta_{\Gamma, \varepsilon}^* - \eta_{\Gamma, \mathrm{ad}}), \Lambda (\theta_\varepsilon^* -\theta_\mathrm{ad})], [L h, L_\Gamma h_\Gamma, M k] \bigr)_{\mathfrak{X} \times \mathscr{H}} 
        \\
        &~ +\bigl( [L u_\varepsilon^*, L_\Gamma u_{\Gamma, \varepsilon}^*, M v_\varepsilon^*], [h, h_\Gamma, k] \bigr)_{\mathfrak{X} \times \mathscr{H}}
        \\
        = &~ \bigl([L p_\varepsilon^*, L_\Gamma p_{\Gamma, \varepsilon}^*, M z_\varepsilon^*] , [h, h_\Gamma, k] \bigr)_{\mathfrak{X} \times \mathscr{H}} +\bigl( [L u_\varepsilon^*, L_\Gamma u_{\Gamma, \varepsilon}^*, M v_\varepsilon^*], [h, h_\Gamma, k] \bigr)_{\mathfrak{X} \times \mathscr{H}}
        \\
        = &~ \bigl( [L( p_\varepsilon^* +u_\varepsilon^*), L_\Gamma (p_{\Gamma, \varepsilon}^* + u_{\Gamma, \varepsilon}^*), M (z_\varepsilon^* + v_\varepsilon^*)], [h, h_\Gamma, k] \bigr)_{\mathfrak{X} \times \mathscr{H}}, 
        \\
        &~ \mbox{ for any $ [\bm{h}, k] = [h, h_\Gamma, k] \in \mathfrak{X} \times\mathscr{H} $ with $\bm{h} = [h, h_\Gamma] $.}
    \end{align*}
    \qed
}

\paragraph{\textbf{Proof of (III-B) of Main Theorem \ref{mainTh03}.}}{
    Let $[\bm{\eta_{0}}, \theta_{0}] = [\eta_0, \eta_{\Gamma, 0}, \theta_0] \in \mathbb{W} \times V_{0}$ with $\bm{\eta}_0 = [\eta_0, \eta_{\Gamma, 0}] $ be the fixed initial triplet. For any $ \varepsilon > 0 $, let $ [ \bm{u}_\varepsilon^*, v_\varepsilon^* ] = [u_\varepsilon^*, u_{\Gamma, \varepsilon}^*, v_\varepsilon^*] \in \mathfrak{X} \times \mathscr{H} $ with $\bm{u}_\varepsilon^* = [u_\varepsilon^*, u_{\Gamma, \varepsilon}^*] $, $ [\bm{\eta}_\varepsilon^*, \theta_\varepsilon^*] = [\eta_\varepsilon^*, \eta_{\Gamma, \varepsilon}^*, \theta_\varepsilon^*] \in \mathfrak{X} \times \mathscr{H} $ with $\bm{\eta}_\varepsilon^* = [\eta_\varepsilon^*, \eta_{\Gamma, \varepsilon}^*] $, and $ [\bm{p}_\varepsilon^*, z_\varepsilon^*] = [p_\varepsilon^*, p_{\Gamma, \varepsilon}^*, z_\varepsilon^*] \in  \mathfrak{Y}$ with $\bm{p}_\varepsilon^* = [p_\varepsilon^*, p_{\Gamma, \varepsilon}^*] $ be as in Main Theorem \ref{mainTh03} (III-A). Then, by Main Theorem \ref{mainTh02} (II-B), we find an optimal control $ [\bm{u}^\circ, v^\circ] = [u^\circ, u_\Gamma^\circ, v^\circ] \in \mathfrak{X} \times \mathscr{H} $ with $\bm{u}^\circ = [u^\circ, u_\Gamma^\circ]$ of (OP)$_0$, and find a zero-convergent sequence $ \{ \varepsilon_n \}_{n = 1}^\infty \subset (0, 1) $, such that:
\begin{subequations}\label{conv3B}
    \begin{align}\label{conv3B-01}
        [L u_n^*, L_\Gamma u_{\Gamma, n}^*,&\, M v_n^*] := [L u_{\varepsilon_n}^*, L_\Gamma u_{\Gamma, \varepsilon_n}^*, M v_{\varepsilon_n}^*]  \to [L u^\circ, L_\Gamma u_\Gamma^\circ, M v^\circ] 
        \nonumber 
        \\
        & \mbox{ weakly in $ \mathfrak{X} \times \mathscr{H} $, as $ n \to \infty $.}
    \end{align}
    Let $ [\bm{\eta}^\circ, \theta^\circ] \in \mathfrak{X} \times \mathscr{H} $ with $\bm{\eta}^\circ = [\eta^\circ, \eta_\Gamma^\circ] $ be the solution to (S)$_0$, for the initial triplet $[\bm{\eta}_{0}, \theta_{0}] = [\eta_0, \eta_{\Gamma, 0}, \theta_0] $ and forcing triplet $ [\bm{u}^\circ, v^\circ] = [u^\circ, u^\circ_\Gamma, v^\circ]$. Then, having in mind \eqref{conv3B-01}, Main Theorem \ref{mainTh01} (I-B), and Remark \ref{Rem.mTh01Conv}, we can find a subsequence of $ \{ \varepsilon_n \}_{n=1}^{\infty} $ (not relabeled) and a function $ \nu^\circ \in L^\infty(Q) $, such that:
\begin{align}\label{conv3B-02}
    [\bm{\eta}_n^*, \theta_n^*] &~ := [\bm{\eta}_{\varepsilon_n}^*, \theta_{\varepsilon_n}^*] \to [\bm{\eta}^\circ, \theta^\circ] \mbox{ in $ [C(\overline{Q}) \times C(\overline{\Sigma})] \times C(\overline{Q}) $, in $ \mathfrak{W} \times \mathscr{V}_0 $,}
    \nonumber
    \\
     &~ \mbox{ and weakly-$*$ in $ L^\infty(0, T; \mathbb{V}) \times L^\infty(0, T; V_0) $,}
\end{align}
\begin{align}\label{conv3B-10}
    [\partial_x \eta_n,&~  \partial_x \theta_n] \to [\partial_x \eta^\circ, \partial_x \theta^\circ] \mbox{ in $ [\mathscr{H}]^2 $,} 
    \nonumber
    \\
    &~ \mbox{and in the pointwise sense a.e. in $ Q $,}
\end{align}
\begin{align}\label{conv3B-03}
    & 
    \begin{cases}
        \mu_n^* := g'(\eta_n^*) + \alpha''(\eta_n^*) f_{\varepsilon_n}(\partial_x \theta_n^*) \to \mu^\circ := g'(\eta^\circ) + \alpha''(\eta^\circ) |\partial_x \theta^\circ| 
        \\
        \quad~ \mbox{weakly-$ * $ in $ L^\infty(0, T; H) $, and}
        \\
        \quad~ \mbox{ in the pointwise sense a.e. in $ Q $,}
        \\[1ex]
        \mu_n^*(t) \to \mu^\circ(t) \mbox{ in $ H $, in the pointwise sense for a.e. $ t \in (0, T) $,}
    \end{cases}
\end{align}
\begin{align}\label{conv3B-05}
    & 
    \begin{cases}
        f_{\varepsilon_n}'(\partial_x \theta_n^*) \to \nu^\circ \mbox{ weakly-$*$ in $ L^\infty(Q) $,}
        \\
        |\nu^\circ| \leq 1 \mbox{ a.e. in $ Q $,}
    \end{cases}
\end{align}
and
\begin{align}\label{conv3B-06}
    \omega_n^* &~ := \alpha'(\eta_n^*) f_{\varepsilon_n}'(\partial_x \theta_n^*) \to \alpha'(\eta^\circ)\nu^\circ \mbox{ weakly-$*$ in $ L^\infty(Q) $, as $ n \to \infty $.} 
\end{align}
\end{subequations}
Besides, from \eqref{conv3B-10}, \eqref{conv3B-05}, Remark \ref{Rem.MG} (\hyperlink{Fact1}{Fact\,1}) and (\hyperlink{Fact2}{Fact\,2}), and \cite[Proposition 2.16]{MR0348562}, one can see that:
\begin{equation}\label{conv3B-07}
    \nu^\circ \in \partial f_0(\partial_x \theta^\circ) = \mathrm{Sgn}^1(\partial_x \theta^\circ) \mbox{ a.e. in $ Q $.}
\end{equation}

Next, let us put:
\begin{equation*}
    \begin{cases}
        [\bm{p}_n^*, z_n^*] = [p_n^*, p_{\Gamma, n}^*, z_n^*]:= [\bm{p}_{\varepsilon_n}^*, z_{\varepsilon_n}^*] = [p_{\varepsilon_n}^*, p_{\Gamma, \varepsilon_n}^*, z_{\varepsilon_n}^*] \mbox{ in $ \mathfrak{X} \times \mathscr{H} $,}
        \\
        A_n^* := \alpha(\eta_n^*) f_{\varepsilon_n}''(\partial_x \theta_n^*) \mbox{ in $ L^\infty(Q) $,}
    \end{cases}
    \mbox{$ n = 1, 2, 3, \dots $.}
\end{equation*}
Then, from \eqref{Thm.5-00}--\eqref{Thm.5-03}, and \eqref{productrule}, it follows that:
\begin{subequations}\label{[p,z]}
\begin{equation}\label{[p,z]-00}
    [L (u_n^* + p_n^*), L_\Gamma (u_{\Gamma, n}^* + p_{\Gamma, n}^*), M (v_n^* + z_n^*)] = [0, 0, 0] \mbox{ in $ \mathfrak{X} \times \mathscr{H} $, $ n = 1, 2, 3, \dots $,}
\end{equation}
\begin{align}\label{[p,z]-01}
    \bigl( -\partial_t & \bm{p}_n^*, \bm{\varphi} \bigr)_{\mathfrak{X}}  +\bigl( \partial_x p_n^*, \partial_x \varphi \bigr)_{\mathscr{H}} +\bigl( \mu_n^* p_n^*, \varphi \bigr)_{\mathscr{H}} +\bigl(\omega_n^* \partial_x z_n^*, \varphi \bigr)_{\mathscr{H}} 
    \nonumber
    \\
    & = \bigl(K  (\eta_n^* -\eta_\mathrm{ad}), \varphi \bigr)_{\mathscr{H}} + \bigl( K_\Gamma  (\eta_{\Gamma, n}^* -\eta_{\Gamma, \mathrm{ad}}), \varphi_\Gamma \bigr)_{\mathscr{H}_\Gamma}, \nonumber
    \\
    & \qquad \mbox{for any $ \bm{\varphi} = [\varphi, \varphi_\Gamma] \in \mathfrak{W} $, $ n = 1, 2, 3, \dots $,}
\end{align}
\begin{align}\label{[p,z]-02}
    \bigl< -\alpha_0 \partial_t & z_n^*, \psi \bigr>_{\mathscr{V}_0} +\bigl( ( -\partial_t \alpha_0) z_n^*, \psi \bigr)_{\mathscr{H}} +\bigl( A_n^* \partial_x z_n^*   +\nu^2 \partial_x z_n^* +\omega_n^* p_n^*, \partial_x \psi \bigr)_{\mathscr{H}}
    \nonumber
    \\
    & = \bigl( \Lambda (\theta_n^* -\theta_\mathrm{ad}), \psi \bigr)_{\mathscr{H}}, \mbox{ for any $ \psi \in \mathscr{V}_0 $, $ n = 1, 2, 3, \dots $,}
\end{align}
and
\begin{equation}\label{[p,z]-03}
    [\bm{p}_n^*(T), z_n^*(T)] = [p_n^*(T), p_{\Gamma, n}^*(T), z_n^*(T)] = [0, 0, 0] \mbox{ in $ \mathbb{X} \times H $, $ n = 1, 2, 3, \dots $.}
\end{equation}
\end{subequations}
Here, invoking the operators $ \mathcal{Q}_{\varepsilon}^* \in \mathscr{L}(\mathfrak{X} \times \mathscr{H}; \mathfrak{Y}) $ and $ \mathcal{R}_T \in \mathscr{L}(\mathfrak{X} \times \mathscr{H}) $ as in Remark \ref{Rem.mTh03}, we apply Theorem \ref{Prop02} to the case when:
\begin{equation*}
    \begin{array}{c}
    \begin{cases}
        [a, b, \mu, \omega, A] = \mathcal{R}_T[\alpha_0, -\partial_t \alpha_0, \mu_n^*, \omega_n^*, A_n^*],
        \\[1ex]
        [\bm{p}_0, z_0] = [p_0, p_{\Gamma, 0}, z_0] = [0, 0, 0],
        \\[1ex]
        [\bm{h}, k] = [h, h_\Gamma, k] 
        \\[1ex]
        \quad = \mathcal{R}_T\bigl[ K (\eta_n^* -\eta_\mathrm{ad}) ,  K_\Gamma (\eta_{\Gamma ,n}^* -\eta_{\Gamma, \mathrm{ad}}), \Lambda (\theta_n^* -\theta_\mathrm{ad})\bigr], 
        \\[1ex]
        [\bm{p}, z] = [p, p_\Gamma, z] 
        \\[1ex]
        \quad = \mathcal{Q}_{\varepsilon_n}^* \bigl[ \mathcal{R}_T [ K (\eta_n^* -\eta_\mathrm{ad}) ,  K_\Gamma (\eta_{\Gamma ,n}^* -\eta_{\Gamma, \mathrm{ad}}) ,  \Lambda (\theta_n^* -\theta_\mathrm{ad}) ] \bigr], 
    \end{cases}
    \mbox{for $ n \in \mathbb{N} $.}
    \end{array}
\end{equation*}
Then, with use of the constant $ \bar{C}_0^* $ as in \eqref{est3-01a}, it is deduced that:
\begin{subequations}\label{[p,z]-10}
\begin{align}\label{Gronwall01}
    \frac{d}{dt} & \bigl( \bigl| (\mathcal{R}_T \bm{p}_n^*)(t) \bigr|_\mathbb{X}^2 +\bigl| \mathcal{R}_T \bigl( \sqrt{\alpha_0} z_n^* \bigr)(t) \bigr|_H^2 \bigr)
    \nonumber
    \\
    & \qquad +\bigl( \bigl|(\mathcal{R}_T \bm{p}_n^*)(t)\bigr|_\mathbb{W}^2 +\nu^2 \bigl|(\mathcal{R}_T z_n^*)(t)\bigr|_{V_0}^2 \bigr)
    \nonumber
    \\
    & \leq \bar{C}_0^* \bigl( \bigl| (\mathcal{R}_T \bm{p}_n^*)(t) \bigr|_\mathbb{X}^2 +\bigl| \mathcal{R}_T \bigl( \sqrt{\alpha_0} z_n^* \bigr)(t) \bigr|_H^2 \bigr) 
    \nonumber
    \\
    & \qquad +\bar{C}_0^* \big( \bigl| \mathcal{R}_T \bigl( K (\eta_n^* -\eta_\mathrm{ad}) \bigr)(t) \bigr|_{V^*}^2 +  \bigl| \mathcal{R}_T \bigl( K_\Gamma (\eta_{\Gamma ,n}^* -\eta_{\Gamma, \mathrm{ad}}) \bigr)(t) \bigr|_{H_\Gamma}^2
    \nonumber
    \\
    & ~ \qquad ~ \qquad + \bigl|\mathcal{R}_T \bigl( \Lambda (\theta_n^* -\theta_\mathrm{ad}) \bigr) (t) \bigr|_{V_0^*}^2 \bigr), 
    \\
    & \mbox{for a.e. $ t \in (0, T) $, $ n = 1, 2, 3, \dots $,}
    \nonumber
\end{align}
and
\begin{align}\label{est_time}
        |\partial_t (\mathcal{R}_T \bm{p}_n^*)|_{\mathfrak{X}}^2 \leq  \bar{C}_1^* & \bigl( |[\mathcal{R}_T \bigl( K (\eta_n^* -\eta_\mathrm{ad}) \bigr), \mathcal{R}_T \bigl( K_\Gamma (\eta_{\Gamma ,n}^* -\eta_{\Gamma, \mathrm{ad}}) \bigr)]|_{\mathfrak{X}}^2 
        \nonumber
        \\
        & +|\mathcal{R}_T \bigl( \Lambda (\theta_n^* -\theta_{\mathrm{ad}}) \bigr)|_{\mathscr{V}_0^*}^2 \bigr),
\end{align}
    with $n$-independent positive constant:
    \begin{equation*}
    \bar{C}_1^*  := 4 (\bar{C}_0^*)^2 e^{\frac{3}{2}\bar{C}_0^* T}.
    \end{equation*}
\end{subequations}
As a consequence of \eqref{est3-01a}, \eqref{conv3B-02}, \eqref{[p,z]-10}, (\hyperlink{A1l}{A1}), and Gronwall's lemma, we can observe that:
\begin{itemize}
    \item[\textmd{$(\hypertarget{star2}{\star\,2})$}]the sequence $ \{ [\bm{p}_n^*, z_n^*] \}_{n = 1}^\infty = \{ [p_n^*, p_{\Gamma, n}^*, z_n^*] \}_{n = 1}^\infty $ is bounded in $ [ C([0, T]; \mathbb{X}) \times$ \linebreak $ C([0, T]; H) ] \cap [ \mathfrak{X} \times \mathscr{V}_0 ]$, and $\{\bm{p}_n^* \}_{n=1}^\infty$ is bounded in $W^{1,2}(0, T; \mathbb{X})$.
\end{itemize}

Furthermore, from \eqref{emb01}, \eqref{f_eps}, \eqref{conv3B-02}, \eqref{conv3B-06}, \eqref{[p,z]-02}, $(\hyperlink{star2}{\star\,2})$, and (\hyperlink{A1l}{A1}), we can derive the following estimate: 
\begin{align}\label{[p,z]-11}
    \bigl| \bigl< -\partial_x & (A_n^* \partial_x z_n^*), \psi \bigr>_{\mathscr{U}_{0}} \bigr| 
    = \bigl| \bigl( A_n^* \partial_x z_n^*, \partial_x \psi \bigr)_{\mathscr{H}} \bigr| 
    \nonumber
    \\
    \leq &~ \bigl| \bigl( \alpha_0 z_n^*, \partial_t \psi \bigr)_{\mathscr{H}} \bigr| +\bigl| \bigl( \nu^2 \partial_x z_n^* +\omega_n^* p_n^*, \partial_x \psi \bigr)_{\mathscr{H}} \bigr| +\bigl| \bigl( \Lambda (\theta_n^* -\theta_\mathrm{ad}), \psi \bigr)_{\mathscr{H}} \bigr|
    \\
    \leq &~ \bar{C}_2^* |\psi|_{\mathscr{U}_0}, \mbox{ for any $ \psi \in C_\mathrm{c}^\infty(Q) $, $ n = 1, 2, 3, \dots $,}
    \nonumber
\end{align}
with $ n $-independent positive constant:
\begin{equation*}
    \bar{C}_2^* := 2 \sup_{n \in \mathbb{N}} \left\{ \begin{array}{c}
        (1 +\nu^2 +|\alpha_0|_{L^\infty(Q)} +|\omega_n^*|_{L^\infty(Q)})  
        \\[1ex]
        \cdot \bigl( \bigl| [\bm{p}_n^*, z_n^*] \bigr|_{\mathfrak{W} \times \mathscr{V}_0} +\bigl| \Lambda (\theta_n^* -\theta_\mathrm{ad}) \bigr|_{\mathscr{H}} \bigr) 
    \end{array} \right\} ~(< \infty).
\end{equation*}
\medskip

Due to \eqref{conv3B-05}, \eqref{conv3B-06}, \eqref{est_time}, \eqref{[p,z]-11}, $ (\hyperlink{star2}{\star\,2}) $, and the compactness theory of Aubin's type (cf. \cite[Corollary 4]{MR0916688}), we can find subsequences of $ \{ [\bm{p}_n^*, z_n^*] \}_{n = 1}^\infty = \{ [p_n^*, p_{\Gamma, n}^*, z_n^*] \}_{n = 1}^\infty \subset \mathfrak{W} \times \mathscr{V}_0 $, $ \{ \omega_n^* \partial_x z_n^* \}_{n = 1}^\infty \subset \mathscr{H} $, and $ \{ -\partial_x (A_n^* \partial_x z_n^*) \}_{n = 1}^\infty \subset \mathscr{U}_0^* $ (not relabeled), together with the respective limits $ [\bm{p}^\circ, z^\circ] = [p^\circ, p_\Gamma^\circ, z^\circ] \in \mathfrak{W} \times \mathscr{V}_0 $ with $\bm{p}^\circ = [p^\circ, p_\Gamma^\circ] $, $ \xi^\circ \in \mathscr{H} $, and $ \zeta^\circ \in \mathscr{U}_0^* $, such that:
\begin{subequations}\label{conv3B_a}
\begin{align}\label{conv3B-10'}
    & \begin{cases}
        [\bm{p}_n^*, z_n^*] \to [\bm{p}^\circ, z^\circ] \mbox{ weakly in $ \mathfrak{W} \times \mathscr{V}_0 $,}
        \\[1ex]
        \bm{p}_n^* \to \bm{p}^\circ \mbox{ in $ \mathfrak{X} $, weakly in $ W^{1, 2}(0, T; \mathbb{V}^*) $,}
        \\
        \quad \mbox{and in the pointwise sense a.e. in $ Q $,}
    \end{cases}
\end{align}
\begin{align}\label{conv3B-11}
        \omega_n^* p_n^* \to \alpha'(\eta^\circ) \nu^\circ p^\circ \mbox{ weakly in $ \mathscr{H} $,}
\end{align}
\begin{align}\label{conv3B-12}
        \omega_n^* \partial_x z_n^* \to \xi^\circ \mbox{ weakly in $ \mathscr{H} $,}
\end{align}
and
\begin{align}\label{conv3B-13}
    -\partial_x (A_n^* \partial_x z_n^*) \to \zeta^\circ \mbox{ weakly in $ \mathscr{U}_0^* $, as $ n \to \infty $.}
\end{align}
\end{subequations}

Now, the properties \eqref{Thm.5-10}--\eqref{Thm.5-12} will be verified through the limiting observations for \eqref{[p,z]}, as $ n \to \infty $, with use of \eqref{conv3B}, \eqref{conv3B-07}, and \eqref{conv3B_a}.
\medskip

Thus, we complete the proof. \qed
}

\section{Appendix}

The objective of the appendix is to give the proofs of three Theorems \ref{Prop01}--\ref{Prop03}, that are stated as a part of auxiliary results in Section 2.

The three Theorems \ref{Prop01}--\ref{Prop03} are proved by means of the time-discretization method. In view of this, we divide the rest part in two Subsections, which are concerned with the auxiliary Lemmas in the time-discretization, and the proofs of Theorems \ref{Prop01}--\ref{Prop03}. 

\subsection{Auxiliary Lemmas in the time-discretization}
Let $ [a, b, \mu, \omega, A] \in \mathscr{S} $ be a fixed quintet of functions, and let $ \delta_*(a) $ be the positive constant as in \eqref{delta(a)}. Let $ [\bm{p}_0, z_0] = [p_0, p_{\Gamma, 0}, z_0] \in \mathbb{W} \times H $ with $ \bm{p}_0 = [p_0, p_{\Gamma, 0}]$ be a fixed initial triplet, and let $ [\bm{h}, k] = [h, h_\Gamma, k] \in \mathfrak{X} \times \mathscr{V}_0^* $ with $ \bm{h} = [h, h_\Gamma] $ be a fixed forcing triplet. 

On this basis, we denote by $ \tau \in (0, 1) $ the constant of time-step size, and consider the following time-discretization scheme for (\hyperlink{P}{P}), denoted by (\hyperlink{DP}{DP})$_\tau$. 
\begin{description}
    \item[\textmd{~~(\hypertarget{DP}{DP})$_\tau$}]Find a sequence $ \{ [\bm{p}_i, z_i] \}_{i = 1}^n = \{ [p_i, p_{\Gamma, i} z_i] \}_{i = 1}^\infty $ of triplets of functions $ [\bm{p}_i, z_i] = [p_i, p_{\Gamma, i}, z_i] \in \mathbb{W} \times V_0 $ with $ \bm{p}_i = [p_i, p_{\Gamma, i}], i = 1, 2, 3, \ldots ,$ such that:
\begin{align}\label{DP_1st}
& \parbox{12cm}{
$ \begin{array}{l}
\displaystyle 
    \frac{1}{\tau} \bigl( \bm{p}_i -\bm{p}_{i -1}, \bm{\varphi} \bigr)_{\mathbb{X}} +(\partial_x p_i, \partial_x \varphi)_H + \bigl( \mu_{i} p_i +\omega_i \partial_x z_i, \varphi \bigr)_{H} = (\bm{h}_i, \bm{\varphi})_{\mathbb{X}},
\\[2ex]
     \qquad \qquad \qquad \mbox{for every $\bm{\varphi} = [\varphi, \varphi_\Gamma] \in \mathbb{W} $, $ i = 1, 2, 3, \dots $,
    }  
\end{array} $
}
\end{align}
\vspace{-5.5ex}
\begin{align}\label{DP_2nd}
\parbox{7ex}{~~} &  
\nonumber
\\
& \parbox{12cm}{
$ \begin{array}{l}
\displaystyle \frac{1}{\tau} \bigl( a_{i} (z_i - z_{i -1}), \psi \bigr)_{H} +(b_{i} z_i, \psi )_{H} + \bigl( A_i \partial_x z_i +\nu^2 \partial_x z_i +p_i \omega_i, \partial_x \psi \bigr)_H
\\[2ex]
    \qquad \displaystyle = \langle k_i, \psi \rangle_{V_0}, \mbox{ for every $ \psi \in V_0 $, $ i = 1, 2, 3, \dots $, 
    } 
\end{array}
 $
}
\end{align}
        starting from the initial triplet $ [\bm{p}_0, z_0] = [p_0, p_{\Gamma, 0}, z_0] \in \mathbb{W} \times H $ with $ \bm{p}_0 = [p_0, p_{\Gamma, 0}] $.
\end{description}
In the context, $ \{ [a_{i}, b_{i}, \mu_i, \omega_i, A_i] \}_{i = 0}^\infty  $ is a bounded sequence in $ W^{1, \infty}(\Omega) \times  L^\infty(\Omega) \times H  \times L^\infty(\Omega) \times L^\infty(\Omega) $, such that:
\begin{subequations}\label{D-given}
    \begin{equation}\label{D-given_a}
        \begin{cases}
            \displaystyle \sup_{i \geq 0} |a_{i}|_{W^{1, \infty}(\Omega)} \leq |a|_{W^{1, \infty}(Q)}, ~ 
            \displaystyle \sup_{i \geq 0} |b_{i}|_{L^\infty(\Omega)} \leq |b|_{L^\infty(Q)}, ~ 
            \\
            \displaystyle \sup_{i \geq 0} |\mu_i|_H \leq |\mu|_{L^\infty(0, T; H)}, ~ 
            \displaystyle \sup_{i \geq 0} |\omega_i|_{L^\infty(\Omega)} \leq |\omega|_{L^\infty(Q)}, 
            \\
            \displaystyle \sup_{i \geq 0} |A_i|_{L^\infty(\Omega)} \leq |A|_{L^\infty(Q)}, ~ \sup_{i \geq 0} |\log A_i|_{L^\infty(Q)} \leq |\log A|_{L^\infty(Q)},
        \end{cases}
    \end{equation} 
    \begin{equation}\label{D-given_b}
        \mbox{$ a_{i} \geq \delta_*(a) $,  a.e. in $ \Omega $, for $ i = 0, 1, 2, \dots $,}
    \end{equation} 
    \begin{equation}\label{D-given_c}
        \left\{ \parbox{7cm}{
            $ [\overline{a}]_\tau \to a $ in $ L^\infty(0, T; C(\overline{\Omega})) $, 
            \\[1ex]
            $ [{a}]_\tau \to a $  in $ C(\overline{Q}) $, 
        } \right. 
    \end{equation}
    \begin{equation}\label{D-given_e}
        \left\{ \parbox{7cm}{
            $ [\overline{\mu}]_\tau \to \mu $ weakly-$ * $ in $ L^\infty(0, T; H) $, 
            \\[1ex]
            $ [\overline{\mu}]_\tau(t) \to \mu(t) $ in $ H $, a.e. $ t \in (0, T) $, 
        }\right. 
    \end{equation}
    \begin{align}\label{D-given_o}
        \bigl[ \partial_t [a]_\tau, \partial_x [\overline{a}]_\tau, & [\overline{b}]_\tau, [\overline{\omega}]_\tau, [\overline{A}]_\tau \bigr]
        \to [\partial_t a, \partial_x a, b, \omega, A] \mbox{ weakly-$*$ in $ [L^\infty(Q)]^5 $,}
        \nonumber
        \\ 
        & \mbox{and in the pointwise sense a.e. in $ Q $, as $ \tau \downarrow 0 $,}
    \end{align}
    and $ \{ [\bm{h}_i, k_i] \}_{i = 0}^\infty = \{ [h_i, h_{\Gamma, i}, k_i] \}_{i=0}^\infty \subset \mathbb{X} \times V_0^* $ with $ \{ \bm{h}_i \}_{i = 0}^\infty = \{ [h_i, h_{\Gamma, i}] \}_{i=0}^\infty $ is a bounded sequence, such that:
\begin{equation}\label{rxForces01}
    \begin{cases}
        \displaystyle K^* := \sup_{\tau \in (0, 1)} \bigl| \bigl[ [\overline{\bm{h}}]_\tau, [\overline{k}]_\tau \bigr] \bigr|_{\mathfrak{X} \times \mathscr{H}} < \infty,
        \\ 
        \bigl[ [\overline{\bm{h}}]_\tau, [\overline{k}]_\tau \bigl] \to [\bm{h}, k] \mbox{ in $ \mathfrak{X} \times \mathscr{H} $, as $ \tau \downarrow 0 $.}
    \end{cases}
\end{equation}
\end{subequations}
\begin{remark}\label{Rem.HA01}
    Notice that it is straightforward to obtain $ \{ [a_{i}, b_{i}, \mu_i, \omega_i, A_i] \}_{i = 0}^\infty  $ and \linebreak $ \{ [\bm{h}_i, k_i] \}_{i = 0}^\infty = \{[h_i, h_{\Gamma, i}, k_i]\}_{i = 0}^\infty $ fulfilling \eqref{D-given}, because the assumptions $ [a, b, \mu, \omega, A] \in \mathscr{S} $ and $ [\bm{h}, k] = [h, h_\Gamma, k] \in \mathfrak{X} \times \mathscr{V}_0^* $ allow us to apply the standard method as in Remark \ref{Rem.t-discrete} (\hyperlink{Fact0}{Fact\,0}). 
\end{remark}

Now, for the solvability of the time-discretization scheme (\hyperlink{DP}{DP})$_\tau$, we prepare the following lemma. 
\begin{lemma}\label{axLem01}
    Let the assume $ [a, b, \mu, \omega, A] \in \mathscr{S} $ and $ [\bm{h}, k] = [h, h_\Gamma, k] \in \mathfrak{X} \times \mathscr{V}_0^* $ with $ \bm{h} = [h, h_\Gamma] $. Let $ a^\circ \in L^\infty(\Omega) $, $ b^\circ \in L^\infty(\Omega) $, $ \mu^\circ \in H $, $ \omega^\circ \in L^\infty(\Omega) $, and $ A^\circ \in L^\infty(\Omega) $ be functions, such that
    \begin{subequations}\label{axLem01-01}
        \begin{equation} \label{axLem01-01_a}
            \begin{cases}
                \displaystyle |a^\circ|_{L^\infty(\Omega)} \leq |a|_{L^\infty(Q)}, ~  |b^\circ|_{L^\infty(\Omega)} \leq |b|_{L^\infty(Q)},
                \\
                |\mu^\circ|_H \leq |\mu|_{L^\infty(0, T; H)}, ~ |\omega^\circ|_{L^\infty(\Omega)} \leq |\omega|_{L^\infty(Q)}, 
                \\
                |A^\circ|_{L^\infty(\Omega)} \leq |A|_{L^\infty(Q)}, ~ |\log A^\circ|_{L^\infty(\Omega)} \leq |\log A|_{L^\infty(Q)},
            \end{cases}
        \end{equation}
        and
        \begin{equation} \label{axLem01-01_b}
            a^\circ \geq \delta_*(a), \mbox{ a.e. in $ \Omega $.}
        \end{equation}
    \end{subequations}
    Additionally, let us assume:
    \begin{equation}\label{axLem01-03}
        0 < \tau < \tau_0 := \frac{\min \{1, \nu^2,  \delta_*(a) \}}{16(1 +|b|_{L^\infty(Q)} +|\mu|_{L^\infty(0, T; H)}^2 +|\omega|_{L^\infty(Q)}^2)}.
    \end{equation} 
    Then, for every pairs of functions $ [\bm{h}^\circ, k^\circ] = [h^\circ, h_\Gamma^\circ, k^\circ] \in \mathbb{X} \times V_0^* $  with $ \bm{h}^\circ = [h^\circ, h_\Gamma^\circ] $ and $ [\bm{p}_0, z_0] = [p_0, p_{\Gamma, 0}, z_0] \in \mathbb{W} \times H $ with $ \bm{p}_0 = [p_0, p_{\Gamma, 0}] $, the following variational system admits a unique solution $ [\bm{p}, z] = [p, p_\Gamma, z] \in \mathbb{W} \times V_0 $ with $ \bm{p} = [p, p_\Gamma] $: 
\begin{equation}\label{ax1st}
\begin{array}{c}
    \displaystyle \frac{1}{\tau} ( \bm{p} -\bm{p}_0, \bm{\varphi} )_{\mathbb{X}}  +(\partial_x p, \partial_x \varphi)_H +\bigl( \mu^\circ p +\omega^\circ \partial_x z, \varphi \bigr)_{H}
\\[1ex]
    \displaystyle  = (\bm{h}^\circ, \bm{\varphi} )_{\mathbb{X}},
\mbox{ for any $ \bm{\varphi} = [\varphi, \varphi_\Gamma] \in \mathbb{W} $,}
\end{array}
\end{equation}
\begin{equation}\label{ax2nd}
\begin{array}{c}
\displaystyle \frac{1}{\tau} ( a^\circ(z -z_0), \psi )_H +( b^\circ z, \psi )_H +\bigl( A^\circ \partial_x z +\nu^2 \partial_x z +p \omega^\circ, \partial_x \psi \bigr)_H 
    \\[1ex]
    = \langle k^\circ, \psi \rangle_{V_0}, \mbox{ for any $ \psi \in V_0 $.}
\end{array}
\end{equation}
\end{lemma}
\paragraph{Proof.}{
    First, for the proof of existence, we  define a (non-convex) functional $ \mathcal{E} : \mathbb{X} \times H \longrightarrow (-\infty, \infty] $, by letting:
    \begin{equation*}
        \begin{array}{c}
            \mathcal{E}(\bm{p}, z) = \mathcal{E}(p, p_\Gamma, z) := \left\{ \begin{array}{ll}
            \multicolumn{2}{l}{
                \displaystyle \frac{1}{2 \tau} \bigl( |\bm{p} -\bm{p}_0|_{\mathbb{X}}^2 +|\sqrt{a^\circ}(z -z_0)|_H^2 \bigr)
            }
            \\[2ex]
            & \displaystyle +\frac{1}{2} \int_\Omega \bigl( |\partial_x p|^2 +|[A^\circ]^{\frac{1}{2}} \partial_x z|^2 +\nu^2 |\partial_x z|^2 \bigr) \, dx 
            \\[2ex]
            & \displaystyle +\frac{1}{2} \int_\Omega \bigl( \mu^\circ |p|^2 +b^\circ |z|^2 \bigr) \, dx +\int_\Omega p \bigl( \omega^\circ \cdot \partial_x z \bigr) \, dx 
            \\[2ex]
            & -(\bm{h}^\circ, \bm{p})_{\mathbb{X}} -\langle k^\circ, z \rangle_{V_0}, 
            \\[2ex]
            & \mbox{ if $  [\bm{p}, z] = [p, p_\Gamma, z] \in \mathbb{W} \times V_0 $ with $\bm{p} = [p, p_\Gamma]$,}
            \ \\
            \ \\[-0ex]
            \infty, & \mbox{otherwise,}
        \end{array} \right.
        \ \\
        \ \\[-0ex]
        \mbox{for any $ [\bm{p}, z] = [p, p_\Gamma, z] \in \mathbb{X} \times H $ with $ \bm{p} = [p, p_\Gamma] $.}
    \end{array}
\end{equation*} 
Then, by using the assumption \eqref{axLem01-03}, Remark \ref{Rem.Prelim01}, and Young's inequality, one can easily check that $ \mathcal{E} $ is a proper lower semi-continuous functional on $ \mathbb{X} \times H $, such that:
\begin{align*}
    & \mathcal{E} (\bm{p}, z)  = \mathcal{E} (p, p_\Gamma, z) \geq \frac{1}{8 \tau} \bigl( |\bm{p}|_{\mathbb{X}}^2 +\delta_*(a) |z|_H^2 \bigr) +\frac{1}{4} \bigl( |\partial_x p|_{H}^2 +\nu^2 |\partial_x z|_{H}^2 \bigr) 
\\
    & \qquad \qquad -\frac{1}{2 \tau} \bigl( |\bm{p}_0|_{\mathbb{X}}^2 +|a|_{L^\infty(Q)} |z_0|_H^2 \bigr) -\left( \frac{1}{2} |\bm{h}^\circ|_{\mathbb{X}}^2  +\frac{2}{\nu^2} |k^\circ|_{V_0^*}^2 \right),
\\
    & \qquad \qquad \qquad  \mbox{for any $ [\bm{p}, z] = [p, p_\Gamma, z] \in \mathbb{W} \times V_0 $ with $\bm{p} = [p, p_\Gamma]$,}
\end{align*}}
    via the following computations:
    \begin{subequations}\label{axLem01-02}
    \begin{align*}
        \frac{1}{2 \tau} \bigl( |\bm{p} & -\bm{p}_0|_{\mathbb{X}}^2 +|\sqrt{a^\circ} (z -z_0)|_H^2 \bigr) 
        \\
        \geq & \frac{1}{4 \tau} \bigl( |\bm{p}|_{\mathbb{X}}^2 +\delta_*(a) |z|_H^2 \bigr) -\frac{1}{2 \tau} \bigl( |\bm{p}_0|_{\mathbb{X}}^2 +|a|_{L^\infty(Q)} |z_0|_H^2 \bigr),
    \end{align*}
        \begin{align}\label{axLem01-02_a}
        \frac{1}{2} \int_\Omega & \mu^\circ |p|^2 \, dx \geq -\frac{1}{2} \bigl|\mu^\circ p \bigr|_H |p|_H  
        \geq -\frac{1}{\sqrt{2}} |\mu^\circ|_H |p|_H \bigl( |p|_H +|\partial_x p|_H \bigr)
        \nonumber
        \\
        & \geq -\frac{1}{4} |\partial_x p|_H^2 -\left( \frac{|\mu^\circ|_H}{\sqrt{2}} +\frac{|\mu^\circ|_H^2}{2} \right) |p|_H^2 
            \nonumber
            \\
            & \geq -\frac{1}{4} |\partial_x p|_H^2 -\left( |\mu|_{L^\infty(0, T; H)}^2 +\frac{1}{4} \right) |p|_H^2,
    \end{align}
    \begin{align}\label{axLem01-02_b}
        \int_\Omega p (\omega^\circ & \cdot \partial_x z) \, dx +\frac{1}{2} \int_\Omega b^\circ |z|^2 \, dx 
        \nonumber
        \\
        \geq & -\frac{\nu^2}{8} |\partial_x z|_{H}^2 -\frac{2}{\nu^2} |\omega|_{L^\infty(Q)}^2|p|_H^2 -\frac{1}{2} |b|_{L^\infty(Q)} |z|_H^2, 
    \end{align}
    and
    \begin{align}\label{axLem01-02_c}
        -(\bm{h}^\circ, \bm{p})_{\mathbb{X}} & -\langle k^\circ, z \rangle_{V_0} \geq -\frac{1}{2} |\bm{p}|_{\mathbb{X}}^2 -\frac{\nu^2}{8} |z|_{V_0}^2 -\frac{1}{2}  |\bm{h}^\circ|_{\mathbb{X}}^2 -\frac{2}{\nu^2} |k^\circ|_{V_0^*}^2.
    \end{align}
    \end{subequations}
    Additionally, when $ \tau \in (0, \tau_0) $, the system \{\eqref{ax1st},\eqref{ax2nd}\} coincides with the stationarity system for $ \min \mathcal{E} $, and hence, the solution to \{\eqref{ax1st},\eqref{ax2nd}\} is immediately obtained, by means of the direct method of calculus of variations (cf. \cite[Theorem 3.2.1]{MR2192832}).
\bigskip

    Next, to prove uniqueness, we assume that there are two solutions $ [\bm{p}^{\ell}, z^{\ell}] = [p^\ell, p_\Gamma^\ell, z^\ell]$ $ \in \mathbb{W} \times V_0 $ with $ \bm{p}^\ell = [p^\ell, p_\Gamma^\ell]  $, $ \ell = 1, 2 $, to the system \{\eqref{ax1st},\eqref{ax2nd}\}. Besides, let us take the difference between the equations \eqref{ax1st} (resp. \eqref{ax2nd}) corresponding to $ \bm{p}^\ell = [p^\ell, p_\Gamma^\ell] $ (resp. $ z^\ell $), $ \ell = 1, 2 $, and put $ \bm{\varphi} = [\varphi, \varphi_\Gamma] = [p^1 -p^2, p_\Gamma^1 -p_\Gamma^2] $ (resp. $ \psi = z^1 -z^2 $). Then, taking the sum of the results, we arrive at
\begin{align*}
    \frac{1}{\tau} & \bigl( |\bm{p}^1 -\bm{p}^2|_{\mathbb{X}}^2 +|\sqrt{a^\circ}(z^1 -z^2)|_H^2 \bigr) +  |\partial_x (p^1 -p^2)|_{H}^2
\\
&  +|[A^\circ]^{\frac{1}{2}}\partial_x (z^1 -z^2)|_{H}^2 +\nu^2 |\partial_x (z^1 -z^2)|_{H}^2 +\int_\Omega \mu^\circ |p^1 -p^2|^2 \, dx
\\
    &  +2 \int_\Omega (p^1 -p^2) \omega^\circ \cdot \partial_x (z^1 -z^2) \, dx +\int_\Omega b^\circ |z^1 -z^2|^2 \, dx= 0. 
\end{align*}
Here, applying \eqref{axLem01-02_a} and \eqref{axLem01-02_b} to the case when:
\begin{align*}
 [\bm{p}, z] = [p, p_\Gamma, z] = [\bm{p}^1 - \bm{p}^2, z^1- z^2] = [p^1 - p^2, p_\Gamma^1 - p_\Gamma^2, z^1 - z^2],
\end{align*}
and invoking \eqref{axLem01-01}, \eqref{axLem01-03}, and Young's inequality, it is inferred that:
\begin{equation*}
    \frac{1}{2 \tau} \bigl( |\bm{p}^1 -\bm{p}^2|_{\mathbb{X}}^2 +\delta_*(a)|z^1 -z^2|_H^2 \bigr) \leq 0, \mbox{ whenever $ \tau \in (0, \tau_0) $.}
\end{equation*}
Since $\delta_*(a) > 0$ (cf. Remark \ref{Rem.delta(a)}), the proof is finished. \hfill \qed
\begin{remark}\label{Rem.DPtau00}
    The existence and uniqueness of the solution to the time-discretization scheme (\hyperlink{DP}{DP})$_\tau$ are verified by applying Lemma \ref{axLem01}, inductively, for every time-steps $ i = 1, 2, 3, \dots .$ Here, we note that we can obtain the solution to the scheme (\hyperlink{DP}{DP})$_\tau$, for any sequence data of forcing $ \{ [\bm{h}_i, k_i] \}_{i = 0}^\infty = \{ [h_i, h_{\Gamma, i}, k_i] \}_{i = 0}^\infty \subset \mathbb{X} \times V_0^* $ with $ \{ \bm{h}_i \}_{i = 0}^\infty =  \{ [h_i, h_{\Gamma, i}] \}_{i = 0}^\infty $, and in particular, we do not need the assumption \eqref{rxForces01} for the solvability of (\hyperlink{DP}{DP})$_\tau$.  
\end{remark}

Next, we prepare the following Lemma, for the limiting observation of the time-discretization scheme as $  \tau \downarrow 0 $. 
\begin{lemma}\label{axLem02}
    Let $ C_0^* $ be the constant given in \eqref{C0*}. Let us assume $ [\bm{p}_0, z_0] = [p_0, p_{\Gamma, 0}, z_0] \in \mathbb{W} \times H $ with $ \bm{p}_0 = [p_0, p_{\Gamma, 0}] $, and assume that $ \{ [a_i, b_i, \mu_i, \omega_i, A_i] \}_{i = 0}^\infty \subset W^{1, \infty}(\Omega) \times L^\infty(\Omega) \times H \times L^\infty(\Omega) \times L^\infty(\Omega) $ is a given sequence satisfying \eqref{D-given_a} and \eqref{D-given_b}. Let $\{[\bm{h}_i, k_i] \}_{i=0}^\infty = \{[h_i, h_{\Gamma, i}, k_i] \}_{i=0}^\infty \subset \mathbb{X} \times V_0^*$ with $\{\bm{h}_i \}_{i=0}^\infty = \{[h_i, h_{\Gamma, i}] \}_{i=0}^\infty $ be a given sequence, and let $ \{ [\bm{p}_i, z_i] \}_{i = 1}^\infty = \{ [p_i, p_{\Gamma, i}, z_i] \}_{i = 1}^\infty$ $ \subset \mathbb{W} \times V_0 $ with $\{ \bm{p}_i \}_{i = 1}^\infty = \{ [p_i, p_{\Gamma, i}] \}_{i = 1}^\infty $ be the solution to the scheme (\hyperlink{DP}{DP})$_\tau$. Then, it holds that:
    \begin{align}\label{I}
        & \frac{1}{\tau} \bigl( |\bm{p}_i|_{\mathbb{}}^2 -|\bm{p}_{i -1}|_{\mathbb{X}}^2 \bigr) 
        +\frac{1}{\tau} \bigl( \bigl| \textstyle{\sqrt{a_i}}z_i \bigr|_H^2 -\textstyle{\sqrt{a_{i -1}}}z_{i -1} \bigr|_H^2 \bigr)
        +|\bm{p}_i|_{\mathbb{W}}^2 +\nu^2 |z_i|_{V_0}^2
        \nonumber
        \\
        & \qquad \leq \frac{C_0^*}{2} \left(\bigl( |\bm{p}_i|_{\mathbb{X}}^2 +|\bm{p}_{i -1}|_{\mathbb{X}}^2 \bigr) 
        +\bigl( \bigl| \textstyle{\sqrt{a_i}}z_i \bigr|_H^2 +\textstyle{\sqrt{a_{i -1}}}z_{i -1} \bigr|_H^2 \right)
        \\
        & \qquad \qquad \qquad +C_0^* \bigl( |\bm{h}_i|_{\mathbb{X}}^2 
        +|k_i|_{V_0^*}^2 \bigr), ~ i = 1, 2, 3, \dots,
        \nonumber
    \end{align}
    and
    \begin{align}\label{II}
        & \frac{1}{\tau} |\bm{p}_i -\bm{p}_{i -1}|_{\mathbb{X}}^2 
        +\bigl( |\partial_x p_i|_{H}^2 -|\partial_x p_{i -1}|_{H}^2 \bigr)
        \nonumber
        \\
        & \qquad \leq C_0^* \tau \bigl( |p_i|_V^2 +\nu^2 |z_i|_{V_0}^2 \bigr) +2\tau |\bm{h}_i|_{\mathbb{X}}^2, ~ i = 1, 2, 3, \dots.
    \end{align}
\end{lemma}
\paragraph{Proof.}{
    Let us fix any integer $ i \in \mathbb{N} $ of the time-step, and let us put $ \bm{\varphi} = [\varphi, \varphi_\Gamma] = \bm{p}_i = [p_i, p_{\Gamma, i}] \in \mathbb{W} $ in \eqref{DP_1st}. Then, by using Young's inequality, it is easily seen that:
    \begin{align}\label{axLem02-01}
        & \frac{1}{2 \tau} \bigl( |\bm{p}_i|_{\mathbb{X}}^2 -|\bm{p}_{i -1}|_{\mathbb{X}}^2 \bigr) 
        +|\partial_x p_i|_H^2
        \leq -\int_\Omega \mu_i |p_i|^2 \, dx -\int_\Omega \omega_i p_i \partial_x z_i \, dx +(\bm{h}_i, \bm{p}_i)_{\mathbb{X}}. 
    \end{align}
    Also, putting $ \psi = z_i \in V_0 $ in \eqref{DP_2nd}, we have:
    \begin{align}\label{axLem02-02}
        \frac{1}{2 \tau} & \bigl( |{\textstyle \sqrt{a_i}} z_i|_H^2 -|{\textstyle \sqrt{a_{i -1}}} z_{i -1}|_H^2 \bigr) +\nu^2 |\partial_x z_i|_H^2 
        \nonumber
        \\
        & \leq \frac{| a|_{W^{1, \infty}(Q)}}{2\delta_*(a)} |{\textstyle \sqrt{a_{i -1}} z_{i -1}}|_H^2 -\int_\Omega b_i |z_i|^2 \, dx -\int_\Omega \omega_i p_i \partial_x z_i \, dx +\langle k_i, z_i \rangle_{V_0},
    \end{align}
    via the computation:
    \begin{align*}
        \frac{1}{\tau} & \bigl( {\textstyle \sqrt{a_{i}} (z_i -z_{i -1}}), z_i \bigr)_H 
        \\
        & \geq \frac{1}{2 \tau} \int_\Omega \bigl( a_{i} |z_i|^2 -a_{i -1} |z_{i -1}|^2 \bigr) \, dx -\frac{1}{2} \int_\Omega \left( \frac{a_i -a_{i -1}}{\tau} \right) |z_{i -1}|^2 \, dx
        \\
        & \geq \frac{1}{2 \tau} \bigl( |{\textstyle \sqrt{a_{i}}} z_i|_H^2 -|{\textstyle \sqrt{a_{i -1}}} z_{i -1}|_H^2 \bigr) -\frac{| a|_{W^{1, \infty}(Q)}}{2\delta_*(a)} |{\textstyle \sqrt{a_{i -1}}}z_{i -1}|_H^2,
    \end{align*}
    with the use of \eqref{D-given_a}, \eqref{D-given_b}, and Young's inequality.
\bigskip 

   Now, the required inequality \eqref{I} will be verified by taking the sum of \eqref{axLem02-01} and \eqref{axLem02-02}, by invoking \eqref{C0*}, and by applying \eqref{axLem01-02} to the case when: 
    \begin{align*}
    & \begin{cases}
        [a^\circ, b^\circ, \mu^\circ, \omega^\circ, A^\circ] = [a_i, b_i, \mu_i, \omega_i, A_i], 
        \\
         [\bm{p}, z] = [p, p_\Gamma, z] = [\bm{p}_i, z_i ] = [p_i, p_{\Gamma, i}, z_i],
        \\
        [\bm{h}^\circ, k^\circ] = [h^\circ, h_\Gamma^\circ, k^\circ] = [\bm{h}_i, k_i] = [h_i, h_{\Gamma, i}, k_i].
        \end{cases}
    \end{align*}

    Next, let us put $\bm{\varphi} = [\varphi, \varphi_\Gamma] = \bm{p}_i = [p_i -p_{i -1}, p_{\Gamma, i} -p_{\Gamma, i -1}] \in \mathbb{W} $ in \eqref{DP_1st}. Then, having in mind \eqref{C0*}, and using \eqref{emb01} and Young's inequality, we will observe that:
    \begin{align*}
        & \frac{1}{\tau} |\bm{p}_i -\bm{p}_{i -1}|_{\mathbb{X}}^2 
        +\frac{1}{2} |\partial_x p_i|_H^2 -\frac{1}{2} |\partial_x p_{i -1}|_H^2
        \\
        & \leq \sqrt{2} |\mu_i|_{H} |p_i|_V |p_i -p_{i -1}|_H +|\omega_i|_{L^\infty(\Omega)} |\partial_x z_i|_{H} |p_i -p_{i -1}|_H +|\bm{h}_i|_{\mathbb{X}}|\bm{p}_i -\bm{p}_{i -1}|_{\mathbb{X}} 
        \\
        & \leq \frac{1}{2\tau} |\bm{p}_i -\bm{p}_{i -1}|_{\mathbb{X}}^2 
        +\frac{C_0^* \tau}{2} \bigl( |p_i|_V^2 +\nu^2 |z_i|_{V_0}^2 \bigr) + \tau|\bm{h}_i |_{\mathbb{X}}.
    \end{align*}
    This inequality directly leads to the required \eqref{II}. 
    \qed
    \bigskip

    Finally,  
    we prove the following Lemma concerned with a time-discrete version of Gronwall's inequality.
    \begin{lemma}\label{axLem03}
        Let $ c \geq 0 $ be a fixed constant, and let $ \tau \in (0, 1) $ be a time-step size satisfying:
        \begin{equation}\label{axLem03-00}
            0 < c \tau < 2.
        \end{equation}
        Let $ 0 < T < \infty $ be a constant of time, and let $ N_{[\frac{T}{\tau}]} \in \mathbb{N} $ be a time-step such that:
        \begin{equation}\label{nT}
            (N_{[\frac{T}{\tau}]} -1) \tau < T \leq N_{[\frac{T}{\tau}]}\tau.
        \end{equation} 
        Let $ \{P_i\}_{i = 0}^\infty \subset [0, \infty) $ and $ \{ Q_i \}_{i = 1}^{\infty} \subset [0, \infty) $  be sequences such that:
        \begin{equation}\label{axLem03-01}
            \frac{1}{\tau} (P_i -P_{i -1}) \leq \frac{c}{2} (P_i +P_{i -1}) +Q_i, ~ i = 1, 2, 3,  \dots.
        \end{equation}
        Then, it is estimated that:
        \begin{align}\label{axLem03-02}
            P_i ~& \leq 2 e^{\frac{3}{2}cT} \left( P_0 +\tau \sum_{i = 1}^{N_{[\frac{T}{\tau}]}} Q_i  \right), ~ i = 1, \dots, N_{[\frac{T}{\tau}]}.
        \end{align}
    \end{lemma}
    \paragraph{Proof.}{
        From the assumptions \eqref{axLem03-00} and \eqref{axLem03-01}, it is easily derived that:
        \begin{align*}
            P_i ~& \leq \frac{1 +\frac{c\tau}{2}}{1 -\frac{c\tau}{2}} \, P_{i -1} +\frac{\tau}{1 -\frac{c\tau}{2}} \, Q_i, ~ i = 1, 2, 3, \dots.
        \end{align*}
        On this basis, we observe that:
        \begin{align*}
            P_1 ~& \leq \frac{1 +\frac{c\tau}{2}}{1 -\frac{c\tau}{2}} P_0  +\frac{\tau}{1 -\frac{c\tau}{2}} Q_1,
            \\
            P_2 ~& \leq \left( \frac{1 +\frac{c\tau}{2}}{1 -\frac{c\tau}{2}}\right)^2 P_0 +\tau \left( \frac{1 +\frac{c\tau}{2}}{(1 -\frac{c\tau}{2})^2} Q_1 +\frac{1}{1 -\frac{c\tau}{2}} Q_2 \right),
            \\
            P_3 ~& \leq \left( \frac{1 +\frac{c\tau}{2}}{1 -\frac{c\tau}{2}}\right)^3 P_0 +\tau \left( \frac{(1 +\frac{c\tau}{2})^2}{(1 -\frac{c\tau}{2})^3} Q_1 +\frac{1 +\frac{c\tau}{2}}{(1 -\frac{c\tau}{2})^2} Q_2 +\frac{1}{1 -\frac{c\tau}{2}} Q_3 \right),
        \end{align*}
        and in general,
        \begin{align}\label{axLem03-03}
            P_i  \leq \left( \frac{1 +\frac{c\tau}{2}}{1 -\frac{c\tau}{2}}\right)^i P_0 +\tau \sum_{j = 1}^{i} \frac{(1 +\frac{c\tau}{2})^{i -j}}{(1 -\frac{c\tau}{2})^{i -j +1}} \, Q_j &\, \leq \left( \frac{1 +\frac{c\tau}{2}}{1 -\frac{c\tau}{2}}\right)^{N_{[\frac{T}{\tau}]}} \left( P_0 +\tau \sum_{i = 1}^{N_{[\frac{T}{\tau}]}}  Q_i \right), 
\nonumber
\\
\mbox{ for } i = 1, \dots, &\, N_{[\frac{T}{\tau}]}.
        \end{align}
        Here, in view of  \eqref{nT}, it is inferred that:
        \begin{align}\label{axLem03-04}
& \left( \frac{1 +\frac{c\tau}{2}}{1 -\frac{c\tau}{2}}\right)^{N_{[\frac{T}{\tau}]}}  \leq  \left( 1 +\frac{1}{~\frac{1}{c \tau} -\frac{1}{2}~} \right) \left( 1 +\frac{1}{~\frac{1}{c \tau} -\frac{1}{2}~} \right)^{N_{[\frac{T}{\tau}]} -1} 
            \nonumber
            \\[1ex]
            &  \qquad \quad \leq 2 \left( 1 +\frac{1}{~\frac{N_{[\frac{T}{\tau}]} -1}{c T} -\frac{1}{2}~} \right)^{N_{[\frac{T}{\tau}]} -1} \leq 2 \left( 1 +\frac{1}{\tilde{N}} \right)^{\tilde{N} \cdot\, cT} \left( 1 +\frac{1}{\tilde{N}} \right)^{\frac{1}{2}cT} 
\\[1ex]
& \qquad \quad \leq 2 e^{\frac{3}{2} cT}, \mbox{ with $ \tilde{N} := \frac{N_{[\frac{T}{\tau}]} -1}{cT} -\frac{1}{2} $.}
\nonumber
        \end{align}
        The estimate \eqref{axLem03-02} is obtained as a straightforward consequence of \eqref{axLem03-03} and \eqref{axLem03-04}. 
        \qed
    \subsection{Proof of Theorems \ref{Prop01}--\ref{Prop03}}

    For efficiency of explanation, we prove the three Theorems \ref{Prop01}--\ref{Prop03} in accordance with the following Steps.
    \begin{description}
        \item[\textmd{Step\,$1$:}]proof of the existence part of Theorem \ref{Prop01}.
        \item[\textmd{Step\,$2$:}]proof of Theorem \ref{Prop02} (I).
        \item[\textmd{Step\,$3$:}]proof of the uniqueness part of Theorem \ref{Prop01}.
        \item[\textmd{Step\,$4$:}]proof of Theorem \ref{Prop02} (II).
        \item[\textmd{Step\,$5$:}]proof of Theorem \ref{Prop03}.
    \end{description}
    
    \paragraph{\boldmath\underline{Step\,$1$: proof of the existence part of Theorem \ref{Prop01}.}}{
        Let $ C_0^* $ be the positive constant given in \eqref{C0*}, and let  $ \tau_0 \in (0, 1) $ be the constant given in \eqref{axLem01-03}. 
    Besides, we assume that the time-step size $ \tau \in (0, 1) $ is so small to satisfy that:
    \begin{equation*}
        0 < \tau \leq \tau_1 := \frac{1}{2 C_0^*} \left( \leq \frac{\tau_0}{2} \right),
    \end{equation*}
    and we set
    \begin{equation*}
        T_{[\tau]} := N_{[\frac{T}{\tau}]}\tau \mbox{ with use of the time-step $ N_{[\frac{T}{\tau}]} \in \mathbb{N} $ as in \eqref{nT}.}
    \end{equation*}
    On this basis, let us apply Lemma \ref{axLem03} to the inequality \eqref{I} in Lemma \ref{axLem02}, under the setting:
    \begin{align*}
        & \begin{cases}
            c = C_0^* ~ (\geq 1),  
            \\
            \displaystyle P_i = |\bm{p}_i|_{\mathbb{X}}^2 
            +|\sqrt{a_i}z_i|_H^2 +\tau \sum_{j = 0}^i \bigl( |\bm{p}_j|_{\mathbb{W}}^2 
            +\nu^2 |z_j|_{V_0}^2 \bigr) 
            \\
            \hspace{8ex} -\tau \bigl( |\bm{p}_0|_{\mathbb{W}}^2 
            +\nu^2 |z_0|_{V_0}^2 \bigr), ~ i = 0, 1, 2, 3, \dots,
            \\[2ex]
            Q_i = C_0^* \bigl( |\bm{h}_i|_{\mathbb{X}}^2 +|k_i|_{V_0^*}^2 \bigr), ~ i = 1, 2, 3, \dots.
        \end{cases}
    \end{align*}
    Then, we have:
    \begin{align*}
        |\bm{p}_i|_{\mathbb{X}}^2 & 
        +|\sqrt{a_i}z_i|_H^2 +\tau \sum_{j = 1}^i \bigl( |\bm{p}_j|_{\mathbb{W}}^2 
        +\nu^2 |z_j|_{V_0}^2 \bigr)
        \\
        & \leq 2 e^{\frac{3}{2}C_0^* T} \left( \bigl( |\bm{p}_0|_{\mathbb{X}}^2 
        +|\sqrt{a_0} z_0|_{H}^2 \bigr) +C_0^* \tau \sum_{i = 1}^{N_{[\frac{T}{\tau}]}} \bigl( |\bm{h}_i|_{\mathbb{X}}^2  
        +|k_i|_{V_0^*}^2 \bigr) \right),
        \\
        & \hspace{8ex} \mbox{for $ i = 1, \dots, N_{[\frac{T}{\tau}]} $,}
    \end{align*}
    which leads to:
    \begin{align}\label{prEx01}
    & \max \left\{ \begin{array}{c}
        \bigl| [\bm{p}]_\tau \bigr|_{C([0, T]; \mathbb{X})}^2, \bigl| [\overline{\bm{p}}]_\tau \bigr|_{L^\infty(0, T; \mathbb{X})}^2, \bigl| [\underline{\bm{p}}]_\tau \bigr|_{L^\infty(0, T; \mathbb{X})}^2,
        \\[1ex]
        \delta_*(a)\bigl| [z]_\tau \bigr|_{C([0, T]; H)}^2, \delta_*(a)\bigl| [\overline{z}]_\tau \bigr|_{L^\infty(0, T; H)}^2, \delta_*(a)\bigl| [\underline{z}]_\tau \bigr|_{L^\infty(0, T; H)}^2,
        \\[1ex]
        \bigl|[\overline{\bm{p}}]_\tau \bigr|_{L^2(0, T_{[\tau]}; \mathbb{W})}^2+ \nu^2 \bigl| [\overline{z}]_\tau \bigr|_{L^2(0, T_{[\tau]}; V_0)}^2
    \end{array} \right\}
    \nonumber
    \\
    \leq & 2 C_0^* e^{\frac{3}{2}C_0^* T} \bigl( |\bm{p}_0|_\mathbb{X}^2 +|\sqrt{a_0} z_0|_{H}^2 +\bigl|[\overline{\bm{h}}]_\tau\bigr|_{L^2(0, T_{[\tau]}; \mathbb{X})}^2 +\bigl|[\overline{k}]_\tau\bigr|_{L^2(0, T_{[\tau]}; {V}_0^*)}^2 \bigr).
    \end{align}
    Also, taking the sum of the inequality in \eqref{II}, for $ i = 1, \dots, N_{[\frac{T}{\tau}]} $, it is estimated that:
    \begin{align}\label{prEx02}
        & \bigl|\partial_t [\bm{p}]_\tau\bigr|_{L^2(0, T_{[\tau]}; \mathbb{X})}^2 +\max \bigl\{ \bigl|\partial_x [p]_\tau \bigr|_{L^\infty(0, T; H)}^2, \bigl|\partial_x [\overline{p}]_\tau \bigr|_{L^\infty(0, T; H)}^2, \bigl|\partial_x [\underline{p}]_\tau \bigr|_{L^\infty(0, T; H)}^2 \bigr\}
        \nonumber
        \\
        & \quad \leq 2|p_0|_V^2 +C_0^* \bigl( \bigl| [\overline{p}]_\tau \bigr|_{L^2(0, T_{[\tau]}; V)}^2 +\nu^2 \bigl| [\overline{z}]_\tau \bigr|_{L^2(0,  T_{[\tau]}; V_0)}^2 \bigr) +2\bigl|[\overline{\bm{h}}]_\tau\bigr|_{L^2(0, T_{[\tau]}; \mathbb{X})}^2
        \nonumber
        \\
        & \quad \leq  
        4 (C_0^*)^2 e^{\frac{3}{2}C_0^*T} \bigl(|\bm{p}_0|_{\mathbb{W}}^2 +|\sqrt{a_0} z_0|_{H}^2 +\bigl|[\overline{\bm{h}}]_\tau\bigr|_{L^2(0, T_{[\tau]}; \mathbb{X})}^2  +\bigl|[\overline{k}]_\tau\bigr|_{L^2(0, T_{[\tau]}; {V}_0^*)}^2\bigr).
    \end{align}
    Meanwhile, since \eqref{DP_2nd} implies that:
    \begin{align*}
        & \left| \frac{1}{\tau} \bigl( a_i(z_i -z_{i -1}), \psi \bigr)_H \right|
        \\
        & \leq \left| A_i \partial_x z_i +\nu^2 \partial_x z_i +p_i \omega_i \right|_H|\partial_x \psi|_H + 2|b_i|_{L^\infty(\Omega)}|z_i|_{V_0}|\psi|_{V_0} +|k_i|_{V_0^*} |\psi|_{V_0}
        \\
        & \leq \left( \left( {\frac{|A|_{L^\infty(Q)} + 2|b_i|_{L^\infty(\Omega)}}{\nu} +\nu} \right) \nu|z_i|_{V_0} +|\omega_i|_{L^\infty(\Omega)} |p_i|_H +|k_i|_{V_0^*} \right) |\psi|_{V_0}
        \\
        & \leq \frac{2}{\min \{ 1, \nu \}} \bigl( 1 +\nu + |b|_{L^\infty(Q)} +|\omega|_{L^\infty(Q)} +|A|_{L^\infty(Q)} \bigr) \bigl( |p_i|_H^2 +\nu^2 |z_i|_{V_0}^2 +|k_i|_{V_0^*}^2 \bigr)^{\frac{1}{2}} |\psi|_{V_0},
\\
        & \mbox{for all $ \psi \in V_0 $, and $ i = 1, \dots, N_{[\frac{T}{\tau}]} $,}
    \end{align*}
    one can deduce from \eqref{C0*} that:
    \begin{align}\label{prEx06}
    \bigl|[\overline{a}]_\tau(t)&\, \partial_t [z]_\tau(t) \bigr|_{V_0^*}^2 \leq C_0^*(1 + \nu + |b|_{L^\infty(Q)} + |\omega|_{L^\infty(Q)} + |A|_{L^\infty(Q)} )^2 \cdot
    \\
    & \cdot \left(\rule{-1pt}{16pt} \right. |[\overline{p}_\tau(t)]|_H^2 + \nu^2|[\overline{z}_\tau(t)]|_{V_0}^2 + |[\overline{k}]_\tau(t)|_{V_0^*}^2 \left. \rule{-1pt}{16pt} \right), \mbox{ for any $t \in [0, T_{[\tau]}]$}.
    \nonumber
    \end{align}
Additionally, integrating the both sides of \eqref{prEx06} over $[0, T_{[\tau]}]$, and invoking Remark \ref{Rem.Prelim02}, we obtain that:   
    \begin{align}\label{prEx04}
         \bigl| \partial_t &\, [z]_\tau \bigr|_{L^2(0, T_{[\tau]}; V_0^*)}^2 
        \leq \frac{(1 +\sqrt{2})^2({\textstyle|[\overline{a}]_\tau|_{L^\infty(Q)}} +|\partial_x {\textstyle[\overline{a}]_\tau}|_{L^\infty(Q)})^2}{\delta_*(a)^4} \bigl| {\textstyle[\overline{a}]_\tau} \, \partial_t [z]_\tau \bigr|_{L^2(0, T_{[\tau]}; V_0^*)}^2
        \nonumber
        \\
        & \leq (C_0^*)^5 (1 + |a|_{W^{1, \infty}(Q)})^2 \bigl( 1 +\nu + |b|_{L^\infty(Q)} +|\omega|_{L^\infty(Q)} +|A|_{L^\infty(Q)} \bigr)^2 \cdot
        \nonumber
        \\
        & \hspace{16ex} \cdot \left( \rule{-1pt}{16pt} \right. 
        \bigl| [\overline{p}]_\tau \bigr|_{L^2(0, T_{[\tau]}; H)}^2 +\nu^2 \bigl| [\overline{z}]_\tau \bigr|_{L^2(0, T_{[\tau]}; V_0)}^2 +\bigl| [\overline{k}]_\tau \bigr|_{L^2(0,  T_{[\tau]}; V_0^*)}^2 \left. \rule{-1pt}{16pt} \right)
        \nonumber
        \\
        & \leq 4(C_0^*)^6 e^{\frac{3}{2} C_0^* T} (1 + |a|_{W^{1, \infty}(Q)})^2(1 +\nu + |b|_{L^\infty(Q)} +|\omega|_{L^\infty(Q)} +|A|_{L^\infty(Q)})^2 \cdot 
        \nonumber
        \\
        & \hspace{16ex}\cdot \left( \rule{-1pt}{16pt} \right. |\bm{p}_0|_{\mathbb{X}}^2 +|\sqrt{a_0} z_0|_{H}^2 +\bigl|[\overline{\bm{h}}]_\tau\bigr|_{L^2(0, T_{[\tau]}; H)}^2  +\bigl|[\overline{k}]_\tau\bigr|_{L^2(0, T_{[\tau]}; {V}_0^*)}^2 \left. \rule{-1pt}{16pt} \right).
    \end{align}

    Now, on account of the estimates \eqref{prEx01}--\eqref{prEx04}, we can say that:
    \begin{description}
        \item[\textmd{$(\hypertarget{star3}{\star\,3})$}]$ \{ [\bm{p}]_\tau \}_{\tau \in (0, \tau_1]} = \bigl\{ \bigl[ [p]_\tau, [p_\Gamma]_\tau \bigr] \bigr\}_{\tau \in (0, \tau_1]} $ is bounded in $ W^{1, 2}(0, T; \mathbb{X}) \cap L^\infty(0, T; \mathbb{W}) $, and 
            $ \{ [\overline{\bm{p}}]_\tau \}_{\tau \in (0, \tau_1]} = \bigl\{ \bigl[ [\overline{p}]_\tau, [\overline{p_\Gamma}]_\tau \bigr] \bigr\}_{\tau \in (0, \tau_1]} $ and $ \{ [\underline{\bm{p}}]_\tau \}_{\tau \in (0, \tau_1]} = \bigl\{ \bigl[ [\underline{p}]_\tau, [\underline{p_\Gamma}]_\tau \bigr] \bigr\}_{\tau \in (0, \tau_1]} $ are bounded in $ L^\infty(0, T; \mathbb{W}) $;
        \item[\textmd{$(\hypertarget{star4}{\star\,4})$}]$ \{ [z]_\tau \}_{\tau \in (0, \tau_1]} $ is bounded in $ W^{1, 2}(0, T; V_0^*) \cap C([0, T]; H) \cap \mathscr{V}_0 $, and $ \{ [\overline{z}]_\tau \}_{\tau \in (0, \tau_1]} $ and $ \{ [\underline{z}]_\tau \}_{\tau \in (0, \tau_1]} $  are bounded in $ L^\infty(0, T; H) \cap \mathscr{V}_0 $.
    \end{description}
    In this light, we can apply the compactness theory of Aubin's type (cf. \cite[Corollary 4]{MR0916688}) with the one-dimensional compact embeddings $ V \subset C(\overline{\Omega}) $ and $V_0 \subset C(\overline{\Omega})$, and we can find a sequence $ \{ \tau_n \}_{n = 2}^\infty \subset (0, \tau_1) $, and a limiting point $ [\bm{p}, z] = [p, p_\Gamma, z] \in \mathfrak{Y} $ with $ \bm{p} = [p, p_\Gamma] $ such that:
    \begin{equation}\label{tau_n}
        \tau_1 > \tau_2 > \tau_3 > \cdots > \tau_n \downarrow 0,  \mbox{ as $ n \to \infty $,}
    \end{equation}
    \begin{subequations}\label{prEx10}
        \begin{align}\label{prEx10a}
            [\bm{p}]_{\tau_n} ~& \to \bm{p} \mbox{ in $ C(\overline{Q}) \times C(\overline{\Sigma}) $, in $ \mathfrak{W} $,}
            \nonumber
            \\
            & \mbox{weakly in $ W^{1, 2}(0, T; \mathbb{X}) $, weakly-$*$ in $ L^\infty(0, T; \mathbb{W}) $,}
            \nonumber
            \\
            & \mbox{and in the pointwise sense a.e. in $ Q $,}
        \end{align}
        \begin{align}\label{prEx10b}
            & \begin{cases}
                [\overline{\bm{p}}]_{\tau_n} \to \bm{p},
                \\
                [\underline{\bm{p}}]_{\tau_n} \to \bm{p},
            \end{cases}
            \mbox{in $ L^\infty(Q) \times L^\infty(\Sigma) $,}
            \nonumber
            \\
            & \hspace{7ex}\mbox{in $ \mathfrak{W} $, weakly-$*$ in $ L^\infty(0, T; \mathbb{W}) $,}
            \nonumber
            \\
            & \hspace{7ex}\mbox{and in the pointwise sense a.e. in $Q$, }
        \end{align}
    \end{subequations}
    \begin{subequations}\label{prEx11}
        \begin{align}\label{prEx11a}
            [z]_{\tau_n} & \to z \mbox{ in $ C([0, T]; V_0^*) $, in $ \mathscr{H} $, weakly in $ \mathscr{V}_0 $,}
            \nonumber
            \\
            & \mbox{weakly in $ W^{1, 2}(0, T; V_0^*) $, weakly-$*$ in $ L^\infty(0, T; H) $,}
            \nonumber
            \\
            & \mbox{and in the pointwise sense a.e. in $ Q $,}
        \end{align}
and
        \begin{align}\label{prEx11b}
            & [\overline{z}]_{\tau_n} \to z \mbox{ and } [\underline{z}]_{\tau_n} \to z \mbox{ in $ L^\infty([0, T]; V_0^*) $, in $ \mathscr{H} $, }
            \nonumber
            \\
            & \qquad \mbox{weakly in $ \mathscr{V}_0 $, weakly-$*$ in $ L^\infty(0, T; H) $,}
            \nonumber
            \\
            & \qquad \mbox{and in the pointwise sense a.e. in $ Q $, }
            \mbox{as $ n\to \infty $.}
        \end{align}
    \end{subequations}
    Furthermore, with \eqref{DP_1st}--\eqref{D-given}, \eqref{prEx10}, \eqref{prEx11}, and Remark \ref{Rem.Prelim02} in mind, it will be inferred that:
    \begin{align}\label{prEx12}
        \begin{cases}
            \partial_t [p]_{\tau_n} +[\overline{\mu}]_{\tau_n}[\overline{p}]_{\tau_n} +[\overline{\omega}]_{\tau_n} \partial_x [\overline{z}]_{\tau_n} -[\overline{h}]_{\tau_n}
            \\
            \qquad \to \partial_t p +\mu p +\omega \partial_x z -h \mbox{ weakly in $ \mathscr{H} $,}
            \\[1ex]
            \partial_t [p_\Gamma]_{\tau_n} -[\overline{h_\Gamma}]_{\tau_n} \to \partial_t p_\Gamma -h_\Gamma \mbox{ weakly in $ \mathscr{H}_\Gamma $,}
        \end{cases}
    \end{align}
   \begin{align}\label{prEx13}
       & \begin{cases}
           [\overline{a}]_{\tau_n} \partial_t [z]_{\tau_n} + [\overline{b}]_{\tau_n}[\overline{z}]_{\tau_n} -[\overline{k}]_{\tau_n} \to a \partial_t z + b z -k
           \\
           \qquad \mbox{ weakly in $ \mathscr{V}_0^* $,}
           \\[1ex]
           \bigl( [\overline{A}]_{\tau_n} +\nu^2 \bigr) \partial_x [\overline{z}]_{\tau_n} +[\overline{p}]_{\tau_n}[\overline{\omega}]_{\tau_n} \to \bigl( A +\nu^2) \partial_x z +p \omega 
           \\
           \qquad  \mbox{ weakly in $ \mathscr{H} $,}
       \end{cases}
       \mbox{as $ n \to \infty $,}
    \end{align}
    \begin{subequations}\label{prEx14}
        \begin{align}\label{prEx14a}
            & \int_s^t \left( \rule{-1pt}{11pt} \bigl( \partial_t [p]_{\tau_n} +[\overline{\mu}]_{\tau_n} [\overline{p}]_{\tau_n} +[\overline{\omega}]_{\tau_n} \partial_x [\overline{z}]_{\tau_n} -[\overline{h}]_{\tau_n} \bigr)(\varsigma), \varphi \right)_H \, d \varsigma
            \nonumber
            \\
            & \qquad +\int_s^t \left( \rule{-1pt}{11pt} \bigl( \partial_t [p_\Gamma]_{\tau_n} -[\overline{h_\Gamma}]_{\tau_n} \bigr)(\varsigma), \varphi_\Gamma \right)_{H_\Gamma} \, d \varsigma +\int_s^t \bigl( \partial_x [\overline{p}]_{\tau_n}(\varsigma), \partial_x \varphi \bigr)_H \, d \varsigma = 0, 
        \end{align}
        and
        \begin{align}\label{prEx14b}
            & \int_s^t \left\langle \rule{-1pt}{11pt} \bigl( [\overline{a}]_{\tau_n} \partial_t [z]_{\tau_n} + [\overline{b}]_{\tau_n}[\overline{z}]_{\tau_n} -[\overline{k}]_{\tau_n} \bigr)(\varsigma), \psi \right\rangle_{V_0} \, d \varsigma
            \nonumber
            \\
            & \qquad +\int_s^t \left( \rule{-1pt}{11pt} \bigl( ([\overline{A}]_{\tau_n} +\nu^2) \partial_x [\overline{z}]_{\tau_n} +[\overline{p}]_{\tau_n}[\overline{\omega}]_{\tau_n} \bigr)(\varsigma), \partial_x \psi \right)_{H} \, d \varsigma = 0,
        \end{align}
        \begin{center}
            for all $ \bm{\varphi} = [\varphi, \varphi_\Gamma] \in \mathbb{W} $, $ \psi \in V_0 $, $ 0 \leq s \leq t \leq T $, and $ n = 1, 2, 3, \dots $.
        \end{center}
    \end{subequations}
    On account of \eqref{prEx12}--\eqref{prEx14}, and the arbitrary choices of $ 0 \leq s \leq t \leq T $, we will verify that the limit $ [\bm{p}, z] = [p, p_\Gamma, z] \in \mathfrak{Y} $ with $ \bm{p} = [p, p_\Gamma] $ will be a solution to the system (\hyperlink{P}{P}), by letting $ n \to \infty $ in \eqref{prEx14}. 
    \qed
    }
    \paragraph{\boldmath\underline{Step\,2: proof of Theorem \ref{Prop02} (I).}}{
        Let $ [\bm{p}, z] = [p, p_\Gamma, z] \in \mathfrak{Y} $ with $ \bm{p} = [p, p_\Gamma] $ be a solution to the system (\hyperlink{P}{P}). Besides, let us put $ \bm{\varphi} = [\varphi, \varphi_\Gamma] = \bm{p}(t) = [p(t), p_{\Gamma}(t)] \in \mathbb{W} $ in \eqref{DP_1st}, put $ \psi = z(t) $ in \eqref{DP_2nd}, and take the sum of results. Then, it is seen that:
        \begin{align*}
            & \frac{1}{2} \frac{d}{dt} \bigl( \bigl| \bm{p}(t) \bigr|_{\mathbb{X}}^2 +\bigl| {\textstyle (\sqrt{a}z)(t)} \bigr|_{H}^2 \bigr) +|\partial_x p(t)|_H^2 +\nu^2 |z(t)|_{V_0}^2
            \\
            & \qquad \leq \bigl( \bm{h}(t), p(t) \bigr)_{\mathbb{X}} +\langle k(t), z(t) \rangle_{V_0} +\frac{1}{2} \int_\Omega \partial_t a(t) |z(t)|^2 \, dx -\int_\Omega \mu(t) |p(t)|^2 \, dx 
            \\
            & \qquad \qquad -\int_\Omega b(t)z(t) \, dx -2 \int_\Omega p(t) \omega(t) \partial_x z(t) \, dx, \mbox{ a.e. $ t \in (0, T) $.}
        \end{align*}
        Here, referring to the computations as in \eqref{axLem01-02}, and using the positive constant $ C_0^* $ as in \eqref{C0*}, we arrive at the conclusion \eqref{est_I-B} in Theorem \ref{Prop02} (I). 
        \qed
    }
    \paragraph{\boldmath\underline{Step\,$3$: proof of the uniqueness part of Theorem \ref{Prop01}.}}{
        For every $ \ell = 1, 2 $, let \linebreak $ [\bm{p}^\ell, z^\ell] = [p^\ell, p_\Gamma^\ell, z^\ell] \in \mathfrak{Y} $ with $ \bm{p}^{\ell} = [p^\ell, p_\Gamma^\ell] $ be the solutions to the system (\hyperlink{P}{P}) for the same initial triplet $ [\bm{p}_0, z_0] = [p_0, p_{\Gamma, 0}, z_0] \in \mathbb{W} \times H $ with $ \bm{p}_0 = [p_0, p_{\Gamma, 0}] $, and the same forcing triplet $ [\bm{h}, k] = [h, h_\Gamma, k] \in \mathfrak{X} \times \mathscr{V}_0^* $ with $ \bm{h} = [h, h_\Gamma]  $. Then, since (\hyperlink{P}{P}) is a linear system, the difference of solutions $ [\bm{p}^1 -\bm{p}^2, z^1 - z^2] = [p^1 -p^2, p_\Gamma^1 - p_\Gamma^2, z^1 - z^2] $ is also a solution to (\hyperlink{P}{P}) for the homogeneous initial triplet $ [0, 0, 0] \in \mathbb{W} \times H $ and homogeneous forcing triplet $ [0, 0, 0] \in \mathfrak{X} \times \mathscr{V}_0^* $. So, from Theorem \ref{Prop02} (I), it immediately follows that:
        \begin{align}\label{prUq01}
            & \frac{d}{dt} \bigl( \bigl| (\bm{p}^1 -\bm{p}^2)(t) \bigr|_{\mathbb{X}}^2 +\bigl| {\textstyle \sqrt{a}(z^1 -z^2)(t)} \bigr|_{H}^2 \bigr) 
            \nonumber
            \\
            & \qquad \leq C_0^* \bigl( \bigl| (\bm{p}^1 -\bm{p}^2)(t) \bigr|_{\mathbb{X}}^2 +\bigl| {\textstyle \sqrt{a}(z^1 -z^2)(t)} \bigr|_{H}^2 \bigr),
            \mbox{ a.e. $ t \in (0, T) $.}
        \end{align}
        The uniqueness result will be verified by applying Gronwall's lemma to \eqref{prUq01} with the assumption \eqref{delta(a)}. 
        \qed
    }

    \begin{remark}\label{Rem.Th01}
        By virtue of the uniqueness result in Theorem \ref{Prop01}, we can also conclude the convergence result of the time-discretization scheme (\hyperlink{DP}{DP})$_\tau$, as $ \tau \downarrow 0 $. More precisely, we can obtain the convergences as in \eqref{prEx10}--\eqref{prEx13} for         any sequence $ \{ \tau_n \}_{n = 1}^\infty $ (subsequence) satisfying \eqref{tau_n}. 
    \end{remark}

    \paragraph{\boldmath\underline{Step\,$4$: proof of Theorem \ref{Prop02} (II).}}{
        As consequences of \eqref{C1*}, \eqref{prEx01}, \eqref{prEx02}, and \eqref{prEx04}, it is inferred that:
        \begin{align}\label{prTh22_01}
            & \begin{cases}
                \bigl| \partial_t [\bm{p}]_\tau \bigr|_{\mathfrak{X}}^2 
                +\bigl| [{p}]_\tau \bigl|_{L^\infty(0, T; V)}^2 
                \\
                \qquad \leq C_1^* \bigl( |\bm{p}_0|_{\mathbb{W}}^2 +|\sqrt{a} z_0|_{H}^2 +\bigl|[\overline{\bm{h}}]_\tau\bigr|_{L^2(0, T_{[\tau]}; H)}^2  +\bigl|[\overline{k}]_\tau\bigr|_{L^2(0, T_{[\tau]}; {V}_0^*)}^2 \bigr),
                \\[2ex]
                \bigl| \partial_t [z]_\tau \bigr|_{\mathscr{V}_0^*}^2 \leq C_2^*  \bigl( |\bm{p}_0|_{\mathbb{X}}^2 +|\sqrt{a} z_0|_{H}^2 +\bigl|[\overline{\bm{h}}]_\tau\bigr|_{L^2(0, T_{[\tau]}; H)}^2  +\bigl|[\overline{k}]_\tau\bigr|_{L^2(0, T_{[\tau]}; {V}_0^*)}^2 \bigr).
            \end{cases}
        \end{align}
        Hence, having in mind \eqref{D-given}, \eqref{prEx10}, \eqref{prEx11}, and Remark \ref{Rem.Th01}, we can verify the estimate \eqref{est_I-B01} in Theorem \ref{Prop02} (II), just by letting $ \tau \downarrow 0 $ in \eqref{prTh22_01}. 
        \qed
    }
    \paragraph{\boldmath\underline{Step\,$5$: proof of Theorem \ref{Prop03}.}}{
        Since \eqref{ASY01-02} implies that:
        \begin{equation*}
            a^n \to a \mbox{ in $ C(\overline{Q}) $, as $ n \to \infty $,}
        \end{equation*}
        we may suppose:
        \begin{equation*}
            a^n \geq \frac{\delta_*(a)}{2}, \mbox{ for $ n = 1, 2, 3, \dots $,}
        \end{equation*}
        without loss of generality. Here, for any $ n \in \mathbb{N} $, we define:
        \begin{equation}\label{C0n}
            C_0^n := \frac{16 \bigl( 1 +| a^n|_{W^{1, \infty}(Q)} +|b^n|_{L^\infty(Q)} +|\mu^n|_{L^\infty(0, T; H)}^2 +|\omega^n|_{L^\infty(Q)}^2 \bigr)}{\min \{1, \nu^2, \frac{\delta_*(a)}{2}\}},
        \end{equation}
        and
        \begin{align}\label{C1n}
            & \begin{cases}
                C_1^n := 4 (C_0^n)^2 e^{\frac{3}{2}C_0^n T},
                \\[1ex]
                \displaystyle C_2^n := 4(C_0^n)^6 e^{\frac{3}{2} C_0^n T}(1 + |a^n|_{W^{1, \infty}(Q)})^2 \cdot
                \\[1ex]
                \displaystyle \qquad \cdot(1 +\nu + |b^n|_{L^\infty(Q)} +|\omega^n|_{L^\infty(Q)} +|A^n|_{L^\infty(Q)})^2. 
            \end{cases}
        \end{align}
        Then, in view of \eqref{C0*}, \eqref{C1*}, \eqref{ASY01}, \eqref{pr03-01}, \eqref{D-given_a}, \eqref{D-given_b}, \eqref{C0n}, \eqref{C1n}, and Remark \ref{Rem.C_0*}, we will infer that:
        \begin{align*}
            & \begin{cases}
                |\bm{p}^n|_{C([0, T]; \mathbb{X})}^2 +|{\textstyle \sqrt{a^n}} z^n|_{C([0, T]; H)}^2 +|\bm{p}^n|_{L^2(0, T; \mathbb{W})}^2 +\nu^2 |z^n|_{\mathbb{V}_0}^2 \leq C_3^*,
                \\[1ex]
                |\partial_t \bm{p}^n|_{\mathfrak{X}}^2 +|p^n|_{L^\infty(0, T; V)}^2 +|\partial_t z^n|_{\mathscr{V}_0^*}^2 \leq C_3^*, 
            \end{cases}
            \\[0ex]
            & \hspace{24ex}\mbox{for $ n = 1, 2, 3, \dots $,}
        \end{align*}
        with use of a uniform positive constant $ C_3^* $:
        \begin{equation*}
            C_3^* := \sup_{n \in \mathbb{N}} \left\{ \bigl( C_1^n +C_2^n \bigr) \bigl( |\bm{p}_0^n|_{\mathbb{W}}^2 +|\sqrt{a^n} z_0^n|_{H}^2 +|\bm{h}^n|_{\mathfrak{X}}^2 +|k^n|_{\mathscr{V}_0^*}^2 \bigr) \right\} < \infty. 
        \end{equation*}
    
    Now, we can say that:
    \begin{description}
        \item[\textmd{$(\hypertarget{star5}{\star\,5})$}]$ \{ \bm{p}^n \}_{n = 1}^\infty = \{ [p^n, p_\Gamma^n \bigr] \}_{n = 1}^\infty $ is bounded in $ W^{1, 2}(0, T; \mathbb{X}) \cap L^\infty(0, T; \mathbb{W}) $;
        \item[\textmd{$(\hypertarget{star6}{\star\,6})$}]$ \{ z^n \}_{n = 1}^\infty $ is bounded in $ W^{1, 2}(0, T; V_0^*) \cap C([0, T]; H) \cap \mathscr{V}_0 $.
    \end{description}
    In this light, we can apply the compactness theory of Aubin's type (cf. \cite[Corollary 4]{MR0916688}) with the one-dimensional compact embeddings $ V \subset C(\overline{\Omega}) $ and $V_0 \subset C(\overline{\Omega})$, and can find a subsequence of $ \{ [\bm{p}^n, z^n]\}_{n = 1}^\infty $ (not relabeled), and a limiting point $ [\bm{p}, z] = [p, p_\Gamma, z] \in \mathfrak{Y} $ with $ \bm{p} = [p, p_\Gamma] $ such that:
        \begin{align}\label{prCD00}
            \bm{p}^n ~& \to \bm{p} \mbox{ in $ C(\overline{Q}) \times C(\overline{\Sigma}) $, in $ \mathfrak{W} $,}
            \nonumber
            \\
            & \mbox{weakly in $ W^{1, 2}(0, T; \mathbb{X}) $, weakly-$*$ in $ L^\infty(0, T; \mathbb{W}) $,}
            \nonumber
            \\
            & \mbox{and in the pointwise sense a.e. in $ Q $,}
        \end{align}
    and
        \begin{align}\label{prCD01a}
            z^n & \to z \mbox{ in $ C([0, T]; V_0^*) $, in $ \mathscr{H} $, weakly in $ \mathscr{V}_0 $,}
            \nonumber
            \\
            & \mbox{weakly in $ W^{1, 2}(0, T; V_0^*) $, weakly-$*$ in $ L^\infty(0, T; H) $,}
            \nonumber
            \\
            & \mbox{and in the pointwise sense a.e. in $ Q $, as $ n \to \infty $.}
        \end{align}
    Therefore, with \eqref{ASY01}, \eqref{pr03-01}, \eqref{prCD00}, \eqref{prCD01a}, and Remark \ref{Rem.Prelim02} in mind, we can see that:
    \begin{align}\label{prCD02}
        \begin{cases}
            \partial_t p^n +\mu^n p^n +\omega^n \partial_x z^n -h^n
            \\
            \qquad \to \partial_t p +\mu p +\omega \partial_x z -h \mbox{ weakly in $ \mathscr{H} $,}
            \\[1ex]
            \partial_t p_\Gamma^n -h_\Gamma^n \to \partial_t p_\Gamma -h_\Gamma \mbox{ weakly in $ \mathscr{H}_\Gamma $,}
        \end{cases}
    \end{align}
   \begin{align}\label{prCD03}
       & \begin{cases}
           a^n \partial_t z^n + b^n z^n -k^n \to a \partial_t z + b z -k \mbox{ weakly in $ \mathscr{V}_0^* $,}
           \\[1ex]
           \bigl( A^n +\nu^2 \bigr) \partial_x z^n +p^n \omega^n 
           \\
           \qquad \to \bigl( A +\nu^2) \partial_x z +p \omega \mbox{ weakly in $ \mathscr{H} $,}
       \end{cases}
       \mbox{as $ n \to \infty $,}
    \end{align}
    \begin{equation}\label{prCD04}
        [\bm{p}(0), z(0)] = \lim_{n \to \infty} [\bm{p}^n(0), z^n(0)] = \lim_{n \to \infty} [\bm{p}_0^n, z_0^n] = [\bm{p}_0, z_0] \mbox{ in $ [C(\overline{\Omega}) \times H_\Gamma] \times V_0^* $,}
    \end{equation}
        \begin{subequations}\label{prCD05}
        \begin{align}\label{prCD05a}
            & \int_s^t \left( \rule{-1pt}{11pt} \bigl( \partial_t p^n +\mu^n p^n +\omega^n \partial_x z^n -h^n \bigr)(\varsigma), \varphi \right)_H \, d \varsigma
            \nonumber
            \\
            & \qquad +\int_s^t \left( \rule{-1pt}{11pt} \bigl( \partial_t p_\Gamma^n -h_\Gamma^n \bigr)(\varsigma), \varphi_\Gamma \right)_{H_\Gamma} \, d \varsigma +\int_s^t \bigl( \partial_x p^n(\varsigma), \partial_x \varphi \bigr)_H \, d \varsigma = 0, 
        \end{align}
        and
        \begin{align}\label{prCD05b}
            & \int_s^t \left\langle \rule{-1pt}{11pt} \bigl( a^n \partial_t z^n + b^n z^n -k^n \bigr)(\varsigma), \psi \right\rangle_{V_0} \, d \varsigma
            \nonumber
            \\
            & \qquad +\int_s^t \left( \rule{-1pt}{11pt} \bigl( (A^n +\nu^2) \partial_x z^n +p^n \omega^n \bigr)(\varsigma), \partial_x \psi \right)_{H} \, d \varsigma = 0,
        \end{align}
        \begin{center}
            for all $ \bm{\varphi} = [\varphi, \varphi_\Gamma] \in \mathbb{W} $, $ \psi \in V_0 $, $ 0 \leq s \leq t \leq T $, and $ n = 1, 2, 3, \dots $.
        \end{center}
    \end{subequations}
    As a consequence of \eqref{prCD02}--\eqref{prCD05}, and the arbitrary choices of $ 0 \leq s \leq t \leq T $, we will verify that the limit $ [\bm{p}, z] = [p, p_\Gamma, z] \in \mathfrak{Y}$ with $\bm{p} = [p, p_\Gamma]$ will be a solution to the system (\hyperlink{P}{P}), by letting $ n \to \infty $ in \eqref{prCD05}. Furthermore, on account of the convergences as in \eqref{prCD00} and \eqref{prCD01a}, and the uniqueness result in Theorem \ref{Prop01}, we will conclude the required convergence \eqref{pr03-02}.

    Thus, the proof of Theorem \ref{Prop03} is finished. 
    \qed    }
    }
}


\end{document}